\documentclass{article}

\def\<{\langle}
\def\>{\rangle}
\def\nl        {\hfill\break}
\def\comp      {\lower1.0ex\hbox{$\mathchar"2017$}}
\def\be{\begin}
\def\e{\end}
\def\isdef{:=}
\def\qrau{Q(R,\alpha,u)}
\def\ra{\rightarrow}
\def\lra{\longrightarrow}
\def\lla{\longleftarrow}
\def\da{\downarrow}
\def\conj{\ol{\phantom{x}}}

\def\arfred#1{\mathopen{\langle\langle}{#1}\mathclose{\rangle\rangle}}
\def\denstein#1{\mathopen{<}{#1}\mathclose{>}}
\def\denstc#1{\mathopen{<}{#1}\mathclose{>_\circ}}
\def\steinberg#1{\mathopen{\{}{#1}\mathclose{\}}}
\def\plane#1{\mathopen{(}{#1}\mathclose{)}}
\def\ol#1{\overline{#1}}
\def\te{\otimes}

\def\F         {{\cal F}}
\def\D         {{\cal D}}
\def\I         {{\cal I}}
\def\J         {{\cal J}}

\def\L         {{\cal L}}
\def\cx        {{\cal X}}
\def\P         {{\cal P}}
\def\R         {{{\bf R}\llap{\vrule height6.5ptwidth0.5pt depth0pt\kern.8em}}}
\def\C         {{{\bf C}\llap{\vrule height6pt width0.5pt depth0pt\kern.45em}}}
      
\def\Z         {{\bf Z}} 

\newcommand{\mbbf}{\mbox {F{\llap {I \kern.13em}}}\hspace{-.4em}}

\def\N {{{\bf N}\llap{\vrule height6.5pt width0.5pt depth0pt\kern.8em}}}

\def\gz{G_z}
\def\gzs{\ol{\gz}}
\def\gzm{G_{z^{-1}}}
\def\zm{z^{-1}}
\def\bs{\backslash}
\def\cee{{\cal C}}
\def\cees{\ol{\cee}}
\def\ces{\ol{C}}
\def\des{\ol{d}}
\def\hes{\ol{h}}
\def\kl{\cee\!\ell}

\newcommand{\stackler}[2]  {\mathrel{\mathop{\kern0pt #2}
                            \limits_{\scriptstyle #1}}}
\newtheorem{altel}{\hspace{-.3em}}[section]

{\begin{trivlist} \item[]{\bf Proof. }}%
{\hbox{ }\nobreak\hfill $\bullet$ \end{trivlist}}  

{\begin{trivlist} \item[]{\bf Corollary.}{#1}}{\end{trivlist}}

{\begin{trivlist} \item[]{\bf #1.}}{\end{trivlist}}

\def\surarrow{\longrightarrow\mkern-15mu\rightarrow}
\def\injarrow{\lhook\joinrel\longrightarrow}
\def\bijarrow{\lhook\joinrel\surarrow}
\def\mapright#1{\mathrel{\leavevmode\vbox
{\baselineskip0pt\lineskip.25pt\lineskiplimit0pt%
\ialign{##\cr\hfil$\scriptstyle{#1}$\hfil\cr$\longrightarrow$\cr}}}}

\def\mapleft#1{\mathrel{\leavevmode\vbox
{\baselineskip0pt\lineskip.25pt\lineskiplimit0pt%
\ialign{##\cr\hfil$\scriptstyle{#1}$\hfil\cr$\longleftarrow$\cr}}}}
\def\mapleftd#1{\mathrel{\leavevmode\vtop
{\baselineskip0pt\lineskip0pt\lineskiplimit0pt%
\ialign{##\cr$\longleftarrow$\cr\hfil$\scriptstyle{#1}$\hfil\cr}}}}
\def\mapdown#1{\Big\downarrow\rlap{$\vcenter{\hbox{$\scriptstyle#1$}}$}}
\def\mapup#1{\vcenter{\hbox{$\scriptstyle#1$}}\Big\uparrow}

\def\isodown#1{\vcenter{\hbox{$\scriptstyle\cong$}}\!\mapdown{#1}}
 
\def\double#1#2{\mathrel{\vcenter{\hbox
{\ooalign{\raise2.5pt\hbox{${#1}$}\crcr${#2}$}}}}}
\def\rightleftmaps#1#2{\double{\mapright{#1}}{\mapleftd{#2}}}

\def\swarrowpl{\raise3.5pt\hbox{$\scriptscriptstyle+$}\mkern-15mu\swarrow}
\def\diagram{\def\normalbaselines{\baselineskip20pt
                       \lineskip3pt\lineskiplimit3pt}\matrix}

\def\ssb#1{{}_{\lower 1pt\hbox{$\scriptstyle{#1}$}}}   
\def\\#1#2{{{\cal #1}\ssb{#2}}}  
\def\ce#1#2#3{{{\bf C}\llap{\vrule height#1pt width#2pt depth0pt\kern#3em}}}

\def\ker{\mathop{\rm Ker}\nolimits}
\def\im{\mathop{\rm Im}\nolimits}
\def\coker{\mathop{\rm Coker}\nolimits}
\def\tot{\mathop{\rm Tot}\nolimits}
\def\rad{\mathop{\rm rad}\nolimits}
\def\span{\mathop{\rm Span}\nolimits}
\def\tr{\mathop{\rm Tr}\nolimits}
\def\aut{\mathop{\rm Aut}\nolimits}
\def\gl{\mathop{\rm GL}\nolimits}
\def\st{\mathop{\rm St}\nolimits}
\def\gq{\mathop{\rm GQ}\nolimits}
\def\Sesq{\mathop{\rm Sesq}\nolimits}
\def\Arf{\mathop{\rm Arf}\nolimits}

\def\ent{\mathop{\rm ent\/}\nolimits}

\def\Hom{\mathop{\rm Hom}\nolimits}
\def\obj{\mathop{\rm Obj}\nolimits}

\begin{document}
\pagenumbering{roman}
\thispagestyle{empty}
\mbox{\vspace{6mm}}
\be{center}
{\large\bf GENERALIZED ARF INVARIANTS AND REDUCED\vspace{2mm}\\
POWER OPERATIONS IN CYCLIC HOMOLOGY}
\e{center}

\newpage
\mbox{  }
\thispagestyle{empty}

\newpage
\thispagestyle{empty}
\mbox{\vspace{6mm}}
\be{center}
{\large\bf GENERALIZED ARF INVARIANTS AND REDUCED\vspace{2mm}\\
POWER OPERATIONS IN CYCLIC HOMOLOGY\vspace{8mm}}
\e{center}
\be{center}
{\large\bf een wetenschappelijke proeve op het gebied van de\vspace{1mm}\\
wiskunde en \vspace{50mm}informatica}
\e{center}
\be{center}
{\large\bf PROEFSCHRIFT \vspace{10mm}}
\e{center}
\be{center}
{\bf ter verkrijging van de graad van doctor\\
aan de Katholieke Universiteit te Nijmegen,\\
volgens besluit van het college van decanen\\
in het openbaar te verdedigen\\
op donderdag 27 september 1990\\
des namiddags te 1.30 uur precies\vspace{6.5mm}\\
door \vspace{2mm}  }
\e{center}
\be{center}
{\bf PAULUS MARIA HUBERT WOLTERS \vspace{2mm}}
\e{center}
\be{center}
{\bf geboren op 1 juli 1961 te Heerlen}
\e{center}

\newpage
\noindent{\bf PROMOTOR: 
\hspace{6.32ex} PROF. DR. J.H.M. \vspace{1,5mm}STEENBRINK \newline
CO-PROMOTOR: \hspace{1,5ex} DR. F.J.-B.J. CLAUWENS\vspace{120mm}}\newline
CIP-GEGEVENS KONINKLIJKE BIBLIOTHEEK, DEN HAAG \vspace{5mm}\newline
Wolters, Paulus Maria Hubert \vspace{5mm}\newline
Generalized Arf invariants and reduced power operations in 
cyclic homology / \newline
Paulus Maria Hubert Wolters.\newline
-[S.l. : s.n.]. - Ill.\newline
Proefschrift Nijmegen. - Met lit. opg. - Met samenvatting in het Nederlands.
\newline
ISBN 90-9003548-6 \newline
SISO 513 UDC 512.66(043.3)\newline
Trefw.: homologische algebra.
\thispagestyle{empty}

\newpage
\noindent{\Large\bf Contents}\vspace{12mm}\nl
Introduction \hfill vii\linebreak
List of some frequently used notations \vspace{4mm}\hfill xi
\be{center}
{\bf Chapter I\\ Algebraic \mbox{\boldmath $K$}-theory of quadratic forms
and the Arf invariant}
\e{center}
\hspace*{3ex}\S1 Collecting the relevant $K$- and $L$-theory\hfill 1\linebreak
\hspace*{3ex}\S2 The Arf invariant   \vspace{4mm}  \hfill          16
\be{center}
{\bf Chapter II\\ New Invariants for \mbox{\boldmath $L$}-groups}
\e{center}
\hspace*{3ex}\S1 Extension of the anti-structure to the ring of formal power series \hfill 27\linebreak
\hspace*{3ex}\S2 Construction of the invariants $\omega_1^{s,h}$ and $\omega_2$ \hfill 28\linebreak
\hspace*{3ex}\S3 Recognition of $\omega_1^h$  \hfill 30\linebreak
\hspace*{3ex}\S4 Computations on the invariant $\omega_2$  \hfill 39\linebreak
\hspace*{3ex}\S5 Examples  \vspace{4mm}\hfill 51
\be{center}
{\bf Chapter III\\ Hochschild, cyclic and quaternionic homology}
\e{center}
\hspace*{3ex}\S1 Definitions and notations \hfill 59\linebreak
\hspace*{3ex}\S2 Reduced power operations \hfill 63\linebreak
\hspace*{3ex}\S3 Morita invariance \hfill 76\linebreak
\hspace*{3ex}\S4 Generalized Arf invariants \vspace{4mm}\hfill 79
\be{center}
{\bf Chapter IV\\ Applications to group rings}
\e{center}
\hspace*{3ex}\S1 Quaternionic homology of group rings \hfill 81\linebreak
\hspace*{3ex}\S2 Managing Coker$(1+\vartheta)$ \hfill 93\linebreak
\hspace*{3ex}\S3 Groups with two ends \hfill 101\linebreak
\hspace*{3ex}\S4 $\Upsilon$ for groups with two ends \vspace{5mm}\hfill 105\linebreak
References \hfill 117\linebreak
Samenvatting \hfill 119\linebreak
Curriculum Vitae \hfill 120

\newpage
\mbox{  }
\thispagestyle{empty}

\newpage

\noindent{\large{\bf Introduction}.} \nl
\vspace{3mm}\nl 
The subject of this thesis belongs to the branch of mathematics named
algebraic $K$-theory of forms, also called algebraic $L$-theory.
The main objective of the theory is to classify quadratic and hermitian 
forms over rings endowed with an 
(anti-) involution up to some notion of similarity defined on the forms.

The $L$-groups were originally designed in a geometrical context by
C.T.C. Wall in \cite{surg} as `surgery obstruction groups':
In order to classify compact manifolds up to (simple) homotopy 
equivalence, one performs `surgery', 
a certain process of cutting, replacing and glueing.
Whether this surgery succeeds depends on an obstruction
in some $L$-group of the group ring of the fundamental group
of the manifold. Conversely, every element of the $L$-group represents
a surgery problem. 
Thus the geometric question is translated into an algebraic one.
We refer to \S9 and particularly to corollary 9.4.1 of {\em loc. cit.}
for the details.

In his paper \cite{Wall;lfound} Wall defines groups $L_i^s(R)$
and $L_i^h(R)$  for arbitrary rings $R$ with anti-structure 
in a purely algebraic way.
These $L$-groups, the ones we will work with,
are in essence Grothendieck- and Whitehead groups
of categories of (non-singular) quadratic modules. We refer to chapter I
for the actual definition.
In view of the preceding it is clear that one is particularly interested
in the case where $R$ is a group ring.
Very little is known about $L$-groups of infinite groups.
The computations in the literature are either concerned with 
$L$-groups of finite groups 
(see e.g. \cite{Kolster,HM}) or, when infinite groups are considered,
the results obtained yield but very limited information (see e.g. \cite{FH}).
Our main aim is to compute $L$-groups by constructing sufficiently good
invariants.
Moreover, it is important to have such invariants at our disposal
to evaluate a particular element in an $L$-group.

To give a detailed description of our achievements in this direction,
we need the following facts.\nl
Let $R$ be a ring equipped with an anti-structure.
\begin{enumerate}
\item 
The anti-structure on $R$ induces involutions on the algebraic $K$-groups
$K_i(R)$, which in turn give rise to the
Tate cohomology groups $H^*(K_i(R))$ ($i=1,2$).
\item (\cite{Wall;lfound})
There is an invariant
$L_i^h\,(R)\lra H^{i}(K_1(R))$, called discriminant.
\item (\cite{Wall;class3})
If $I$ is a two-sided ideal of $R$ which is invariant under 
the anti-structure
and $R$ is complete with respect to the $I$-adic topology, then
$L_i^h(R)$ is isomorphic to $L_i^h(R/I)$ and there is an exact sequence 
$$\cdots\ra L_{i}^s(R)\ra L_{i}^s(R/I)\ra H^{i}(V)
\ra L_{i-1}^s(R) \ra L_{i-1}^s(R/I)\ra\cdots.$$
Here $V$ denotes the kernel of the map $K_1(R)\ra K_1(R/I)$.
\item (\cite{Giffen;k2})
There is a homomorphism
$L_i^s(R)\lra H^{i+1}(K_2(R))$, called
Hasse-Witt invariant.
\end{enumerate}
The basic idea in our construction of new invariants 
is to extend the anti-structure on the ring 
$R$ in a rather exotic way to the ring of formal power series in one
variable $R[[T]]$. 
By combining 1 to 4 in this situation, we construct 
invariants for the $L$-groups $L_i^s(R)$ and $L_i^h(R)$ 
taking values in the Tate cohomology group $H^i(K_1(R[[T]]))$
{\em viz.}
$$\omega_1^s\colon L_i^s(R)\ra H^i(K_1(R[[T]])) \hbox{ \hspace{5ex} } 
\omega_1^h\colon L_i^h(R)\ra H^i(K_1(R[[T]])).$$
Further we obtain a secondary invariant $\omega_2$ which
lives on the kernel of $\omega_1^s$ and takes values in 
some quotient of the cohomology group $H^{i+1}(K_2(R[[T]]))$.
The invariant $\omega_1^h$ 
is a generalization of the discriminant homomorphism.
To be precise, the discriminant homomorphism of 2 coincides with the 
composition of $\omega_1^h$ and the natural
projection $H^i(K_1(R[[T]]))\ra H^i(K_1(R))$.

In some of the cases where the results of
\cite{Bass-Murphy,Swan,Maazen-Stienstra,Keune,Clauwens;k}
on the computation of $K$-groups can be applied,
we are able to compute the value groups of our invariants.
For the value group of $\omega_1^h$ we do this for commutative
rings and for group rings.
We compute the value group of $\omega_2$
by assuming $R$ to carry a structure of partial $\lambda$-ring, 
by using the techniques
of \cite{Clauwens;k} to compute $K$-groups of $\lambda$-rings.
The notion of partial $\lambda$-ring occurs more or less implicit in
\cite{Joyal,Joyal;vec}. In the cases we are interested in, 
this notion is much weaker than the notion of 
$\lambda$-ring.

We confine our investigations to the `Arf-part' of $L_{2i}^*(R)$.
Roughly speaking, the `Arf-part' of $L_{2i}^*(R)$ is
the subgroup consisting of differences of similarity classes of 
quadratic forms whose underlying bilinear forms are standard forms.
We refer to section 2 of chapter I for a precise 
formulation of the definition. 
Whereas the discriminant and the Hasse-Witt invariant merely 
reveal information on the underlying bilinear form of a quadratic 
form, the new invariants $\omega_1^s$, $\omega_1^h$ and $\omega_2$ 
detect the structure of the Arf-groups.

It turns out that $\omega_1^h$ is not unfamiliar to us in the sense that
its restriction to
the Arf-group is essentially the Arf invariant
constructed by F. Clauwens in his paper \cite{Clauwens;arf}.
Clauwens proved that this Arf invariant is injective for 
group rings of finite groups. We show that $\omega_1^h$
is an isomorphism whenever $R$ is the group ring of a finite group.
The restriction of $\omega_2$ to the Arf-part of the kernel of $\omega_1^s$
yields what we call the secondary Arf invariant.
As an example we use these invariants to compute an $L$-group
of the polynomial ring in two variables over the integers.
Further we determine the Arf-group of the group presented by
$\langle X,Y,S\mid S^2=(XS)^2=(YS)^2=1,\quad XY=YX\rangle.$
This example is typical in so far as it demonstrates how a suitable 
representation of the group under consideration enables us to compute the 
Arf-group.

The fact that an instance of cyclic homology \cite{LQ}, endowed with 
some kind of squaring operation, emerges
from the computations on the value group of $\omega_2$, 
gives the impulse to search for generalized Arf invariants,
taking values in cyclic homology theory.  
We conceived the idea that some quotient of
the quaternionic homology groups J.-L. Loday
\cite{Loday} came up with, might serve as suitable
value groups for new Arf invariants.
In order to find legitimate value groups for these nascent invariants,
some effort is needed to construct appropriate operations for
the various homology theories.
The theory concerned with what we call
`reduced power operations' on (low-dimensional) Hochschild-
cyclic- and quaternionic homology groups,
is rather interesting in its own right.\nl
Ultimately we obtain an Arf invariant called $\Upsilon$ which takes
values in a quotient of the quaternionic homology group $HQ_1$.
It is a generalization of all Arf invariants mentioned so far:
when we specialize to the case of  commutative rings,
we retain 
both the ordinary and the secondary Arf invariant at the same time.
Based upon the techniques exposed in \cite{Loday} we can compute the group
$HQ_1$ for group rings. This
enables us to comprehend the value group of $\Upsilon$ to a certain extent.
This group turns out to be a rather complicated direct limit of groups,
nevertheless it is quite manageable in concrete situations and we are 
able to prove that $\Upsilon$ is injective for all 
groups that possess an infinite cyclic subgroup of finite index.
Thus we have augmented the class of groups for which we are able to
determine the Arf-group.
It is not yet clear how good this invariant really is:
on one hand there are examples of groups which do not belong to the
above-mentioned class, for which $\Upsilon$ is 
injective, while on the other hand there are examples of groups for 
which we cannot determine whether $\Upsilon$ is injective. 
\vspace{3mm}
\nl
We touch upon the contents of each chapter now.
The first section of chapter I is the reflection of an 
attempt to familiarize the reader with some of the notations, definitions 
and results from algebraic $K$- and $L$-theory, needed to define and study
the Arf-groups.
In the second section we define the Arf-groups and the 
Arf invariant. 
Then we give a presentation of the Arf-groups due to F. Clauwens 
\cite{Clauwens;arf}. Further we mention the
fact that the Arf invariant is injective for finite groups.
We conclude the first chapter with a few examples of Arf-groups of 
(infinite) groups.

In the second section of chapter II we construct the new invariants
$\omega_1^s$, $\omega_1^h$
and $\omega_2$  by extending a given anti-structure 
to the ring of formal power series $R[[T]]$ as described above.
In section 3 we compute the value group $H^0(K_1(R[[T]]))$ of 
$\omega_1^s$ and $\omega_1^h$ for commutative rings, and group rings.
Then we show that the restriction of $\omega_1^h$ to the Arf-groups coincides
with the Arf invariant in these cases.
Further we compute the cohomology group $H^1(K_2(R[[T]]))$
for rings $R$ that possess a structure of partial $\lambda$-ring, 
in section 4.
The final section of this chapter is devoted to the example of the 
polynomial ring in two variables over the integers, $\Z[X,Y]$
and the example of the group we mentioned earlier.

As we mentioned before, the surfacing of operations on
cyclic homology groups in the computations on $\omega_2$,
is the motivation for studying operations on
Hochschild, cyclic and quaternionic homology in chapter III.
We give the definitions of the various homologies
and mention the most important examples in the first section.
In the next section we construct operations on these homologies 
needed to produce well-defined Arf invariants.
Section 3 treats Morita invariance as a preparation for
the definition of the invariant $\Upsilon$ 
in the final section.

In chapter IV we apply all of the preceding to the case where $R$ is a
group ring. In the first section we compute the quaternionic homology
group $HQ_1$ of group algebras. We use this in the second section to
get an idea of what the value group of $\Upsilon$ looks like in this case.
In section 3 we characterize groups having two ends (or equivalently groups
having an infinite cyclic subgroup of finite index), in a for our purposes
convenient way, as pull-backs of finite groups and infinite cyclic
or dihedral groups.
In section 4 we show that $\Upsilon$ is injective for groups with two ends.

\nl
Most of the results in chapter II are already contained in \cite{rap}.
\newpage

\noindent{\Large\bf List of some frequently used notations}\vspace{12mm}\nl
notation \hfill page\vspace{2mm}\linebreak
$(R,\alpha,u)$ \hfill 1\linebreak
$M_n(R)$ \hfill 1\linebreak
$\P(R)$ \hfill 2\linebreak
$D_\alpha$ \hfill 2\linebreak
$\Hom_R(-,-)$ \hfill 2\linebreak
$\eta_{\alpha,u}$ \hfill 2\linebreak
$T_{\alpha,u}$ \hfill 3\linebreak
$b$ (bilinearization) \hfill 4\linebreak
$\qrau$ \hfill 4\linebreak
$E(R)$ \hfill 5\linebreak
$\gl(R)$ \hfill 5\linebreak
$K_1$ \hfill 5,6\linebreak
$\aut$ \hfill 6\linebreak
$\gq(R)$ \hfill 6\linebreak
$\Sigma_{2n}$ \hfill 6\linebreak
$t_{\alpha,u}$ \hfill 7\linebreak
$U_{2n}$ \hfill 7\linebreak
$t_{\alpha}$ (= $t$) (involution)\hfill \vspace{0.5mm}9,15\linebreak
$L_1^\cx(R,\alpha,u)$ \hfill 10\linebreak
$B(R)$ \hfill 10\linebreak
$B\qrau$ \hfill 10\linebreak
$K_0$ \hfill \vspace{0.5mm} 11\linebreak
$\widetilde{K_0}$ \hfill 11\linebreak
$\delta$ (discriminant)\hfill \vspace{0.5mm}11\linebreak
$L_0^\cx(R,\alpha,u)$\hfill 12\linebreak
${\cal N}_{k}(R)$\hfill 12\linebreak
$H^{0,1}$\hfill 14\linebreak
$\st$\hfill 14\linebreak
$e_{ij}(-)$\hfill 14\linebreak
$x_{ij}(-)$\hfill 14\linebreak
$K_2$\hfill 15\linebreak
$\Arf^\cx(R,\alpha,u)$ \hfill 16\linebreak
$\Lambda_m(R)$ \hfill 16\linebreak
$\Gamma_m(R)$ \hfill 16\linebreak
$\plane{-,-}$ \ (Arf-element) \hfill 16\linebreak
$\omega$ (Arf invariant) \hfill 20\linebreak
$\kappa(R)$ \hfill 20\linebreak
$L^{s,h}(G)$ \hfill 20\linebreak
$\Arf^{s,h}(G)$ \hfill 20\linebreak
$K(G)$ \hfill 21\linebreak
$\kl(G)$ \hfill 21\linebreak
$K_i(R,I)$ \hfill 28\linebreak
$\omega_1^s$ \hfill 28\linebreak
$\omega_1^h$ \hfill 29\linebreak
$\omega_2$ \hfill 29\linebreak
$C(R)$ \hfill 30\linebreak
$R_{{\rm ab}}$ \hfill 34\linebreak
$w_{ij}(-)$ \hfill 39\linebreak
$h_{ij}(-)$ \hfill 39\linebreak
$\denstein{-,-}$ (Dennis-Stein symbol)\hfill 39\linebreak
$\steinberg{-,-}$ (Steinberg symbol) \hfill 39\linebreak
$H_i(R)$ \hfill 61\linebreak
$HC_i(R)$ \hfill 61\linebreak
$HQ_i(R)$ \hfill 62\linebreak
$\theta_p$ \hfill 63,64\linebreak
$\nu_R$ \hfill 72\linebreak
$\mu_R$ \hfill 72\linebreak
$\vartheta_R$ \hfill 72\linebreak
$\tr$ \hfill 76\linebreak
$\Upsilon$ (generalized Arf invariant) \hfill 79\linebreak
$\gz$ \hfill \vspace{0.5mm}83\linebreak
$\gzs$ \hfill 83\linebreak
$J_\#$ \hfill 90\linebreak
$\J(G)$ \hfill 96\linebreak
$\L(c)$ \hfill 97\linebreak
$\wp(-,-)$ \hfill 106

\newpage
{\Large {\bf \be{center}
Chapter I \vspace{4mm}\\
Algebraic \mbox{\boldmath $K$}-theory of quadratic
forms\\ and the Arf invariant. 
\e{center}}}
\vspace{6mm}

\section{Collecting the relevant $K$- and $L$-theory.}
\setcounter{altel}{0}
\pagenumbering{arabic}

The material in this first section related to quadratic forms
and L-groups has primarily been extracted from \cite{Giffen;k2} and
\cite{Wall;lfound}, while the facts concerning the algebraic $K$-groups
$K_0$, $K_1$ and $K_2$ have been taken from \cite{Bass} and \cite{Milnor}.

First of all we need the notion of ring with anti-structure.
The word ring will always mean associative ring with identity, written 1.

\be{defi}
An anti-automorphism $\alpha$ of a ring $R$ is a ring isomorphism 
$\alpha\colon R\ra R^\circ$, where $R^\circ$ denotes the opposite ring of 
$R$.
A ring with anti-structure $(R,\alpha,u)$,
consists of a ring $R$,
equipped with  an anti-automorphism 
$\alpha$ of $R$
and a unit $u\in R$ such that
$\alpha(u)u=1$ and $\alpha^2(r)=uru^{-1}$ for every $r\in R$.
\e{defi}
\be{nitel}{Remark}
Let $(R,\alpha,u)$ be a ring with anti-structure.
\begin{itemize}
\item[$\cdot$]
If $u$ is central in $R$, then $\alpha$ is an anti-involution,  
i.e. an anti-automorphism of order at most 2.
\item[$\cdot$]
If $\alpha$ is the identity, 
then $R$ must be commutative and $u^2=1$.
The converse is not necessarily true.
\end{itemize}
\e{nitel}
We give a few examples of rings with anti-structure including the most
important ones.
\begin{itemize}
\item[$\cdot$]
$R$ commutative, $\alpha$ the identity, $u=\pm 1$. 
\end{itemize}
Let $(R,\alpha,u)$ be a ring with anti-structure. 
Then the anti-structure on $R$ can be extended to
\begin{itemize}
\item[$\cdot$]
the group algebra $R[G]$, for every group $G$, by the formula
$$\sum r_ig_i \mapsto \sum \alpha(r_i)g_i^{-1}.$$
\item[$\cdot$]
the ring $M_n(R)$ of $(n\times n)$-matrices over $R$ by
$$A\mapsto A^\alpha,$$
where $(A^\alpha)_{ij}\isdef \alpha(A_{ji})$.
Thus $A^\alpha$ is the conjugate transpose of $A$.
\item[$\cdot$]
the polynomial ring in one variable $R[T]$ by
$$\sum r_iT^i\mapsto \sum \alpha(r_i)(1-T)^i.$$
\item[$\cdot$]
the polynomial ring in one variable $R[T]$ by
$$\sum r_iT^i\mapsto \sum \alpha(r_i)(-T)^i.$$
\e{itemize}

\be{defi}\label{defda}
Denote by $\P(R)$ the category of finitely generated projective 
right $R$-modules and $R$-homomorphisms.
The anti-automorphism $\alpha$ enables us to define
a contravariant functor $D_\alpha\colon \P(R)\ra \P(R)$
as follows.\nl
For every $P\in \obj\,\P(R)$ we define \nl
$D_\alpha P\isdef \Hom_R(P,R)$
equipped with a right $R$-module structure by\nl
$(gr)(p)\isdef\alpha^{-1}(r)g(p)$,
for every $g\in \Hom_R(P,R)$, $p\in P$ and $r\in R$. \nl
For every $f\in \Hom_R(P,Q)$ we define\nl
$D_\alpha f\in \Hom_R(D_\alpha Q,D_\alpha P)$ by
$(D_\alpha f)(h)\isdef h\comp f$ for all $h\in D_\alpha Q $.
\e{defi}
\be{punt}\label{obsfree}
We make the following observations.
\begin{enumerate}
\item
If $M$ is free with basis $e_1,\ldots,e_m$, then $D_\alpha M$
is free with basis $e_1^*,\ldots,e_m^*$.
Here $e_i^*\in D_\alpha M$ is determined by
$e_i^*(e_j)=\delta_{ij}$ (Kronecker delta).
One calls $e_1^*,\ldots,e_m^*$ the basis dual to 
$e_1,\ldots,e_m$. 
\item
If $M$ is free with basis $e_1,\ldots,e_m$,
$N$ is free with basis $f_1,\ldots,f_n$ and
$\phi\in \Hom_R(M,N)$ has $(n\times m)$-matrix
$A$ with respect to these bases,
then $A^\alpha$ is the matrix of $D_\alpha\phi$
with respect to the dual bases.
Just as in the case of square matrices $(A^\alpha)_{ij}=\alpha(A_{ji})$
Note that $(AB)^\alpha=B^\alpha A^\alpha. $
\item
Suppose $e_1,\ldots,e_m$ and $f_1,\ldots,f_m$ are 
both bases of $M$ and $X$ is the base-change matrix.
If $A$ is the matrix of $\phi\in \Hom_R(M,D_\alpha M)$ 
with respect to $e_1,\ldots,e_m$ and its dual,
then $X^\alpha AX$ is the matrix of $\phi$ with respect to
$f_1,\ldots,f_m$ and its dual.
\end{enumerate}
\e{punt}
\be{lemma}\label{lemmaeta} \cite[section 1]{Giffen;k2}
The map $\eta_{\alpha,u}\colon 1_{\P(R)}\ra D_\alpha^2$ defined by \nl
$(\eta_{\alpha,u}P)(p)(g)\isdef u^{-1}\alpha(g(p))$ 
for every $P\in \obj\,\P(R)$, $p\in P $ and $g\in D_\alpha P $ 
is a natural equivalence.
\e{lemma}
\be{proof}
Although the proof is rather straightforward we give some
of the arguments because they might be instructive.
\begin{itemize}
\item[$\cdot$]
$(\eta_{\alpha,u}P)(p)\in D_\alpha^2 P $: 
for every $r\in R$ we have
\begin{eqnarray*}
(\eta_{\alpha,u}P)(p)(gr)&=&u^{-1}\alpha(gr(p))\\
                         &=&u^{-1}\alpha(\alpha^{-1}(r)g(p))\\
                         &=&u^{-1}\alpha(g(p))r\\
                         &=&(\eta_{\alpha,u}P)(p)(g)r
\end{eqnarray*}
\item[$\cdot$]
$\eta_{\alpha,u}P\in \Hom_R(P,D_\alpha^2 P)$: 
for every $r\in R$ we have
\begin{eqnarray*}
(\eta_{\alpha,u}P)(pr)(g)&=&u^{-1}\alpha(g(pr))\\
                         &=&u^{-1}\alpha(g(p)r)\\
                         &=&u^{-1}\alpha(r)\alpha(g(p))\\
                         &=&\alpha^{-1}(r)u^{-1}\alpha(g(p))\\
                         &=&\alpha^{-1}(r)(\eta_{\alpha,u}P)(p)(g)\\
                         &=&((\eta_{\alpha,u}P)(p)r)(g)
\end{eqnarray*}
\item[$\cdot$]
$\eta_{\alpha,u}$ is natural: for every $\phi\in \Hom_R(P,Q)$ and
$ h\in D_\alpha Q $ the diagram
$$\diagram{
P&{\buildrel \eta_{\alpha,u}P \over {\hbox to 25pt{\rightarrowfill}}}&
D_\alpha^2P\cr
\mapdown{f}&&\mapdown{D_\alpha^2f}\cr
Q&{\buildrel \eta_{\alpha,u}Q \over {\hbox to 25pt{\rightarrowfill}}}
&D_\alpha^2Q\cr}$$
commutes since
\begin{eqnarray*}
(D_\alpha^2 f)((\eta_{\alpha,u}P)(p))(h)&=&(\eta_{\alpha,u}P)(p)(D_\alpha
f(h))\\
 &=&u^{-1}\alpha(h(f(p)))\\
 &=&(\eta_{\alpha,u}Q)(f(p))(h)
\end{eqnarray*}
\item[$\cdot$]
$\eta_{\alpha,u}P$ is an isomorphism: there exists a canonical
isomorphism \nl
$ D_{\alpha}(P\oplus Q)\cong D_\alpha P\oplus D_\alpha Q $,
so we may assume that $P$ is free with basis $e_1,\ldots,e_m$ say.
From the definition of $\eta$ we deduce
$(\eta_{\alpha,u}P)(e_i)=e_i^{**}u$.
\end{itemize}
The rest is clear.
\e{proof}
\be{cor}{}
If $M$ is free with basis $e_1,\ldots,e_m$,
then $\eta_{\alpha,u}(M)\colon M\ra D_\alpha^2 M$ has matrix $uI_m$
with respect to $e_1,\ldots,e_m$ and $e_1^{**},\ldots,e_m^{**}$.
\e{cor}
\be{nota} 
From now on we write
$P^\alpha $ instead of $D_\alpha P $ and $f^\alpha $ instead
of $D_\alpha f$.
\e{nota}
\be{prop}\label{propadju}
The map $T_{\alpha,u}=T_{\alpha,u}(P,Q)\colon
 \Hom_R(Q,P^\alpha)\ra \Hom_R(P,Q^\alpha)$ 
defined by 
$$T_{\alpha,u}(f)\isdef f^\alpha\comp\eta_{\alpha,u}P$$
is a natural isomorphism and 
$T_{\alpha,u}(P,Q)\comp T_{\alpha,u}(Q,P)=1_{\Hom_R(P,Q^\alpha)}$.
In other words $T_{\alpha,u}$ defines a self-adjunction of
the functor $D_\alpha$.
\e{prop}
\be{proof}
As in \cite[Proposition 1.2]{Giffen;k2}
\e{proof}
\be{lemma}
If $M$ is free with basis $e_1,\ldots,e_m$ and
$\phi\in \Hom_R(M,M^\alpha)$ has matrix $A$ with respect to this 
basis and its dual, then $T_{\alpha,u}(\phi)$
has matrix $A^\alpha u$ with respect to the same bases.
\e{lemma}
\be{proof}
Immediate by the corollary to definition~\ref{lemmaeta} and the second
observation of~\ref{obsfree}.
\e{proof}
We are now in a position to introduce the notion of quadratic module.
\be{defi}\label{defnonsing}
In the case that $P=Q$ in proposition~\ref{propadju} we obtain a group 
endomorphism $T_{\alpha,u}\colon \Hom_R(P,P^\alpha)\lra \Hom_R(P,P^\alpha)$
satisfying $T_{\alpha,u}^2=1$.\nl
A quadratic, to be precise $(\alpha,u)$-quadratic, $R$-module
is a pair $(P,[\phi])$
consisting of a module $ P\in \obj\,\P(R) $ and the class
$[\phi]\in \coker(1-T_{\alpha,u})$ of
an element $\phi\in \Hom_R(P,P^\alpha)$.\nl
The quadratic module $(P,[\phi])$ is called non-singular
if the image $b_{[\phi]}$ of $[\phi]$ under the
`bilinearization-map' $ b\colon \coker(1-T_{\alpha,u})\ra
\ker(1-T_{\alpha,u})$, induced by the homomorphism
$1+T_{\alpha,u}\colon \Hom_R(P,P^\alpha)\ra \Hom_R(P,P^\alpha)$,
is an isomorphism.
\e{defi}
\be{nitel}{Remark}
\be{itemize}
\item[$\cdot$]
If $2$ is invertible in $R$, then $b$ is an isomorphism,
with inverse determined by 
$\phi\mapsto [\frac{1}{2}\phi].$
Thus there is a 1-1 correspondence between non-singular quadratic forms 
and symmetric non-singular bilinear forms,
i.e. elements of ${\rm Iso}(P,P^\alpha)\cap\ker(1-T_{\alpha,u})$. 
\item[$\cdot$]
In the literature one denotes by $\Sesq(P)$ the additive group
of sesquilinear forms on $P$ i.e. biadditive maps 
$\phi\colon P\times P \lra R $ satisfying 
$\phi(p_1r_1,p_2r_2)=\alpha^{-1}(r_1)\phi(p_1,p_2)r_2$
for every $p_1,p_2\in P$ and $r_1,r_2\in R$.
In the case that $R$ is commutative and $\alpha$ is the identity,
$\Sesq(P)$ is
the group of $R$-bilinear maps.
There is a bijective correspondence $\Sesq(P)\longleftrightarrow
\Hom_R(P,P^\alpha)$ by associating to an element $\phi\in \Sesq(P)$
the map $f\in \Hom_R(P,P^\alpha)$ defined by
$f(p_1)(p_2)\isdef \phi(p_1,p_2)$ for every $p_1,p_2\in P$.
\e{itemize}
\e{nitel}
We proceed to define the various categories of quadratic modules.
Along the way we shall briefly recall the relevant definitions
and facts from algebraic K-theory.\nl
The following categories and functors will all be `categories
with product' as in \cite[Ch.VII, \S1]{Bass}.

\be{defi}\label{deffunctors}
\be{itemize}
\item
Let $\qrau$ denote the category with \nl
objects: non-singular quadratic (right) $R$-modules,\nl
morphisms: $(P,[\phi])\ra(Q,[\psi])$ are the isomorphisms
$f\colon P\ra Q$ satisfying
$[f^\alpha\psi f]=[\phi]$,\nl
product: $$(P,[\phi])\perp(Q,[\psi])\isdef
(P\oplus Q, [(\pi_P)^\alpha\phi\pi_P+(\pi_Q)^\alpha\psi\pi_Q]),$$
where $\pi_P\colon P\oplus Q\ra P$ and $\pi_Q\colon P\oplus Q\ra Q$
are the natural projections.
\item
Let $\overline{\P(R)}$ denote the category with\nl
objects: objects of $\P(R)$,\nl
morphisms: isomorphisms of $\P(R)$,\nl
product: product of $\P(R)$.
\item
Now on one hand we have the forgetful functor 
$F\colon\qrau\ra\overline{\P(R)}$, which is of course product
preserving. While on the other hand there is the
so-called hyperbolic functor
$H\colon\overline{\P(R)}\ra\qrau$
defined by 
$$H(P)\isdef(P\oplus P^\alpha,[\upsilon]),
\quad H(f)\isdef f\oplus(f^\alpha)^{-1},$$ where
$\upsilon\colon P\oplus P^\alpha\ra(P\oplus P^\alpha)^\alpha$
is determined by
$(\upsilon(p,g))(p',g')\isdef g(p')$.
$H$ is product preserving as well.
The objects $H(P)$ are called hyperbolic.
\item
A product preserving functor $G\colon {\cal C}\ra {\cal D}$ is called
cofinal if for each object $A$ of ${\cal D}$ there exist objects
$B$ of ${\cal D}$ and $C$ of ${\cal C}$,
such that $A\perp B\cong G(C)$.\nl
A subcategory ${\cal C}$ of a category ${\cal D}$ is called cofinal if the
inclusion functor is cofinal.
\e{itemize}
\e{defi}
\be{lemma} \cite[theorem 3]{Wall;phil}\label{lemmahcof}.
For every $(P,[\phi])\in \obj\,\qrau$ there exists an isomorphism
$(P,[\phi])\perp (P,-[\phi])\cong H(P)$. 
Consequently $H$ is cofinal.
\e{lemma}
\be{proof}
It is not hard to verify that the morphism
$\xi\colon P\oplus P\ra P\oplus P^\alpha$
defined by
$\xi(p_1,p_2)\isdef
(p_1-b_{[\phi]}^{-1}(\phi(p_1-p_2)),b_{[\phi]}(p_1-p_2))$
does the job.
We refer to {\em loc. cit.} for a detailed proof.
\e{proof}
\be{defi}\label{defk1}
As usual $\gl(R)$ denotes the direct limit of the general linear
groups $\gl_n(R)$ consisting of invertible $n\times n$-matrices
over $R$, with respect to the embeddings
$\gl_n(R)\hookrightarrow \gl_{n+1}(R)$ defined by
\[(A)\mapsto\pmatrix{A&0\cr0&1\cr}
\mbox{ \ for all \ } (A)\in \gl_n(R).\]
A matrix is called elementary if it differs from the identity 
matrix at no more than one off-diagonal position. Denote by
$E_n(R)$ resp. $E(R)$ the subgroup of $\gl_n(R)$ resp. $\gl(R)$
generated by all elementary matrices.
According to the Whitehead lemma \cite[\S3]{Milnor}
$E(R)$ coincides with the 
commutator subgroup of $\gl(R)$.
By definition $K_1R\isdef \gl(R)/E(R)$. We use the additive notation
in the abelian group $K_1R$.
\e{defi}

There is a general procedure for defining the Whitehead group
$K_1{\cal C}$ of a category ${\cal C}$ with product,
but we do not need it for our purposes.
It follows from lemma~\ref{lemmahcof} that the $H(R^n)$ are cofinal in
$\qrau$. According to \cite[Ch.VII, \S2.3]{Bass} we may just as well
define $K_1\qrau$ as follows under these circumstances.

\be{defi} 
$K_1\qrau$ is
the commutator quotient of the direct limit
$$\lim_{\lra}\,\aut(H(R^n))$$
where the limit is taken with respect to the canonical embeddings \nl
$\aut(H(R^n))\lra\aut(H(R^n)\perp H(R))\cong\aut(H(R^{n+1})).$
\e{defi}
\be{remark}
Analogously $K_1(R)$ is the Whitehead group of both
$\P(R)$ and $\overline{\P(R)}$. Since the free modules $R^n$ are cofinal
in both categories, the groups $K_1(\P(R))$ and $K_1(\overline{\P(R)})$
both coincide with the commutator quotient of the direct limit 
\[\lim_{\lra}\,\aut(R^n)\]
where the limit is taken with respect to the canonical embeddings
$\aut(R^n)\lra \aut((R^n)\perp(R))\cong\aut(R^{n+1})$.
Upon choosing a basis for $R^n$ we may identify $\aut(R^n)$ with 
$\gl_n(R)$ and consequently
$K_1\P(R)\cong K_1(\overline{\P(R)})\cong K_1R$.
\e{remark}
\be{punt}\label{defgq}
Let us return to $\qrau$.
We choose a basis for $R^n$ and the dual basis for $(R^n)^\alpha$.
Since the matrix of $\upsilon$ with respect to these bases,
takes the form
\[\Sigma_{2n}\isdef\left(\begin{array}{cc}0&I_n\\0&0\end{array}\right)\]
we may identify $\aut(H(R^n))$ with the subgroup of $\gl_{2n}(R)$
consisting of all matrices 
\[\left(\begin{array}{cc}A&B\\C&D\end{array}\right)\in \gl_{2n}(R)
\quad\mbox{\ (here $A,B,C$ and $D$ are $n\times n$-matrices)}\] satisfying
\[\left(\begin{array}{cc}A&B\\C&D\end{array}\right)^\alpha
  \left(\begin{array}{cc}0&I_n\\0&0\end{array}\right)
  \left(\begin{array}{cc}A&B\\C&D\end{array}\right)-
  \left(\begin{array}{cc}0&I_n\\0&0\end{array}\right)=
  X-X^\alpha u\]
for some $(2n\times 2n)$-matrix $X.$
This subgroup of $\gl_{2n}(R)$ is called the general quadratic group
and is denoted by $\gq_{2n}(R)$.
As a consequence $K_1(\qrau)$ can be identified with the commutator
quotient of the group 
\[\gq(R)\isdef\lim_{\lra}\,\gq_{2n}(R),\]
where the limit is taken with respect to the embeddings
\[\gq_{2n}(R)\hookrightarrow \gq_{2(n+1)}(R)\mbox{ defined by}
\left(\begin{array}{cc}A&B\\C&D\end{array}\right)\mapsto
\left(\begin{array}{cccc}A&0&B&0\\
                         0&1&0&0\\
                         C&0&D&0\\
                         0&0&0&1
\end{array}\right).\]
\e{punt}
\be{defi}\label{defantit}
For every $n\in N$ we define
$t_{\alpha,u}\colon \gl_{2n}(R)\ra \gl_{2n}(R)$
by 
\[t_{\alpha,u}(X)=U_{2n}^{-1}X^\alpha U_{2n} 
\mbox{ for every }  X\in \gl_{2n}(R),
\mbox{ here } U_{2n}\isdef\left(\matrix{0&I_n\cr uI_n&0\cr}\right).\]
Explicitly: for every  
\[\left(\begin{array}{cc}A&B\\C&D\end{array}\right)\in \gl_{2n}(R)\]
we have
\[t_{\alpha,u}\left(\begin{array}{cc}A&B\\C&D\end{array}\right)=
\left(\begin{array}{cc}D^{\alpha^{-1}}&u^{-1}B^\alpha\\
C^\alpha u&A^\alpha\end{array}\right).\]
Note that $D^{\alpha^{-1}}=u^{-1}D^\alpha u$  
since $\alpha^2(r)=uru^{-1}$ for every $r\in R.$ 
Furthermore, $t_{\alpha,u}$ is an anti-involution since
\begin{eqnarray*}
t_{\alpha,u}^2(X)&=&U_{2n}^{-1}(U_{2n}^{-1}X^\alpha U_{2n})^\alpha U_{2n}\\
&=&U_{2n}^{-1}U_{2n}^\alpha X^{\alpha\alpha}(U_{2n}^{-1})^\alpha U_{2n}\\
&=&U_{2n}^{-2}uXu^{-1}U_{2n}^2\\
&=&X
\end{eqnarray*}
\e{defi}
\be{prop}\label{gqkar}
The following statements are equivalent:
\begin{description}
\item{(a)}
\[\left(\begin{array}{cc}A&B\\C&D\end{array}\right)\in
\gl_{2n}(R)\] belongs to $\gq_{2n}(R)$
\item{(b)}
\[\left(\begin{array}{cc}A&B\\C&D\end{array}\right)\in
\gl_{2n}(R)\] and
\[\left(\begin{array}{cc}A^\alpha C&A^\alpha D-1\\
B^\alpha C&B^\alpha D\end{array}\right)=X-X^\alpha u  \mbox{ \ for some }X\]
\item{(c)}
\[\left(\begin{array}{cc}A&B\\C&D\end{array}\right)\in
\gl_{2n}(R)\] and
\[\left\{\begin{array}{l}
A^\alpha D+C^\alpha uB=1\\
A^\alpha C+C^\alpha uA=0\\
B^\alpha D+D^\alpha uB=0\\
\mbox{the diagonal entries of $A^\alpha C$ and $B^\alpha D$ belong to}\\
\{x-\alpha(x)u\mid x\in R\}
\end{array}\right.\]
\item{(d)}
\[\left(\begin{array}{cc}A&B\\C&D\end{array}\right)\in \gl_{2n}(R)\] and
\[\left(\begin{array}{cc}A&B\\C&D\end{array}\right)^{-1}=
t_{\alpha,u}\left(\begin{array}{cc}A&B\\C&D\end{array}\right)\]
and the diagonal entries of $A^\alpha C$ and $B^\alpha D$ belong to
$\{x-\alpha(x)u\mid x\in R\}$
 \item{(e)}
\[\left\{\begin{array}{l}
A^\alpha D+C^\alpha uB=1\\
A^\alpha C+C^\alpha uA=0\\
B^\alpha D+D^\alpha uB=0\\
DA^\alpha +Cu^{-1}B^\alpha=1\\
DC^\alpha +Cu^{-1}D^\alpha=0\\
BA^\alpha +Au^{-1}B^\alpha=0\\
\mbox{the diagonal entries of $A^\alpha C$ and $B^\alpha D$ belong to}\\
\{x-\alpha(x)u\mid x\in R\}
\end{array}\right.\]
\end{description}
\end{prop}
\be{proof}\nl
{(a)}$\Leftrightarrow$ {(b)}:\nl
Immediate by writing out the condition in ~\ref{defgq}.\nl
{(b)}$\Leftrightarrow$ {(c)}:\nl 
From
\[\left(\begin{array}{cc}A^\alpha C&A^\alpha D-1\\
B^\alpha C&B^\alpha D\end{array}\right)=X-X^\alpha u  \mbox{ \ for some }X.\]
it follows that 
\[\left\{\begin{array}{l}
A^\alpha D-1=-(B^\alpha C)^\alpha u=-C^\alpha uB\\
0=A^\alpha C+(A^\alpha C)^\alpha u=A^\alpha C+C^\alpha uA\\
\mbox{the diagonal entries of $A^\alpha C$ belong to
$\{x-\alpha(x)u\mid x\in R\}$}\\
0=B^\alpha D+(B^\alpha D)^\alpha u=B^\alpha D+D^\alpha uB\\
\mbox{the diagonal entries of $B^\alpha D$ belong to
$\{x-\alpha(x)u\mid x\in R\}$}\\
\end{array}\right .\]
and vice versa.\nl
{(c)}$\Leftrightarrow$ {(d)}:\nl
The identity $A^\alpha D+C^\alpha uB=1$ holds if and only if
$D^{\alpha^{-1}}A+u^{-1}B^\alpha C=1$.
Combined with the other equations of statement {\bf (c)} this reads
\be{eqnarray*}
\left(\begin{array}{cc}1&0\\0&1\end{array}\right)&=&
\left(\begin{array}{cc}D^{\alpha^{-1}}&u^{-1}B^\alpha\\
C^\alpha u&A^\alpha\end{array}\right)
\left(\begin{array}{cc}A&B\\C&D\end{array}\right)\\
&=&\left(t_{\alpha,u}
\left(\begin{array}{cc}A&B\\C&D\end{array}\right)\right)
\left(\begin{array}{cc}A&B\\C&D\end{array}\right).
\e{eqnarray*}
The rest is obvious.\nl
{(d)}$\Leftrightarrow$ {(e)}:\nl
Immediate by writing out the equations
\[\left(\begin{array}{cc}1&0\\0&1\end{array}\right)=
\left(\begin{array}{cc}A&B\\C&D\end{array}\right)
\left(t_{\alpha,u}\left(\begin{array}{cc}A&B\\C&D\end{array}\right)\right)\]
and
\[\left(\begin{array}{cc}1&0\\0&1\end{array}\right)=
\left(t_{\alpha,u}\left(\begin{array}{cc}A&B\\C&D\end{array}\right)\right)
\left(\begin{array}{cc}A&B\\C&D\end{array}\right).\]
\e{proof}
\be{punt}
The product preserving functors 
$$F\colon\qrau\ra\overline{\P(R)},$$
$$H\colon\overline{\P(R)}\ra\qrau$$ 
and 
$$D_\alpha\colon\P(R)\ra\P(R)$$
of definition~\ref{deffunctors} and~\ref{defda} induce homomorphisms
$$F_*\colon K_1\qrau\ra K_1R,$$
$$H_*\colon K_1R\ra K_1\qrau$$
and 
$$t=t_\alpha \colon K_1R\ra K_1R.$$
Now $F_*$ is determined by 
$$F_*([X])=[X] \mbox{ for every }X\in \gq(R),$$
$H_*$ by
\[H_*([X])=\left[\left(\begin{array}{cc}X&0\\0&(X^\alpha)^{-1}\end{array}
\right)\right]\mbox{ for every } X\in \gl(R).\]
$t$ by
$$t([X])=[X^\alpha]  \mbox{ for every } X\in \gl(R).$$
Note that $t$ is an involution since 
$$t^2([X])=[X^{\alpha\alpha}]=[uXu^{-1}]=[X].$$
\e{punt}
\be{lemma}
$F_*\comp H_*=1-t$.
\e{lemma}
\be{proof}
For every $X\in \gl(R)$ we have
$$\pmatrix{X&0\cr 0&(X^\alpha)^{-1}\cr}=
\pmatrix{X(X^\alpha)^{-1}&0\cr 0&1\cr}
\pmatrix{X^\alpha&0\cr 0&(X^\alpha)^{-1}\cr}.$$
But according to \cite[\S2]{Milnor} the class of
$$\pmatrix{X^\alpha&0\cr 0&(X^\alpha)^{-1}\cr}$$
is trivial in $K_1(R)$. In view of the preceding this proves the assertion.
\e{proof}
\be{defi} \cite{Wall;lfound} 
A subgroup $\cx$ of $K_1R$ is called involution invariant if $t(\cx)=\cx$.
For every involution invariant subgroup $\cx$ of $K_1R$ define
$$L_1^\cx(R,\alpha,u)\isdef{F_*^{-1}(\cx)\over H_*(\cx)}.$$
\e{defi}
\be{defi}
\be{itemize}
\item
Let $B(R)$ denote the category with \nl
objects:
$(M,e)$ where $M$ is a free right $R$-module and 
$e=[e_1,\ldots,e_{2m}]$ is an equivalence class of bases 
of $M$;
two bases being equivalent when the base-change-matrix belongs to $E(R)$,
i.e. it represents $0\in K_1(R)$,\nl
morphisms:
isomorphisms preserving classes,\nl
product: 
$(M,e)\perp(N,f)\isdef(M\oplus N,ef)$
where \nl $e=[e_1,\ldots,e_{2m}]$,
$f=[f_1,\ldots,f_{2n}]$ and
$ef=[e_1,\ldots,e_{2m},f_1,\ldots,f_{2n}]$.
\item
Let $B\qrau$  denote the category with \nl
objects:
$(M,[\phi],e)$,\nl where
$(M,[\phi])\in \obj\,\qrau$ and
$(M,e)\in \obj\,B(R)$,\nl
morphisms:
isomorphisms preserving both structures,\nl
product: 
obvious.
\item
Again there is a product-preserving functor 
$H_b\colon B(R)\ra B\qrau$ defined by
$$H_b(M,e)\isdef(M\oplus M^\alpha,[\upsilon],ee^*) \quad
H_b(f)\isdef f\oplus(f^\alpha)^{-1}.$$
Here $e^*=[e_1^*,\ldots,e_{2m}^*]$ and $\upsilon$ is as before. 
\e{itemize}
\e{defi}
\be{lemma}\label{hbcofi}
$H_b$ is cofinal.
\e{lemma}
\be{proof}
Let $(M,\theta,e)$ be an object of $B\qrau$.
Lemma~\ref{lemmahcof} supplies a $\qrau$-isomorphism
$$\xi\colon(M,\theta,e)\perp(M,-\theta,e)\ra H_b(M,e).$$
Let $\gamma$ be the element $\xi$ determines in $K_1R$.
By choosing the class $f$ of bases of $M\oplus M$
in such a way that
$$\xi\colon (M\oplus M,\theta\perp -\theta,f)\ra H_b(M,e)$$
represents $-\gamma\in K_1R$, we obtain a $B\qrau$-isomorphism
\[\xi\perp\xi\colon
(M,\theta,e)\perp(M,-\theta,e)\perp(M\oplus M,\theta -\theta,f)
\ra H_b(M\oplus M,ee).\]
This proves the assertion. 
\e{proof}
\be{defi}
Let $({\cal C},\perp)$ be a category with product.
The Grothendieck group $K_0{\cal C}$ of ${\cal C}$
is defined as the abelian group given by the following presentation:\nl
generators: classes $[A]$ of isomorphic objects $A$ of ${\cal C}$. We
assume that these classes form a set.\nl
relations: $[A]+[B]=[A\perp B]$.
\e{defi}
\be{punt}\label{k0pres}
Lemma~\ref{hbcofi} implies that
\begin{itemize} 
\item[$\cdot$]
each element of $K_0B\qrau$ can be written in the form
$[A]-[B]$ where $A\in B\qrau$ and $B$ is hyperbolic.
\item[$\cdot$]
the equality $[A]-[B]=[A']-[B']$ holds in $K_0B\qrau$
if and only if there exists a hyperbolic object
$C$ such that $A\perp B'\perp C\cong A'\perp B\perp C$.
\end{itemize}
\end{punt}
\be{defi}
Define $\widetilde{K_0}B\qrau$ as the kernel of the rank-map
$$rk\colon\widetilde{K_0}B\qrau\ra\Z$$ induced by the map 
$$B\qrau\ra\Z \mbox{ \ given by \ } 
(M,\theta,[e_1,\ldots,e_{2m}])\mapsto 2m.$$
\e{defi}
\be{defi}
The map $B\qrau\ra K_1R$ determined by
$$(M,\theta,e)\mapsto[\hbox{a `matrix' of }
b_\theta\hbox{ with respect to }e\hbox{ and }e^*]$$
induces a homomorphism 
$\delta\colon K_0B\qrau\ra K_1R$, called discriminant.
\e{defi}
\be{remark}
$b_\theta $ determines a matrix with respect to $e$ and $e^*$
only up to elementary matrices. It is therefore legitimate to speak about
the class of this `matrix' in $K_1R$.
\e{remark}
\be{remark}\label{hypnontriv}
Further we ought to mention the fact that $\delta$ is a priori non-trivial
on hyperbolic objects:\nl
given a hyperbolic object $H_b(M,e)=(M\oplus M^\alpha,[\upsilon],ee^*)$ in
$B\qrau$, the matrix $\Sigma_{2m}$ of $\upsilon$ actually 
(not only up to elementary matrices) takes the form
\[\Sigma_{2m}=\left(\begin{array}{cc}0&I_m\\0&0\end{array}\right)
\mbox{ (no matter what $e$ looks like).}\]
Hence $b_{[\upsilon]}$ has matrix
$U_{2m}=\left(\matrix{0&I_m\cr uI_m&0\cr}\right)$.
The class of this matrix in $K_1R$ is not necessarily trivial.
\e{remark}
\be{defi} \cite[\S3]{Wall;lfound}\label{deftau}
Define a homomorphism
$\tau\colon K_1(R)\ra \widetilde{K_0}B\qrau$ as follows :\nl
Suppose we are given an $x\in K_1(R)$.
Choose $(M,\theta,e)\in B\qrau$ and
$\gamma\in\aut(M)$ in such a way that the matrix determined by $\gamma$
represents $x$ in $K_1(R)$.
Define $\tau([x])\isdef[(M,\theta,\gamma(e)]-[M,\theta,e]$ where
$\gamma(e)=[\gamma(e_1),\ldots,\gamma(e_{2m})]$.
It is not hard to check that $\tau$ is a  well-defined homomorphism.
\e{defi}
\be{lemma}
$\delta\comp\tau=1+t$.
\e{lemma}
\be{proof}
Using the third observation of ~\ref{obsfree} we obtain
$$\delta\comp\tau([A])=[A^\alpha BA]-[B]=[A^\alpha A]=(1+t)([A])\quad
\mbox {for all } \quad A\in\gl(R),$$
where $B$ is a `matrix' of 
$b_\theta$ and $\theta$ is as in the construction of $\tau$.
\e{proof}
\be{defi}
For every involution invariant subgroup
$\cx$ of $K_1R$ define
$$L_0^\cx(R,\alpha,u)\isdef{\delta^{-1}(\cx)\over\tau(\cx)}$$
here  $\delta\colon\widetilde{K_0}B\qrau\ra K_1(R)$ is the 
restriction of the 
discriminant.
\e{defi}
\be{nota}
Write $L_\varepsilon^s$ instead of $L_\varepsilon^{\{0\}}$
and $L_\varepsilon^h$ instead of $L_\varepsilon^{K_1(R)}$
for $\varepsilon=0,1$.
\e{nota}
\be{punt}\label{obslgroup}
Let $(R,\alpha,u)$ be a ring with anti-structure and $\cx$ an 
involution invariant subgroup of $K_1(R).$
Every element $l$ of $L_0^\cx(R,\alpha,u)$ 
can be written in the form 
$$[M,[\phi],e]-[M',[\phi'],e'],$$ 
with $rk([M,[\phi],e])=rk([M',[\phi'],e'])=2m$ say.
Let $$\Gamma([M,[\phi],e])\; \hbox{ resp. }\; \Gamma([M',[\phi'],e'])$$ 
denote the matrix of $\phi$ resp. $\phi'$ with respect 
to a basis in the class  $e$ resp. $e'$.
Since  the quadratic modules 
$(M,[\phi])$ and $(M',[\phi'])$ are non-singular,
it follows from definition~\ref{defnonsing} that these matrices 
belong to 
${\cal N}_{2m}(R),$
where $${\cal N}_{k}(R)
\isdef\{\Gamma\in M_k(R)\mid \Gamma+\Gamma^\alpha u\in \gl_k(R)\}.$$
We associate to $l$ the difference  
$$[\Gamma([M,[\phi],e])]-[\Gamma([M',[\phi'],e'])]$$
of classes with respect to the following relations:
\be{itemize}
\item[$\diamond$]
For all  
$\Gamma_1,\Gamma_1'\in{\cal N}_{2m_1}(R)$ and  
$\Gamma_2,\Gamma_2'\in{\cal N}_{2m_2}(R),$ 
$$[\Gamma_1]-[\Gamma_1']+[\Gamma_2]-[\Gamma_2']= 
[\Gamma_1\perp\Gamma_2]-[\Gamma_1'\perp\Gamma_2'].$$
where $\perp$ is determined by
$$\pmatrix{A&B\cr C&D\cr}\perp\pmatrix{A'&B'\cr C'&D'\cr}=
\pmatrix{A&0&B&0\cr0&A'&0&B'\cr C&0&D&0\cr 0&C'&0&D'\cr}.$$
This follows from the definition of the product in $B\qrau$.
\item[$\diamond$]
For all $\Xi\in M_{2m}(R)$ 
$$[\Gamma]=[\Gamma+\Xi-\Xi^\alpha u]. $$
This is clear in view of definition~\ref{defnonsing} and the 
observations of~\ref{obsfree}.
\item[$\diamond$]
For all $\Delta\in \gl_{2m}(R)$ with $[\Delta]\in\cx$ 
$$[\Gamma]=[\Delta^\alpha\Gamma\Delta].$$ 
This is a consequence of definition~\ref{deffunctors} and the 
observations of~\ref{obsfree}. 
\e{itemize}
Conversely, for all $\Gamma,\Gamma'\in{\cal N}_{2m}(R)$
we associate to $[\Gamma]-[\Gamma']$ the element 
$$[R^{2m},[\phi],st]-[R^{2m},[\phi'],st]\in L_0^\cx(R,\alpha,u).$$ 
Here $st$ denotes the standard basis of $R^{2m}$ and 
$\phi$ resp. $\phi'$
is the homomorphism which has matrix $\Gamma$ resp. $\Gamma'$ with 
respect to this standard basis.

Thus we have established a bijective correspondence between 
elements of $L_0^\cx(R,\alpha,u)$ and
differences of classes of elements of 
${\cal N}_{2m}(R)$ under the given relations.
Regarding the first item of \ref{k0pres} 
we may thus write every element of 
$L_0^\cx(R,\alpha,u)$ as a difference $[\Gamma]-[\Sigma_{2m}]$, with 
$\Gamma\in{\cal N}_{2m}(R)$.

Finally, we interpret the second item of \ref{k0pres} as follows.
For all $\Gamma\in{\cal N}_{2m}(R)$ and
$\Gamma'\in{\cal N}_{2m'}(R)$,
$$[\Gamma]-[\Sigma_{2m}]=[\Gamma']-[\Sigma_{2m'}]
\quad\mbox{ in }\quad L_0^\cx(R,\alpha,u)$$
if and only if there exist 
$n\in\N$, $\Xi\in M_{2(n+m+m')}$ and $\Delta\in \gl_{2(n+m+m')}$ such that
$$\Gamma\perp\Sigma_{2(n+m')}=
\Delta^\alpha(\Gamma'\perp\Sigma_{2(n+m)})\Delta +\Xi-\Xi^\alpha u
\quad\mbox{ and }\quad [\Delta]\in\cx.$$
\e{punt}

We conclude this section by stating some definitions and facts from algebraic
$K$- and $L$-theory needed in the sequel.

\be{thm} \label{exakring}{\rm \cite[Theorem 3]{Wall;lfound}}
Given an abelian group $A$ and an involution $t\colon A\ra A$
the Tate-cohomology groups $H^n(A;t)$ are defined by
$$H^n(A;t)\isdef{\ker(1-(-1)^nt)\over \im(1+(-1)^nt)}. $$
Suppose ${\cx_1}\subset {\cx_2}$ are involution invariant 
subgroups of $K_1(R)$,
then there exists an exact sequence 
$$
\halign{\quad\hfil$#$\hfil&\hfil$#$\hfil
&\hfil$#$\hfil&
 \hfil$#$\hfil&\hfil$#$\hfil
&\hfil$#$\hfil&\hfil$#$\hfil\cr
 H^1({\cx_2}/{\cx_1})&\mapright{\tilde{\tau}}&
L_0^{\cx_1}(R,\alpha,u)&
 \lra&L_0^{\cx_2}(R,\alpha,u)&\mapright{\tilde{\delta}}&
H^0({\cx_2}/{\cx_1})\cr
 \mapup{}&&&&&&\mapdown{}\cr
 L_1^{\cx_2}(R,\alpha,u)&&&&&&L_1^{\cx_1}(R,\alpha,-u)\cr
 \mapup{}&&&&&&\mapdown{}\cr
 L_1^{\cx_1}(R,\alpha,u)&&&&&&L_1^{\cx_2}(R,\alpha,-u)\cr
 \mapup{}&&&&&&\mapdown{}\cr
 H^0({\cx_2}/{\cx_1})&\lla&L_0^{\cx_2}(R,\alpha,-u)&
 \lla&L_0^{\cx_1}(R,\alpha,-u)&\lla&H^1({\cx_2}/{\cx_1})\cr}
$$
Here $\tilde{\tau}$ resp. $\tilde{\delta}$ is induced by $\tau$ resp.
$\delta$.
\e{thm}
\be{thm} {\rm \cite{w2}} {\sf Morita invariance.}\nl
If $(R,\alpha,u)$ is a ring with anti-structure and the
matrix ring $M_n(R)$ is equipped with the conjugate transpose anti-structure,
then $L_\varepsilon^*(M_n(R),\alpha,uI_n)$ is isomorphic to
$L_\varepsilon^*(R,\alpha,u)$.
\e{thm}
\be{thm} {\rm \cite{w2}} {\sf Scaling.}\nl
If $(R,\alpha,u)$ is a ring with anti-structure and $v$ is a unit in $R$,
then $L_\varepsilon^*(R,\alpha,u)$ is isomorphic to
$L_\varepsilon^*(R,\alpha',u')$, where $\alpha'(r)\isdef v\alpha(r)v^{-1}$
and $u'\isdef v\alpha(v^{-1})u$.
\e{thm}
\be{thm} \label{iadiciso}{\rm \cite[Lemma 5]{Wall;class3}}
Suppose $I$ is a two-sided ideal of $R$ such that
$R$ is complete in the $I$-adic topology.
If $\alpha(I)=I$, then $R/I$ can be equipped with an anti-structure
in an obvious way and the projection $R\ra R/I$ induces an isomorphism
$L_\varepsilon^h(R)\lra L_\varepsilon^h(R/I)$.
\e{thm}
\be{defi} \cite{Milnor}
Denote by $e_{ij}(a)\in E_n(R)$ the elementary matrix having the element
$a\in R$ at the $(i,j)$-entry.\nl
For $n\geq 3$ let $\st_n(R)$ be the group with the following presentation \nl
generators: one generator $x_{ij}(a)$ for every $e_{ij}(a)\in E_n(R)$\nl  
relations:\[ x_{ij}(a)x_{ij}(b)=x_{ij}(a+b)\]
\[[x_{ij}(a),x_{kl}(b)]=\left\{\begin{array}{l}
1 \mbox{\hspace{6ex} if } \; i\neq l,\; j\neq k\\
x_{il}(ab)  \mbox{\hspace{0.5ex} if }\;  j=k,\; i\neq l. \end{array}\right.\] 
The Steinberg group of $R$ denoted by $\st(R)$ is by definition the 
direct limit $$\lim_{\lra}\st_n(R),$$ where the limit is taken with
respect to the embeddings $\st_n(R)\hookrightarrow \st_{n+1}(R)$
coming from the embeddings 
$E_n(R)\hookrightarrow E_{n+1}(R)$
of definition~\ref{defk1}.
Since the relations for the $x_{ij}$ in $\st_n(R)$
also hold for the $e_{ij}$ in $E_n(R)$,
there is a natural homomorphism $\phi\colon \st_n(R)\ra E_n(R)$,
taking generators $x_{ij}(a)$ to $e_{ij}(a)$,
which in the limit gives rise to a homomorphism
$E(R)\ra \st(R)$.
The kernel of this last homomorphism 
is by definition the $K$-group $K_2R.$
\e{defi}
\be{lemma} \cite[theorem 5.1]{Milnor} 
$K_2R$ is the center of the Steinberg group.
\e{lemma}
\be{defi}
Denote by $\gl_{2\infty}(R)$ the direct limit of the groups
$\gl_{2n}(R)$ with respect to the embeddings
\[\gl_{2n}(R)\hookrightarrow \gl_{2(n+1)}(R)\mbox{ defined by }
\left(\matrix{A&B\cr C&D\cr}\right)\mapsto 
\left(\matrix{A&0&B&0\cr0&1&0&0\cr C&0&D&0\cr 0&0&0&1\cr}\right)\]
Similarly one defines $E_{2\infty}(R)$ and
correspondingly $\st_{2\infty}(R).$
\e{defi}
\be{punt}\label{involstek12} \cite[corollary 1.7]{Giffen;k2}
The anti-involutions $t_{\alpha,u}$ on the $\gl_{2n}(R)$ give rise
to anti-involutions on the $E_{2n}(R)$ which in turn lift to 
anti-involutions of $\st_{2n}(R).$ See definition~\ref{defantit} for formulas.
These provide for the following commutative diagram with exact rows
and vertical arrows (anti)-involutions:
$$\diagram{
0\lra&K_2(R)&\lra&\st_{2\infty}(R)&\lra&\gl_{2\infty}(R)&\lra&K_1(R)&\lra0\cr
&\mapdown{t_\alpha}&&\mapdown{t_{\alpha,u}}&&\mapdown{t_{\alpha,u}}&
&\mapdown{t_\alpha}&\cr
0\lra&K_2(R)&\lra&\st_{2\infty}(R)&\lra&\gl_{2\infty}(R)&\lra&K_1(R)&\lra0\cr
}$$
\e{punt}
\be{defi} \label{defgk2i}
Following \cite{Giffen;k2} one can construct a homomorphism
$$G\colon L_0^s(R)\ra H^1(K_2(R);t)$$ as follows: \nl
Let $$l=[\Gamma]-[\Sigma_{2m}]\in L_0^s(R)$$ and 
$X=\Gamma+\Gamma^\alpha u$.
Then $$U_{2m}^{-1}X\in E(R)
\mbox{ \ and \ } X^\alpha u=X.$$
Hence $$t_{\alpha,u}(U_{2m}^{-1}X)=
t_{\alpha,u}(X)\cdot t_{\alpha,u}(U_{2m}^{-1})=
U_{2m}^{-1}X^\alpha U_{2m}U_{2m}=U_{2m}^{-1}X^\alpha u=U_{2m}^{-1}X.$$
Now choose a lift $\gamma\in \st(R)$ of $U_{2m}^{-1}X$
and define $$G(l)\isdef[\gamma^{-1}t_{\alpha,u}\gamma]\in H^1(K_2(R);t).$$
It's not hard to check that $G$ is a well-defined homomorphism.
\e{defi}

\newpage
\section{The Arf invariant.}
\setcounter{altel}{0}
\setcounter{equation}{0}

In this section we define the main object of study in this thesis:
the Arf-groups. 

Suppose we are given a ring with anti-structure $(R,\alpha,u)$ and an 
involution invariant subgroup $\cx$ of $K_1(R).$
We will analyse the subgroup of $L_0^\cx(R,\alpha,u)$
consisting of all differences of classes of forms whose underlying
bilinear form is standard.

\be{defi}
Recall the considerations of \ref{obslgroup} and define 
$\Arf^\cx(R,\alpha,u)$ 
as the subgroup of $L_0^\cx(R,\alpha,u)$
generated by all elements
$$\plane{A,B}\isdef\left[\pmatrix{A&I_m\cr 0&B\cr}\right]
-\left[\pmatrix{0&I_m\cr0&0\cr}\right],$$
where $A,B\in\Lambda_m(R)\isdef\{X\in M_m(R)\mid X+X^\alpha u=0\}.$\nl
Note that 
$$\pmatrix{A&I_m\cr 0&B\cr}\quad\mbox{ and }\quad\pmatrix{0&I_m\cr 0&0\cr}$$
both belong
to ${\cal N}_{2m}(R).$\nl
Further we define $\Gamma_m(R)\isdef\{X-X^\alpha u\mid X\in M_m(R)\}$.
\e{defi}

\be{lemma}\label{lemmadiag}
All elements $\plane{A,B}$ of $\Arf^\cx(R,\alpha,u)$ can be 
written in the form:
\[\plane{A,B}=\sum_{i=1}^m\plane{A_{ii},B_{ii}}.\]
\e{lemma}
\be{proof}
Recall \ref{obslgroup}.
Since $A+A^\alpha u=B+B^\alpha u=0$ we find
\be{eqnarray*}
(A,B)&=&\left[\pmatrix{A&I_m\cr 0&B\cr}\right]-
        \left[\pmatrix{0&I_m\cr 0&0\cr}\right]\\
&=&\sum_{i=1}^m\left[\pmatrix{A_{ii}&1\cr0&B_{ii}\cr}\right]-
               \left[\pmatrix{0&1\cr0&0\cr}\right]\\
&=&\sum_{i=1}^m(A_{ii},B_{ii})
\e{eqnarray*}
\e{proof}
\be{prop}\label{proparfrel}\nl
Suppose we are given $A,B\in\Lambda_m(R)$ and $A',B'\in\Lambda_{m'}(R).$ 
Then
$$\plane{A,B}=\plane{A',B'} \mbox{ \ in \ }\quad  \Arf^\cx(R,\alpha,u),$$ 
if and only if there exist 
$$n\in\N \quad \mbox{ \ and \ } \quad
\pmatrix{X&Y\cr Z&T\cr}\in \gl_{2(n+m+m')}(R) \quad\mbox{ with \ }\quad
\left[\pmatrix{X&Y\cr Z&T\cr}\right]\in\cx,$$
such that
$$\be{array}{l}
A'=X^\alpha AX+X^\alpha Z+Z^\alpha BZ \pmod{\Gamma_{n+m+m'}(R)},\\ \\
B'=Y^\alpha AY+Y^\alpha T+T^\alpha BT \pmod{\Gamma_{n+m+m'}(R)}
\mbox{ \ and }\\ \\
t_{\alpha,u}\pmatrix{X&Y\cr Z&T\cr}=\pmatrix{X&Y\cr Z&T\cr}^{-1}.
\e{array}$$
Here $A,B,A',B'$ are considered to be elements of $M_{n+m+m'}(R)$,
by the embeddings $ M_k(R)\hookrightarrow M_{k+1}(R)$ defined by
$$\pmatrix{C\cr}\injarrow\pmatrix{C&0\cr0&0\cr}.$$  
\e{prop}
\be{proof}
Regarding the final assertion of \ref{obslgroup} 
it suffices to make the following
statements. Define $k\isdef n+m+m'$.
\[\left(\begin{array}{cc}X&Y\\Z&T\end{array}\right)^\alpha
\left(\begin{array}{cc}A&I_k\\0&B\end{array}\right)
\left(\begin{array}{cc}X&Y\\Z&T\end{array}\right)
\mbox{ takes the form }
\left(\begin{array}{cc}A'&I_k\\0&B'\end{array}\right)\]
(mod $\Gamma_{2k}(R)$) precisely when the difference
\[\left(\begin{array}{cc}X^\alpha AX+X^\alpha Z+Z^\alpha BZ
&X^\alpha AY+X^\alpha T+Z^\alpha BT\\
Y^\alpha AX+Y^\alpha Z+T^\alpha BT&Y^\alpha AY+Y^\alpha T+T^\alpha BT
\end{array}\right)-
\left(\begin{array}{cc}A'&I_k\\0&B'\end{array}\right)\]
belongs to $\Gamma_{2k}(R).$
From the fact that the matrices $A,B,A',B'$ in this expression belong
to $\Lambda_{2k}(R)$ we deduce:
\[\left\{\begin{array}{l}X^\alpha T+Z^\alpha uY=1\\
X^\alpha Z+Z^\alpha uX=0\\
Y^\alpha T+T^\alpha uY=0\end{array}\right.\]
This is equivalent to
\[\left(t_{\alpha,u}\left(\begin{array}{cc}X&Y\\Z&T\end{array}\right)\right)
\left(\begin{array}{cc}X&Y\\Z&T\end{array}\right)=1.\]
\e{proof}
We will give a presentation for the groups $\Arf^\cx(R,\alpha,u)$ in 
the next theorem.
Although our definition of the Arf- and $L$-groups 
is a priori quite different from the one in \cite{Clauwens;arf},
the presentation is nearly the same. 
We refer to \cite{Bak} for a comparison of the various $L$-groups.
Moreover this presentation is not quite the same as the one in 
\cite{Clauwens;arf},
because our $u$ is not necessarily central in $R$. 
At least not yet.
\be{thm} \label{thmarfgr}{\rm Compare \cite{Clauwens;arf}}.\nl
As abelian group $Arf^\cx(R,\alpha,u)$ has the following presentation:
\halign{\hfil#&\quad#\hfil&\quad#\hfill\cr
generators:\phantom{1)}
&$\plane{a,b}$&where $a,b\in \Lambda_1(R)$\cr 
relations:
1)&$\plane{a,b_1+b_2}=\plane{a,b_1}+\plane{a,b_2}$&
for all $a,b_1,b_2\in\Lambda_1(R)$\cr
2)&$\plane{a_1+a_2,b}=\plane{a_1,b}+\plane{a_2,b}$&
for all $a_1,a_2,b\in\Lambda_1(R)$\cr
3)&$\plane{a,b}=\plane{b,uau^{-1}}$&
for all $a,b\in\Lambda_1(R)$\cr
4)&$\plane{a,b}=0$&for all $a\in\Lambda_1(R),\;\;b\in \Gamma_1(R)$\cr
5)&$\plane{a,\alpha(x)bx}=\plane{xa\alpha^{-1}(x),b}$&
for all $a,b\in\Lambda_1(R),\;\;x\in R$\cr
6)&$\plane{a,b}=\plane{a,ba\alpha^{-1}(b)}.$&for all $a,b\in\Lambda_1(R)$\cr
7)&$\sum_{i=1}^n\plane{(X^\alpha Z)_{ii},(Y^\alpha T)_{ii}}=0$& 
  if $\pmatrix{X&Y\cr Z&T\cr}\in \gl_{2n}(R),$\cr
  &$t_{\alpha,u}\pmatrix{X&Y\cr Z&T\cr}=\pmatrix{X&Y\cr Z&T\cr}^{-1}$&
    and $\left[\pmatrix{X&Y\cr Z&T\cr}\right]\in\cx.$\cr}
\e{thm}
\be{proof}
$\Arf^\cx(R,\alpha,u)$ is generated by the $\plane{a,b}$ because of
lemma~\ref{lemmadiag}.
To prove the relations,
we will now exploit proposition~\ref{proparfrel}.\nl
Let $A,B\in\Lambda_m(R)$.\nl
Choosing
$$\pmatrix{X&Y\cr Z&T\cr}=\pmatrix{I_m&u^{-1}B\cr0&I_m\cr}$$
in proposition~\ref{proparfrel} 
yields 
\be{eqnarray*}
\plane{A,B}&=&\plane{A,(u^{-1}B)^\alpha Au^{-1}B+(u^{-1}B)^\alpha+B}\\
&=&\plane{A,B^\alpha uAu^{-1}B+B^\alpha u+B}\\
&=&\plane{A,BAB^{\alpha^{-1}}}.
\e{eqnarray*}
Taking $m=1$ this proves {\em 6}.\nl
Choosing $$\pmatrix{X&Y\cr Z&T\cr}=\pmatrix{I_m&0\cr A&I_m\cr}$$
in proposition~\ref{proparfrel} 
yields 
$$\plane{A,B}=\plane{A+A+A^\alpha BA,B}=\plane{A^\alpha BA,B}.$$
As a consequence $\plane{a,0}=\plane{0,a}=0$ 
for all $a\in\Lambda_1(R)$, which proves {\em 4}.\nl
Let $A',B',C',D'\in\Lambda_m(R)$ and $X'\in M_m(R)$.\nl
Choosing $$A=\pmatrix{{A'}&0\cr 0&{C'}\cr},\qquad 
B=\pmatrix{{B'}&0\cr 0&{D'}\cr}$$
and $$\pmatrix{X&Y\cr Z&T\cr}=
\pmatrix{I_m&0&0&{X'}\cr0&I_m&-u^{-1}{X'}^\alpha&0\cr
         0&0&I_m&0\cr0&0&0&I_m\cr}$$ in
proposition~\ref{proparfrel} yields
\be{eqnarray*}
\lefteqn{\left(\pmatrix{{A'}&0\cr 0&{C'}\cr},
\pmatrix{{B'}&0\cr 0&{D'}\cr}\right)}\hspace{2ex}\\
&=&\left(\pmatrix{{A'}&0\cr 0&{C'}\cr},\pmatrix{0&-u{X'}\cr {X'}^\alpha&0\cr}
\pmatrix{{A'}&0\cr 0&{C'}\cr}\pmatrix{0&{X'}\cr -u^{-1}{X'}^\alpha&0\cr}+\right.\\
& &\left.\pmatrix{0&-u{X'}\cr {X'}^\alpha&0\cr}+\pmatrix{{B'}&0\cr 0&{D'}\cr}\right)\\
&=&\left(\pmatrix{{A'}&0\cr 0&{C'}\cr},
\pmatrix{u{X'}{C'}u^{-1}{X'}^\alpha&0\cr0&{X'}^\alpha{A'}{X'}\cr}+
\pmatrix{{B'}&0\cr 0&{D'}\cr}\right)
\e{eqnarray*}
Hence 
\be{equation}
\plane{{A'},{B'}}+\plane{{C'},{D'}}=\plane{{A'},{B'}+
u{X'}{C'}u^{-1}{X'}^\alpha}+
\plane{{C'},{X'}^\alpha {A'}{X'}+{D'}}.\label{xeq}
\e{equation}
First choose ${C'}=u^{-1}{B'}u$, ${D'}=0$ and ${X'}=1$ to obtain
$$\plane{{A'},{B'}}=\plane{u^{-1}{B'}u,{A'}},$$ which proves {\em 3}.\nl
Then choose ${A'}={D'}$ and ${X'}=1$ to obtain
$$\plane{{A'},{B'}}+\plane{{C'},{A'}}=\plane{{A'},{B'}+u{C'}u^{-1}}$$
which by {\em 3} is equivalent to
$$\plane{{A'},{B'}}+\plane{{A'},u{C'}u^{-1}}=\plane{{A'},{B'}+u{C'}u^{-1}}.$$
This proves {\em 1}.\nl
Note that {\em 2} follows from and {\em 1} and {\em 3}.\nl
In order to verify {\em 5} we use  {\em 1, 2, 3} and {\em 4} to see that 
equation~\ref{xeq} comes down to
$$\plane{{A'},u{X'}{C'}u^{-1}{X'}^\alpha}=
\plane{{C'},{X'}^\alpha {A'}{X'}}.$$ 
But since $\plane{{A'},u{X'}{C'}u^{-1}{X'}^\alpha}=
\plane{{X'}{C'}{X'}^{\alpha^{-1}},{A'}},$
this proves {\em 5}.\nl
Note that all choices for $\pmatrix{X&Y\cr Z&T\cr}$ we have made so far
satisfy the conditions of proposition~\ref{proparfrel}.\nl
Finally suppose $\pmatrix{X&Y\cr Z&T\cr}$ agrees with the conditions
of {\em 7}.
To prove the theorem it suffices to show that
$$(X^\alpha AX+X^\alpha Z+Z^\alpha BZ,Y^\alpha AY+Y^\alpha T+T^\alpha BT)=
(A,B)+(X^\alpha Z,Y^\alpha T)$$
modulo the relations {\em 1} to {\em 6}. 
This is accomplished by using the relations for $X,Y.Z$ and $T$
listed in proposition~\ref{gqkar}.
We equate 
$$(X^\alpha AX,Y^\alpha AY)=(YX^\alpha AXY^{\alpha^{-1}},A)=
(XY^{\alpha^{-1}},A)$$
and in the same fashion
$$(Z^\alpha BZ,T^\alpha BT)=(ZT^{\alpha^{-1}},B).$$
Further
\be{eqnarray*}
\lefteqn{(X^\alpha AX,Y^\alpha T)+(X^\alpha Z,Y^\alpha AY)}\\
&=&(A,X^{\alpha^2}Y^\alpha TX^\alpha)+(YX^\alpha ZY^{\alpha^{-1}},A)\\
&=&(u^{-1}X^{\alpha^2}Y^\alpha TX^\alpha u,A)+
(YX^\alpha ZY^{\alpha^{-1}},A)\\
&=&(XY^{\alpha^{-1}},A)
\e{eqnarray*}
and analogously
$$(Z^\alpha BZ,Y^\alpha T)+(X^\alpha Z,T^\alpha BT)=(ZT^{\alpha^{-1}},B).$$
Finally we have
\be{eqnarray*}
\lefteqn{(X^\alpha AX,T^\alpha BT)+(Z^\alpha BZ,Y^\alpha AY)}\\
&=&(A,X^{\alpha^2}T^\alpha BTX^\alpha)+
              (YZ^\alpha BZY^{\alpha^{-1}},A)\\
&=&(A,X^{\alpha^2}T^\alpha BTX^\alpha)+
              (A,uYZ^\alpha BZY^{\alpha^{-1}}u^{-1})\\
&=&(A,X^{\alpha^2}T^\alpha BTX^\alpha)+
              (A,(1-X^{\alpha^2}T^\alpha)B(1-TX^\alpha))\\
&=&(A,B).
\e{eqnarray*}
This completes the proof.
\e{proof}
\be{thm}
There is a well-defined homomorphism, called Arf invariant
$$\omega\colon\Arf^\cx(R,\alpha,u)\ra R/\kappa(R),$$
defined by 
$$\plane{A,B}\mapsto\left[Tr(A^\alpha B)\right].$$
Here $\kappa(R)$ denotes the additive subgroup of $R$ generated by
$$\{x+x^2,y+\alpha(y)\mid x,y\in R\}$$
Observe that $xy-yx,2x\in\kappa(R)$ for all $x,y\in R$.
\e{thm}
\be{proof}
Analogous to the proof of \cite[theorem 2]{Clauwens;arf}.
\e{proof} 
\be{defi}
For every group $G$ we define 
$$L^{s,h}(G)\isdef L_0^{s,h}(\mbbf\,_2[G],\alpha,1)$$
and correspondingly
$$\Arf^{s,h}(G)\isdef\Arf^{s,h}(\mbbf\:_2[G],\alpha,1),$$
where $\alpha$ is determined by $\alpha(g)\isdef g^{-1}$ for all $g\in G$.
Further we define 
$$K(G)\isdef
\frac{\mbbf\,_2[G]}{\kappa(\mbbf\,_2[G])}.$$
\e{defi}
\be{remark} 
See also \cite{Clauwens;arf}.\label{remarfgrel}
From the presentation of theorem~\ref{thmarfgr} 
we deduce that $\Arf^{s,h}(G)$ 
is generated by all
$\plane{g,h}$ with $g,h\in{}_2G\isdef\{x\in G\mid x^2=1\}$ and that
the following relations hold:
\be{eqnarray*}
\plane{g,h}&=&\plane{h,g}\\
\plane{g,h}&=&\plane{xgx^{-1},xhx^{-1}}\quad \mbox{ for all }\quad x\in G\\
\plane{g,h}&=&\plane{g,hgh}.
\e{eqnarray*}
The value group $K(G)$ of the Arf invariant $\Arf^{s,h}(G)\lra K(G),$
is in fact the $\mbbf\,_2$-vectorspace generated by 
the quotient set $\kl(G)\isdef G/\!\sim$, 
where $\sim$ denotes the equivalence relation on
$G$ generated by $g\sim g^{-1}$, $g\sim hgh^{-1}$ and $g\sim g^2$.
\e{remark}
\be{thm}
The Arf invariant $\Arf^{s,h}(G)\ra K(G)$ is injective 
whenever $G$ is a finite group.
\e{thm}
\be{proof}
We refer to \cite{Clauwens;arf} for the proof.
\e{proof}
We will revert to these theorems later on.
\be{lemma}\label{centrel}
Let $a$, $b$ and $c$ be elements of order two in a group $G$ and assume
that $c$ commutes with $a$ and $b$.
Then the relation 
$$\plane{a,bc}=\plane{a,b}$$
holds in $\Arf^{s,h}(G)$.
\e{lemma}
\be{proof}
$\plane{a,bc}=\plane{a,bcabc}=\plane{a,bab}=\plane{a,b}.$
\e{proof}
\be{nitel}{Example}
Let $G$ be the group with presentation
$$\langle X,S\mid X^{12}=S^2=1,\, SXS=X^5\rangle.$$
So $G$ is a semidirect product of the group of order $2$ and the cyclic
group of order $12$.    
\be{prop}
The elements $\plane{1,1},\,\plane{X^2S,S}$ form a basis for
$\Arf^{s,h}(G)$.
\e{prop}
\be{proof}
A little computation yields $\kl(G)=\left\{[1],[X]\right\}$.
The Arf invariant is injective and maps $\plane{1,1}$ to $[1]$ 
and $\plane{X^2S,S}$ to
$[X^2]=[X]$, hence the assertion is true.
\e{proof}
\e{nitel}
The following example is meant to illustrate how tricky 
manipulations with the relations in $\Arf^{s,h}(G)$ can be.
\be{nitel}{Example}
Let $G$ be the group with presentation
$$G\isdef\langle X,Y,S\mid S^2=(XS)^2=Y^{12}=1,\quad SYS=Y^5,
\quad XY=YX\rangle.$$
This group fits into the short exact sequence
$$1\lra C\times C_{12}\lra G\lra C_2\lra 1,$$
where $C_2$ has generator $S$, $C_{12}$ has generator $Y$ and
$C$ is the infinite cyclic group generated by $X$.
Actually $G$ is a semidirect product of $C_2$ and $C\times C_{12}$.
We show that the elements
$\plane{S,SX^2Y^2}$ and $\plane{SX,SX^3Y^2}$ of $\Arf^{s,h}(G)$ coincide.
The Arf invariant $\omega$ maps both elements to the class of $X^2Y^2$ in
$K(G)$.
We equate
\be{eqnarray*}
\plane{S,SX^2Y^2}&=&\plane{S,SX^4Y^4}\\
                 &=&\plane{S,SX^2Y^{8}}\\
                 &=&\plane{S,SXY^{4}}\\
                 &=&\plane{SXY^2SSXY^2,SXY^2SXY^{4}SXY^2}\\
                 &=&\plane{SX^2Y^4,SX}\\
                 &=&\plane{SX,SX^2Y^4}\\
                 &=&\plane{SX,SX^3Y^8}\\
                 &=&\plane{SX,SX^5Y^4}\\
                 &=&\plane{SX,SX^3Y^2}.
\e{eqnarray*}
Since the Arf invariant maps both 
$\plane{S,SY^2}$ and $\plane{SX,SXY^2}$ to the class of $Y^2$ in $K(G)$,
one might conjecture that these elements are equal too,
but this is false.
\e{nitel}
\be{nitel}{Example}
Let $G$ be the group with presentation
$$G\isdef\langle Y,S\mid S^2=(YS)^4=(Y^2S)^2=1\rangle.$$
This group is actually an extension of the infinite cyclic group 
by the
dihedral group $D_4$:
$$\diagram{
1\lra C\lra &G\lra D_4&\lra 1\cr
&\be{array}{l}
S\longmapsto \sigma\\
Y\longmapsto \sigma\tau\vspace{1mm}
\e{array}&\cr}$$
Here $C$ is the infinite cyclic group generated by $Y^2$ and
$D_4$ is the dihedral group with presentation
$$D_4=\langle\sigma,\tau\mid \sigma^2=(\sigma\tau)^2=\tau^4=1\rangle.$$
\be{prop}
The set 
$$\left\{\plane{1,1}\right\}\cup
\left\{\plane{Y^{4i+2}S,S}\mid i>0\right\}$$
constitutes a basis for $\Arf^{s,h}(G)$.
\e{prop}
\be{proof}
The elements of order 2 in $G$ are $Y^{2i}S$, $(YS)^2$ and $Y^{2i}S(YS)^2$.
Note that $(YS)^2$ is central in $G$.
So we may use lemma~\ref{centrel} to see that $\Arf^{s,h}(G)$ is 
generated by elements of the form
$\plane{Y^{2i}S,Y^{2j}S}$.
The identities 
\be{eqnarray*}
\plane{Y^{2i}S,Y^{2j}S}&=&
\plane{Y^{-2k}Y^{2i}SY^{2k},Y^{-2k}Y^{2j}SY^{2k}}\\
&=&\plane{Y^{2i-4k}S,Y^{2j-4k}S},\\
\plane{Y^{2i}S,Y^{2}S}&=&
\plane{Y^{2i}S,Y^{2}S(YS)^2}\\
&=&\plane{Y^{2i}S,YSY^{-1}}\\
&=&\plane{Y^{2i-1}SY,S}\\
&=&\plane{Y^{2i-2}S(SY)^2,S}\\
&=&\plane{Y^{2i-2}S,S},\\
\plane{Y^{4i}S,S}&=&
\plane{Y^{2i}SSY^{2i}S,S}\\
&=&\plane{Y^{2i}S,S},\\
\plane{Y^{2i}S,S}&=&
\plane{SY^{2i}SS,S}\\
&=&\plane{Y^{-2i}S,S}
\e{eqnarray*}
show that $\left\{\plane{1,1}\right\}\cup
\left\{\plane{Y^{4i+2}S,S}\mid i>0\right\}$ is a set of 
generators for $\Arf^{s,h}(G)$.
We use the Arf invariant $\Arf^{s,h}(G) \ra K(G)$
to prove that these elements are independent.
It is easy to verify that $$\kl(G)=\left\{[1]\right\}\cup
\left\{[Y^{2i+1}]\mid i>0\right\}$$
by writing down all generating relations in $\kl(G)$.\nl
The Arf invariant maps $\plane{1,1}$ to $[1]$ and $\plane{Y^{4i+2}S,S}$ to
$[Y^{4i+2}]=[Y^{2i+1}]$.
This proves the assertion.
\e{proof}
\e{nitel}
\be{nitel}{Example}
Let $G$ be the group with presentation
$$G\isdef\langle X,Y,S\mid S^2=(XS)^2=(YS)^4=(Y^2S)^2=1,\quad XY=YX\rangle.$$
This group is actually an extension of the free abelian group $A$ of rank 2
by the dihedral group $D_4$:
$$\diagram{1\lra&A&\lra &G&\mapright{\pi}&D_4&\lra 1\cr}$$
where 
$\pi(S)\isdef\sigma$,
$\pi(X)\isdef1$,
$\pi(Y)\isdef\sigma\tau$ and
$A$ is generated by $X$ and $Y^2$.
\be{prop}
$\Arf^{s,h}(G)$ is generated by 
\be{eqnarray*}
\left\{\plane{1,1}\right\}&\cup&
\left\{\plane{X^{2i+1}Y^{2j}S,S}\mid i\geq 0\right\}\\&\cup&
\left\{\plane{X^{2i}Y^{4j+2}S,S}\mid j\geq 0\right\}\\&\cup&
\left\{\plane{X^{2i+1}Y^{4j+2}S,XS}\mid j\geq 0\right\}.
\e{eqnarray*}
\e{prop}
\be{proof}
The elements of order 2 in $G$ are $X^iY^{2j}S$, $(YS)^2$ and 
$X^iY^{2j}S(YS)^2$.
Note that $(YS)^2$ is central in $G$ again.
So $\Arf^{s,h}(G)$ is generated by elements of the form
$\plane{X^{i}Y^{2j}S,X^{k}Y^{2l}S}$.
We may assume that $k,l\in\{0,1\}$ by the identity
$$\plane{X^iY^{2j}S,X^kY^{2l}S}
=\plane{X^{i-2m}Y^{2j-4n}S,X^{k-2m}Y^{2l-4n}S}.$$
We may even assume that $l=0$ by the relation
\be{eqnarray*}
\plane{X^iY^{2j}S,X^kY^{2}S}&=&
\plane{X^iY^{2j}S,X^kY^{2}S(YS)^2}\\
&=&\plane{X^iY^{2j}S,X^kYSY^{-1}}\\
&=&\plane{X^iY^{2j-1}SY,X^kS}\\
&=&\plane{X^iY^{2j-1}(SY)^3,X^kS}\\
&=&\plane{X^iY^{2j-2}S,X^kS}.
\e{eqnarray*}
When $k=0$ we may assume that $i$ or $j$ is odd:
\be{eqnarray*}
\plane{X^{2i}Y^{4j}S,S}&=&
\plane{X^iY^{2j}SSX^iY^{2j}S,S}\\
&=&\plane{X^iY^{2j}S,S}
\e{eqnarray*}
In this situation we may assume that one odd exponent is positive:
\be{eqnarray*}
\plane{X^iY^{2j}S,S}&=&\plane{SX^iY^{2j}SS,S}\\
&=&\plane{X^{-i}Y^{-2j}S,S}
\e{eqnarray*}
When $k=1$ we may assume that $i$ and $j$ are odd:
\be{eqnarray*} 
\plane{X^{2i}Y^{2j}S,XS}&=&
\plane{Y^{2j}S,X^{-2i+1}S}\\
&=&\plane{X^{2i-1}Y^{2j}S,S}\\
\plane{X^{2i+1}Y^{4j}S,XS}&=&\plane{X^{i+1}Y^{2j}SXSX^{i+1}Y^{2j}S,XS}\\
&=&\plane{X^{i+1}Y^{2j}S,XS}
\e{eqnarray*} 
And finally, we may assume that j is positive:
\be{eqnarray*} 
\plane{X^{i}Y^{2j}S,XS}&=&
\plane{XSX^{i}Y^{2j}SXS,XS}\\
&=&\plane{X^{-i+2}Y^{-2j}S,XS}
\e{eqnarray*} 
This proves the proposition.
\e{proof}
\e{nitel}
\be{nitel}{Example}
Let $G$ be the group with presentation
$$G\isdef\langle X,Y,Z\mid X^2=Y^2=Z^2=(XY)^3=(YZ)^{3}=(XZ)^3=1\rangle.$$
This group is known as the affine Weyl group $\widetilde{A_2}$.\nl
Define $U\isdef XYZY$, $V\isdef YXZX$ and $W\isdef ZXYX$.
Then 
$UVW=1$ and
$U$, $V$ and $W$ commute.
The subgroup $H$ of $G$ generated by $U$, $V$ and $W$ is normal since,
$$\be{array}{ll}
XUX=U^{-1}&XVX=W^{-1}\\
YUY=W^{-1}&YVY=V^{-1}\\
ZUZ=V^{-1}&ZWZ=W^{-1}.
\e{array}$$
Further $G/H\cong S_3=\langle x,y\mid x^2=y^2=(xy)^3=1\rangle$.
These groups fit into the short exact sequence
$$1\lra H\lra G\rightleftmaps{}{\alpha} S_3\lra 1,$$
which splits by $\alpha(x)\isdef X$ and $\alpha(y)\isdef Y$.
Thus $G$ is actually a semidirect product of $S_3$ and $H$.
\be{prop}
The elements
$$\plane{1,1},\quad\plane{X,Y},\quad\plane{Y,Z},\quad\plane{X,Z}\quad
\mbox{ and }\quad\plane{XU^i,X}\quad\mbox{$ i>0 $\, odd}$$
form a basis for $\Arf^{s,h}(G)$. 
\e{prop}
\be{proof}
We merely sketch the proof.\nl
The elements of order two in $G$ are
$XU^i$, $YV^i$ and $XYXW^i$.
So in $\Arf^{s,h}(G)$ one has the following types of elements.
\be{enumerate}
\item
$\plane{XU^i,XU^j}$
\item
$\plane{XU^i,YV^j}$
\item
$\plane{XU^i,XYXW^j}$
\item
$\plane{YV^i,YV^j}$
\item
$\plane{YV^i,XYXW^j}$
\item
$\plane{XYXW^i,XYXW^j}$
\e{enumerate}
We prove that all of these elements are actually of the desired type by
using the relations
\be{enumerate}
\item[ ]
$XU^i=V^iXV^{-i}=W^iXW^{-i},$
\item[ ]
$YV^i=U^iYU^{-i}=W^iYW^{-i},$
\item[ ]
$XYXYV^iXYX=XU^{-i}.$
\e{enumerate}
\be{enumerate}
\item
Conjugation by $W^{-j}$ yields $\plane{XU^i,XU^j}=\plane{XU^{i-j},X}$.\nl
And further
\be{enumerate}
\item[ ]
$\plane{XU^i,X}=\plane{XU^iXXU^i,X}=\plane{XU^{2i},X},$ 
\item[ ]
$\plane{XU^i,X}=\plane{XXU^iX,X}=\plane{XU^{-i},X}.$
\e{enumerate}
\item
Conjugation by $U^{-j}$ yields $\plane{XU^i,YV^j}=\plane{XU^{i+2j},Y}$.
But because
$$\plane{XU^i,Y}=\plane{U^{-1}WXU^iW^{-1}U,U^{-1}WYW^{-1}U}=
\plane{XU^{i+3},Y},$$
only the elements 
\be{enumerate}
\item[ ]
$\plane{X,Y}$,
\item[ ]
$\plane{XU,Y}=\plane{YZY,Y}=\plane{Y,Z}$ \ and 
\item[ ]
$\plane{XU^{-1},Y}=\plane{XYZYX,Y}=\plane{Z,YXYXY}=\plane{X,Z}$ remain. 
\e{enumerate}
\item
Conjugation by $X$ yields $\plane{XU^i,XYXW^j}=\plane{XU^{-i},YV^{-j}}$.
\item
Conjugation by $XYX$ yields $\plane{YV^i,YV^j}=\plane{XU^{-i},XU^{-j}}$.
\item
Conjugation by $Y$ yields $\plane{YV^i,XYXW^j}=\plane{YV^{-i},XU^{-j}}$.
\item
Conjugation by $X$ yields $\plane{XYXW^i,XYXW^j}=\plane{YV^{-i},YV^{-j}}$.
\e{enumerate}
We give a list of generating relations in $\kl(G)$.
\be{enumerate}
\item[$\cdot$]
$U^iV^j\sim U^{-i}V^{-j}\sim U^{2i}V^{2j}\sim U^{j-i}V^j\sim U^iV^{i-j}
\sim U^jV^i$ 
\item[$\cdot$]
$XU^iV^j\sim XV^{j}\sim U^{j}V^{2j}\sim U^{j}V^{-j}$
\item[$\cdot$]
$YU^iV^j\sim YU^{i}\sim U^{2i}V^{i}\sim U^{i}V^{-i}$
\item[$\cdot$]
$XYXU^iV^j\sim YU^{j-i}V^j\sim U^{i-j}V^{j-i}$ 
\item[$\cdot$]
$YXU^iV^j \sim XYU^{j-i}V^j$ 
\item[$\cdot$]
$XYU^iV^j \sim XYU^{i+1}V^{j-1}\sim XYU^{i+1}V^{j+2}$ 
\e{enumerate}
The Arf invariant maps
$$\cases{
\plane{1,1} & to\hspace{3ex} $[1]$\cr
\plane{X,Y} & to\hspace{3ex} $[XY]$\cr
\plane{X,Z} & to\hspace{3ex} $[XZ]=[XYU]$\cr
\plane{Y,Z} & to\hspace{3ex} $[YZ]=[XYU^{-1}]$\cr
\plane{XU^i,X} & to\hspace{3ex} $[U^i]$ \hspace{1ex} $i$ is 
positive and odd.\cr}$$
From the list of relations we see that these images are independent,
which proves the proposition
\e{proof}
\e{nitel}
\noindent We will review some of these examples in chapter IV.

\newpage
{\Large {\bf \be{center}
Chapter II \vspace{4mm}\\
New Invariants for \mbox{\boldmath $L$}-groups.
\e{center}}}
\vspace{6mm}

\setcounter{section}{0}
\section{Extension of the anti-structure to the ring of formal power series.}
\setcounter{altel}{0}
\setcounter{equation}{0}

To construct new invariants we start by extending a given anti-structure
on a ring $R$ to the ring of formal power series $R[[T]]$, in a highly
non-trivial manner. \nl
The fact that the projection $R[[T]]\ra R$ induces an isomorphism of 
the associated $L$-groups, enables us to
build new invariants. 

\be{defi}
Suppose we are given a ring with antistructure
$(R,\alpha,u)$.\nl
For every $n\in\N\cup\{\infty\}$ we define
\be{eqnarray*}
R_n&\isdef&\cases{ 
R[T]/(T^{n+1})\,, & the truncated polynomial ring, if $n\in\N$\cr
R[[T]]\,, & the ring of formal power series, if $n=\infty,$\cr}\\
\I_n&\isdef& TR_n,
\mbox{ the two-sided ideal of $R_n$ generated by the class of $T$},\\
u_n&\isdef& u(1+T).
\e{eqnarray*}
Note that the class of $T$ in $R_n$ is also denoted by $T$.
Now we extend the anti-structure on $R$ to an anti-structure on $R_n$ by 
the formula
\[
\alpha\left(\sum a_kT^k\right)\isdef
\sum\alpha(a_k)\left({-T\over1+T}\right)^k.
\]
\e{defi}
\be{lemma}
For every $n\in\N\cup\{\infty\}$
\begin{enumerate}
 \item
  $(R_n,\alpha,u_n)$ is a ring with antistructure.
 \item
 $\I_n$ is an involution invariant two-sided ideal of $R_n$, i.e. 
 $\alpha(\I_n)=\I_n$.
 \item
 $R_n$ is complete in the $\I_n$-adic topology.
 \item
 The projection $R_n\ra R$ splits and $\alpha$ respects this splitting.
\end{enumerate}
\e{lemma}
\be{proof}
The proof is trivial and therefore omitted.
\e{proof}

\newpage
\section{Construction of the invariants $\omega_1^{s,h}$ and $\omega_2$.}
\setcounter{altel}{0}
\setcounter{equation}{0}

\be{punt}
In algebraic $K$-theory on has functors
$$K_i\colon \quad\mbox{{\sl category of ideals\/}}\ra 
\quad\mbox{{\sl category of abelian groups}}$$
for every $i\in \N$.
The {\sl category of ideals\/} is the category with \nl
objects: pairs $(R,I)$ consisting of a ring $R$ and a two-sided ideal $I$ 
of $R$\nl
morphisms: $f\colon (R,I)\ra (S,J)$ are the ringhomomorphisms $f\colon R\ra S$
satisfying $f(I)\subseteq J$.\nl
The groups $K_i(R)\isdef K_i(R,R)$ are the ones we already came across
in the first chapter.
For every pair $(R,I)$ there exists a long exact sequence
$$\cdots\ra K_{i+1}(R/I)\ra K_i(R,I)\ra K_i(R)\ra K_i(R/I)\ra\cdots$$
\e{punt}
\be{punt}
Let $(R,\alpha,u)$ be a ring with anti-structure and
$(R_n,\alpha,u_n)$ the associated extension.
Since the projection $R_n\ra R$ splits, we have
$$K_i(R_n)\cong K_i(R)\oplus K_i(R_n,\I_n)$$
by the functoriality of the $K_i$.
The involutions $t_\alpha$ on $K_1(R_n)$ 
and $K_2(R_n)$ induced by $\alpha$ respect this splitting.
Consequently, the Tate cohomology groups split accordingly:
$$H^{0,1}(K_i(R_n))\cong H^{0,1}(K_i(R))\oplus H^{0,1}(K_i(R_n,\I_n)).$$
\e{punt}
\be{thm}\label{thminvomega1}
The following periodic sequence is exact.
$$ \halign{\quad\hfil$#$\hfil&\hfil$#$\hfil&\hfil$#$\hfil&
 \hfil$#$\hfil&\hfil$#$\hfil&\hfil$#$\hfil&\hfil$#$\hfil\cr
 H^1(K_1(R_n,\I_n))&\buildrel\tilde\tau\over\lra&L_0^s(R_n,\conj,u_n)&
 \lra&L_0^s(R,\conj,u)&\buildrel \omega_1^s\over\lra&H^0(K_1(R_n,\I_n))\cr
 \uparrow&&&&&&\da\cr
 L_1^s(R,\conj,u)&&&&&&L_1^s(R_n,\conj,-u_n)\cr
 \uparrow&&&&&&\da\cr
 L_1^s(R_n,\conj,u_n)&&&&&&L_1^s(R,\conj,-u)\cr
 \uparrow&&&&&&\da\cr
 H^0(K_1(R_n,\I_n))&\lla&L_0^s(R,\conj,-u)&
 \lla&L_0^s(R_n,\conj,-u_n)&\lla&H^1(K_1(R_n,\I_n))\cr}$$
Here $\tilde{\tau}$ is induced by the homomorphism $\tau$ of 
definition~\ref{deftau}
and $\omega_1^s$ is induced by the discriminant homomorphism.
\e{thm}
\be{proof}
From theorem~\ref{exakring} of chapter I we obtain the following 
commutative diagram with exact rows
$(\varepsilon =0,1)$
$$
\begin{array}{ccccc}
L_{1-\varepsilon}^h(R_n)&\ra H^{1-\varepsilon}(K_1R)\ra
&L_\varepsilon^{K_1(R_n,\I_n)}(R_n)&\ra L_\varepsilon^h(R_n)\ra
&H^\varepsilon(K_1R)\\
\da&\|&\da&\da&\|\\
L_{1-\varepsilon}^h(R)&\ra H^{1-\varepsilon}(K_1R)\ra&L_\varepsilon^s(R)&\ra L_\varepsilon^h(R)\ra&H^\varepsilon(K_1R)
\end{array}
$$
Theorem~\ref{iadiciso} of chapter I implies that 
$L_\varepsilon^h(R_n)\ra L_\varepsilon^h(R)$ is an isomorphism.
Consequently $L_\varepsilon^{K_1(R_n,\I_n)}(R_n)$ is isomorphic to 
$L_\varepsilon^s(R)$
by applying the five lemma to the diagram above.
When we insert this in the sequence of theorem~\ref{exakring} of chapter I
applied to the ring $R_n$
with $\cx_1={0}$ and $\cx_2=K_1(R_n,\I_n)$,
we obtain the desired  periodic exact sequence.
\e{proof}
\be{defi}\label{defomegaend}
Define $$\omega_1^h\colon L_0^h(R,\alpha,u)\lra H^0(K_1(R_n,\I_n);t_\alpha)$$
as the composition of homomorphisms
$$L_0^h(R,\alpha,u)\cong
L_0^h(R_n,\alpha,u_n)\mapright{\tilde{\delta}}
H^0(K_1(R_n);t_\alpha)\lra H^0(K_1(R_n,\I_n);t_\alpha),$$
where $\tilde{\delta}$ is induced by the discriminant homomorphism $\delta$.
Notice that $\omega_1^s$ factors through $\omega_1^h$.\nl
Define $d$ as the composition of homomorphisms
$$H^1(K_1(R_n,\I_n))\mapright{\widetilde{\tau}}
L_0^s(R_n,\conj,u_n)\mapright{G}H^1(K_2(R_n,\I_n)),$$
where $G$ denotes the homomorphism of definition~\ref{defgk2i}
of the first chapter.
\e{defi}
\be{lemma}\label{expld}
The map $d$ can explicitly be given by
$$d([X])=[\gamma^{-1}t_{\alpha,u}\gamma], \mbox{ for all }X\in GL(R),$$
where $\gamma\in \st(R)$ is a lift of $(t_{\alpha,u}X)X\in E(R)$.
\e{lemma}
\be{proof}
Immediate by the definitions of $G$ and $\tau$.
\e{proof}
\be{thm}
The homomorphism $G$ induces a homomorphism
$$\omega_2\colon \ker(\omega_1^s)\ra \coker(d).$$
\e{thm}
\be{proof} This is clear now in view of the exact sequence of 
theorem~\ref{thminvomega1} and 
definition~\ref{defomegaend}.
\e{proof}

\newpage
\section{Recognition of $\omega_1^h$.}
\setcounter{altel}{0}
\setcounter{equation}{0}

We now proceed to analyse $\omega_1^h$. It will turn out that $\omega_1^h$
is strongly related to the Arf invariant of the first chapter.
\be{prop}\label{prophk1}
Let $(R,\alpha,u)$ be a ring with anti-structure.
For all $\plane{a,b}\in \Arf^h(R,\alpha,u)$
\[\omega_1^h(\plane{a,b})=
\left[1+\frac{\alpha(a)bT^2}{1+T}\right]\in H^0(K_1(R_n,\I_n))\]
\e{prop}
\be{proof}
We may take 
$$\left[\pmatrix{a&1\cr0&b\cr}\right]-\left[\pmatrix{0&1\cr0&0\cr}\right]
\in L_0^h(R_n,\alpha,u_n)$$
as a lift of $\plane{a,b}\in\Arf^h(R,\alpha,u)\subseteq L_0^h(R,\alpha,u)$.
\begin{eqnarray*}
\omega_1^h(\plane{a,b})
&=&\left[\pmatrix{a&1\cr0&b\cr}+
         \pmatrix{\alpha(a) &0\cr1&\alpha(b)\cr}u(1+T)\right]\\
&&-\left[\pmatrix{0&1\cr0&0\cr}+\pmatrix{0&0\cr1&0\cr}u(1+T)\right]\\
&=&\left[\pmatrix{\alpha(a) uT&1\cr u(1+T)&\alpha(b) uT\cr}
         \pmatrix{0&1\cr u(1+T)&0\cr}^{-1}\right]\\
&=&\left[\pmatrix{1&\frac{\alpha(a)T}{1+T}\cr\alpha(b)uT&1\cr}\right]\\
&=&\left[\pmatrix{1&\frac{-\alpha(a)T}{1+T}\cr0&1\cr}
\pmatrix{1&\frac{\alpha(a)T}{1+T}\cr\alpha(b)uT&1\cr}
\pmatrix{1&0\cr-\alpha(b)uT&1\cr}\right]\\
&=&\left[1-\frac{\alpha(a)\alpha(b)uT^2}{1+T}\right]\\
&=&\left[1+\frac{\alpha(a) bT^2}{1+T}\right]\\
\end{eqnarray*}
\e{proof}
For the time being we will assume that $R$ is commutative 
and write $\conj$ instead of $\alpha$.
\be{lemma}
$$q\colon H^0(R;\conj)\lra H^0(R;\conj)$$
defined by $$[x]\longmapsto [x^2]$$ is a homomorphism.
\e{lemma}
\be{proof}
$q$ is well-defined:
\be{itemize}
\item[$\cdot$]
$x^2=\ol{x}^2$ for all $x\in R$, satisfying $x=\ol{x}$.
\item[$\cdot$]
$(x+\ol{x})^2=x^2+x\ol{x}+\ol{x}x+\ol{x}^2=(x^2+x\ol{x})+\ol{(x^2+x\ol{x})}$, 
for all $x\in R$.
\e{itemize}
$q$ is a homomorphism, since for all $x,y\in R$ satisfying $x=\ol{x}$, 
$y=\ol{y}$
\be{eqnarray*}
q([x+y])&=&[(x+y)^2]\\
&=&[x^2+xy+yx+y^2]\\
&=&[x^2+xy+\ol{xy}+y^2]\\
&=&[x^2+y^2]\\
&=&q([x])+q([y]).
\e{eqnarray*}
\e{proof}
\be{defi}
Define $C(R)\isdef\coker(1+q)$.
\e{defi}
\be{prop}
If $n$ is even $(\neq0)$ or $n=\infty$, then
$$\lambda\colon H^0(K_1(R_n,\I_n);t_\alpha)\lra C(R)$$ defined below is an 
isomorphism.
\e{prop}
\be{proof}
We denote by $1+\I_n$ the multiplicative group of units in $R_n$, which
are congruent to $1$ modulo $\I_n$.
According to \cite[theorem 3.2]{Bass-Murphy} 
the homomorphism $(1+\I_n)\ra K_1(R_n,\I_n)$ determined by 
the composition 
$$(1+\I_n)\subset(R_n)^*=\gl_1(R_n)\ra K_1(R_n,\I_n)$$ is an isomorphism.
Since this isomorphism respects the involutions we may and will identify
$H^0(K_1(R_n,\I_n);t)$ and $H^0(1+\I_n;\conj)$.\nl
Define $Z\isdef\{f\in1+\I_n\mid f=\overline{f}\}$
and $B\isdef\{g\overline{g}\mid g\in1+\I_n\}$.\nl
If $$f\equiv1+aT+bT^2\pmod{T^3}$$ for certain $a,b\in R$, then
$$\ol{f}\equiv 1-\ol aT+(\ol a+\ol b)T^2\pmod{T^3}$$ and
$$f\ol{f}\equiv 1+(a-\ol a)T+(\ol a-a\ol a+b+\ol b)T^2\pmod{T^3}.$$
So $f\in Z$ implies $a=\overline{b}-b$.
It is easy to verify that the map $Z \ra C(R)$ defined by 
$f\mapsto [b\ol b]$
vanishes on $B$ and induces a homomorphism
$$\lambda\colon H^0(1+\I_n;\conj) \ra C(R).$$
Define $$\mu\colon C(R)\ra H^0(1+\I_n;\conj)$$ by 
$$[z]\mapsto [1+zT^2/(1+T)].$$
First note that $1+zT^2/(1+T)\in Z$. We will prove that $\mu$ is 
well-defined.
If $[z]=0$ in $C(R)$, there exist $x,y\in R$ with $y=\ol y$, such that 
$z=x+\ol x+y+y^2$.\nl
Define 
$$f\isdef 1+zT^2/(1+T) \quad\mbox{ and }\quad  g\isdef 1+yT-(x+y)T^2,$$
then
$$g\ol g=1-(x+\ol x+y+y^2)T^2 \quad \mbox{ and } \quad 
fg\ol g\equiv1\pmod{T^3}.$$
So we may assume $f\equiv1 \pmod{T^3}$.\nl
We assert that $[h]=1$ for all $h\in Z$ satisfying $h\equiv1\pmod{T^3}$.\nl
By induction we assume $k>0$ and
$$h\equiv1+aT^{2k+1}+bT^{2k+2}\pmod{T^{2k+3}},$$ 
for certain $a,b\in R$.
Now 
$$\ol h\equiv1-\ol aT^{2k+1}+((2k+1)\ol a+\ol b)T^{2k+2}\pmod{T^{2k+3}}.$$
So $h\in Z$ implies $(2k+1)a=\overline b-b$ and $\overline a=-a.$\nl
Defining $$g\isdef 1+(b+ka)T^{2k+1}-(k+1)bT^{2k+2},$$
yields
\be{eqnarray*}
g\ol g&\equiv&1+((ka+b)-(k\ol a+\ol b))T^{2k+1}+\\
      &      &((2k+1)(k\ol a+\ol b)-(k+1)(b+\ol b))T^{2k+2}\\
&\equiv&1-aT^{2k+1}-bT^{2k+2}\pmod{T^{2k+3}}
\e{eqnarray*}
and
$$hg\overline g\equiv1\pmod{T^{2k+3}}.$$
By induction we find $[h]=1$.\nl
Thus $\mu$ is well-defined.
Finally we prove that $\mu=\lambda^{-1}$:\nl
For all $[z]\in C(R),$
$$\lambda\mu([z])=\lambda(1+zT^2/(1+T))=[z\ol z]=[z^2]=[z].$$ 
For all $f\isdef1+aT+bT^2+\cdots\in Z,$
$$\mu\lambda([f])=\mu([b\ol b])=\left[1+b\ol bT^2/(1+T)\right],$$
But since
$$f^{-1}(1+b\ol bT^2/(1+T))(1+\ol bT)\ol{(1+\ol bT)}\equiv1\pmod{T^3}$$
we may apply the same argument as before to see that $\mu\lambda([f])=[f].$
\e{proof}
\be{thm}
The composition of homomorphisms
$$\Arf^h(R,\conj,u)\subseteq 
L_0^h(R,\conj,u)\mapright{\omega_1^h}H^0(K_1(R_n,\I_n))\mapright{\lambda}C(R)
\lra R/\kappa(R),$$
is just the Arf invariant $\Arf^h(R,\conj,u)\ra R/\kappa(R)$
defined in section 2 of the first chapter.
Here $C(R)\ra R/\kappa(R)$ is induced by inclusion.
\e{thm}
\be{proof}
In view of proposition~\ref{prophk1} we have
\begin{eqnarray*}
\lambda\omega_1^h(\plane{a,b})
&=&\lambda\left(\left[1+\frac{\ol abT^2}{1+T}\right]\right) \\
&=&[\overline{a}b].
\end{eqnarray*}
The rest is clear.
\e{proof}
From now on $R$ is not necessarily commutative.
Let $(R,\conj,u)$ be a ring with anti-structure.
We wish to prove that the Arf invariant 
$$\Arf^h(R,\conj,u)\ra R/\kappa(R),$$ 
we dealt with in section 2 of chapter I, 
factors through the invariant
$$\omega_1^h\colon \Arf^h(R,\conj,u)\lra H^0(K_1(R_2,\I_2)).$$
Here follows an attempt to uncover the connection between 
$$R/\kappa(R)$$ and the Tate cohomology group $$H^0(K_1(R_2,\I_2)),$$
in the non-commutative case.
Let us fix the following notations.
\be{enumerate}
\item[$\cdot$]
$A$ is the truncated polynomial ring $R_2$. 
\item[$\cdot$]
$\I$ is the two-sided ideal of $A$ generated by $T$,
\item[$\cdot$]
$\conj\colon A\ra A$ is the extension of $\conj$ on $R$ to $A$ determined by
$$T\mapsto {-T\over 1+T}=-T+T^2,$$
i.e. $\ol{a+bT+cT^2}=\ol{a}-\ol{b}T+(\ol{b}+\ol{c})T^2$.
\item[$\cdot$]
$1+\I$ denotes the multiplicative group of units in $A$ which are 
congruent to $1$ modulo $\I$.
\item[$\cdot$]
We write $W=W(A,\I)$ for the subgroup of $1+\I$ generated by
the set $\{(1+ax)(1+xa)^{-1}\,\mid a\in A, x\in \I\}$.
According to \cite[theorem 2.1]{Swan} 
$W$ is the kernel of 
the surjection $1+\I\ra K_1(A,\I)$.
We will identify $K_1(A,\I)$ and $(1+\I)/W$.
\item[$\cdot$]
For all $r,s\in R$ we define $[r,s]\isdef rs-sr$.
And $R_{{\rm ab}}\isdef R/[R,R]$ the quotient of $R$ as an additive group 
by the subgroup generated by all $[r,s]$.
This is actually the Hochschild homology group $H_0(R)$.
\e{enumerate}
As we saw in section 2 of chapter II the anti-automorphism $\conj$ of $A$
induces an involution $t$ on the relative $K$-group $K_1(A,\I)$.
We want to investigate the structure of the Tate cohomology groups 
$H^0(K_1(A,\I))\cong H^0((1+\I)/W)$.
We proceed to take a close look at the group $W$.
\be{lemma}
Every element of $W$ has the form
$$1+\left(\sum_i[u_i,v_i]\right)T+
\left(\sum_k[r_k,s_k]+\sum_{i}u_iv_i[u_i,v_i]+\sum_{i<j}[u_i,v_i][u_j,v_j]
\right)T^2.$$
\e{lemma}
\be{proof}
Substituting $a=a_0+a_1T$ and $x=x_1T+x_2T^2$ in the expression 
$(1+ax)(1+xa)^{-1}$ yields
\be{eqnarray*}
\lefteqn{(1+(a_0+a_1T)(x_1T+x_2T^2))(1+(x_1T+x_2T^2)(a_0+a_1T))^{-1}}\\
&=&(1+a_0x_1T+(a_0x_2+a_1x_1)T^2)(1+x_1a_0T+(x_1a_1+x_2a_0)T^2)^{-1}\\
&=&(1+a_0x_1T+(a_0x_2+a_1x_1)T^2)
(1-x_1a_0T+((x_1a_0)^2-x_1a_1-x_2a_0)T^2)\\
&=&1+[a_0,x_1]T+([a_0,x_2]+[a_1,x_1]+[x_1,a_0]x_1a_0)T^2.
\e{eqnarray*}
When $a_0=0$ we obtain elements like $$1+[r,s]T^2$$ and modulo such elements
we find expressions of the form
$$1+[u,v]T+uv[u,v]T^2.$$
Note that $$(1+[u,v]T+uv[u,v]T^2)^{-1}=1+[v,u]T+vu[v,u]T^2.$$
Thus $W$ is generated by 
$$\left\{1+[u,v]T+uv[u,v]T^2,\,1+[r,s]T^2\,\mid r,s,u,v\in R\right\}.$$
Writing out a product of such elements yields the desired result.
\e{proof}

We also need the Hochschild homology group $H_1(R)$.
We refer to chapter III for the definitions.
The Hochschild homology group $H_1(R)$ and
the cyclic homology group $HC_1(R)$ are defined as: 
$$H_1(R)\isdef
\frac{\ker(b\colon R\te R\ra R)}{\im(b\colon R\te R\te R\ra R)}$$
$$HC_1(R)\isdef
\frac{\ker(b\colon R\te R\ra R)}{\im(b\colon R\te R\te R\ra R)+\im(1-x)}\,,$$
where 
$$b(u\te v)\isdef[u,v],\quad
b(u\te v\te w)\isdef uv\te w-u\te vw+wu\te v$$
$$\mbox{ and }\quad x\colon R\te R\ra R\te R
\mbox{ \ is defined by \ } x(u\te v)\isdef-v\te u.$$

\be{punt}
Define $\theta\colon R\te R\ra  R_{{\rm ab}}$ by
$$\theta(\sum_i u_i\te v_i)\isdef \sum_{i<j}[u_i,v_i][u_j,v_j]+
\sum_i u_iv_i[u_i,v_i].$$
$\theta$ is well-defined in the sense that the right-hand side
does not depend on the order of summation in $\sum_i u_i\te v_i$.
Observe that
$$\theta(x+y)=\theta(x)+\theta(y)+b(x)\cdot b(y).$$
for all $x,y\in R\te R$.
So the restriction of $\theta$ to $\ker(b)$ is a homomorphism.
Furthermore it is easy to verify that $\theta$ vanishes on $\im(b)$
and $\im(1-x)$.
Consequently $\theta$ induces a homomorphism 
$\theta'\colon HC_1(R)\ra  R_{{\rm ab}}$.
\e{punt}
In view of the preceding it is clear that the sequence
$$\diagram{
HC_1(R)&\mapright{\theta'}& R_{{\rm ab}}&\lra& K_1(A,\I)&\lra
& R_{{\rm ab}}&\lra 0\cr
&&[s]&\longmapsto&[1+sT^2]&&&\cr    
&&&&[1+aT+bT^2]&\longmapsto&[a]&\cr}$$
is exact.\nl
The anti-automorphism $\conj\colon R\ra R$ induces 
an involution on $ R_{{\rm ab}}$:
$$\ol{[u,v]}=[\ol{v},\ol{u}]$$
$$\left[\ol{\ol{r}}\right]=[uru^{-1}]=[r]\quad\mbox{ in } R_{{\rm ab}}.$$
Furthermore $\im(\theta')$ is invariant under this involution.
When we equip $ R_{{\rm ab}}$ on the left-hand side with this involution and $ R_{{\rm ab}}$
on the right-hand side with the involution $[a]\mapsto[-\ol{a}]$,
we obtain the short exact sequence of groups with involutions
$$0\lra\coker(\theta')\lra K_1(A,\I)\lra  R_{{\rm ab}}\lra0$$
which gives rise to the six-term exact sequence
$$\diagram{
H^0(\coker(\theta'))&\lra&H^0(K_1(A,\I))&\lra& H^1(R_{{\rm ab}})\cr
\mapup{\delta}&&&&\mapdown{}\cr
H^0(R_{{\rm ab}})&\lla&H^1(K_1(A,\I))&\lla&H^1(\coker(\theta')).\cr}$$
We compute the differential map 
$\delta\colon H^0( R_{{\rm ab}})\ra H^0(\coker(\theta'))$.
\be{lemma}\label{lemzeshrand}
If $[a]\in H^0( R_{{\rm ab}})$, i.e. $\ol{a}-a=b(x)$ for some $x\in R\te R$,
then $$\delta([a])=[a+a\ol{a}+\theta(x)].$$
\e{lemma}
\be{proof}
The element $1+aT$ is a lift of $a$ in $K_1(A,\I)$.
And in $K_1(A,\I)$ we have
\be{eqnarray*}
(1+aT)\ol{(1+aT)}&=&(1+aT)(1-\ol{a}T+\ol{a}T^2)\\
&=&1+(a-\ol{a})T+(\ol{a}-a\ol{a})T^2\\
&=&(1+(a-\ol{a})T+(\ol{a}-a\ol{a})T^2)(1+b(x)T+\theta(x)T^2)\\
&=&1+((a-\ol{a})(\ol{a}-a)+\ol{a}-a\ol{a}+\theta(x))T^2\\
&=&1+((a-\ol{a})(\ol{a}-a)+\ol{a}-a\ol{a}+\theta(x))T^2.
\e{eqnarray*}
But this is the image of 
\be{eqnarray*}
[(a-\ol{a})(\ol{a}-a)+\ol{a}-a\ol{a}+\theta(x)]&=&
[\ol{a}+\ol{a}a+\theta(x)]\\
&=&[a+a\ol{a}+\theta(x)]
\e{eqnarray*}
in $H^0(\coker(\theta')).$
\e{proof}
Now we specialize to the case that $R$ is the group ring $\Z[G]$
of an arbitrary group $G$.
\be{lemma}
$\theta'=0$
\e{lemma}
\be{proof}
Every cycle of $HC_1(R)$ can be written as 
$$\left[\sum_i g_i\te h_i\right],$$ by using the relation
$g\te h+h\te g=0$.
The condition for this element to be a cycle reads $\sum g_ih_i=\sum h_ig_i$.
Such an cycle can be decomposed as a sum of cycles of the form
$$[g\te h] \quad\mbox{with} \quad gh=hg$$
or of the form
$$\left[\sum_i^n g_i\te h_i\right] \quad \mbox{with}\quad
g_ih_i=\cases{h_{i+1}g_{i+1}& for $i<n$\cr
h_1g_1& for $i=n$\cr}.$$
The homomorphism $\theta'$ is obviously zero on elements of the first type.
As far as the second type is concerned we have the following identities
in $R_{{\rm ab}}$
\be{eqnarray*}
\theta'\left(\left[\sum g_i\te h_i\right]\right)&=&
\sum_{i<j}[g_i,h_i][g_j,h_j]+\sum_ig_ih_i[g_i,h_i]\\
&=&\sum_{i<j}g_ih_ig_jh_j+\sum_{i<j}h_ig_ih_jg_j+\sum_ig_ih_ig_ih_i+\\
& &-\sum_{i<j}h_ig_ig_jh_j-\sum_{i<j}g_ih_ih_jg_j-\sum_ig_ih_ih_ig_i\\
&=&\sum_{i<j}g_ih_ig_jh_j+\sum_{i<j}g_ih_ig_jh_j+\sum_ig_ih_ig_ih_i+\\
& &-\sum_{i<j}g_jh_jh_ig_i-\sum_{i<j}g_ih_ih_jg_j-\sum_ig_ih_ih_ig_i\\
&=&\sum_{i,j}g_ih_ig_jh_j-\sum_{i,j}g_jh_jh_ig_i\\
&=&\left(\sum g_ih_i\right)^2-
\left(\sum g_ih_i\right)\left(\sum h_ig_i\right)\\
&=&0
\e{eqnarray*}
This proves the lemma.
\e{proof}
The next move is to figure out what $\delta\colon 
H^0( R_{{\rm ab}})\lra H^0( R_{{\rm ab}})$ looks like in this case.
Suppose we are given an element $[a]\in H^0( R_{{\rm ab}})$.
Then we may assume that $a=\sum g_i$ by using the fact that 
$[g+g^{-1}]=0$ in $H^0( R_{{\rm ab}})$.
The condition for $a$ to be a cycle reads 
$$\sum g_i-g_i^{-1}=\sum h_j-h_j',$$
where $h_j\in G$ and $h_j'$ is a conjugate of $h_j$.
From this we conclude that every $g_i$ is conjugated to some $g_j^{-1}$.
Note that  $[g+h^{-1}g^{-1}h]=[g+g^{-1}]=0$ in $H^0( R_{{\rm ab}})$.
Thus it suffices to consider the case that 
$a=g$ where $g=h^{-1}g^{-1}h$.
We follow lemma~\ref{lemzeshrand}.
Now $g^{-1}-g=[h^{-1},gh]$, so 
\be{eqnarray*}
\delta([g])&=&[g+gg^{-1}+\theta([h^{-1}\te gh])]\\ 
           &=&[g+1+h^{-1}gh(g^{-1}-g)]\\ 
           &=&[g+1+g^{-1}(g^{-1}-g)]\\ 
           &=&[g+g^{-2}]\\ 
           &=&[g+g^{2}].
\e{eqnarray*}
As a consequence we have
$$\coker(\delta)=
\frac{\{a\in \Z[G]_{{\rm ab}}\mid a=\ol{a}\}}%
{\span\{g-h^{-1}gh,g_1+g_1^{-1},g_2+g_2^{2}\mid g_2\sim g_2^{-1}\}}.$$
Our main conclusion is that in the case of a group ring the invariant
$$\omega_1^h\colon \Arf^h(R,\conj,u)\lra H^0(K_1(A,\I))$$
factors through an injective homomorphism
$$\coker(\delta\colon H^0(R_{{\rm ab}})\ra H^0(R_{{\rm ab}}))\injarrow 
H^0(K_1(A,\I))$$ 
and that there is a homomorphism
$$\coker(\delta)\lra R/\kappa(R).$$

\newpage
\section{Computations on the invariant $\omega_2.$}
\setcounter{altel}{0}
\setcounter{equation}{0}

In order to study the invariant $\omega_2$, we wish to compute
the cokernel of the homomorphism
$$d\colon H^1(K_1(R_n,\I_n);t_\alpha)\ra H^1(K_2(R_n,\I_n);t_\alpha).$$
We confine our inquiries to the case where $R$ is commutative,
for then we have the following theorem at our disposal.
\be{thm}
Let $R$ be a commutative ring with identity and $I$ an ideal 
contained in the Jacobson radical of $R$.
Then $K_2(R,I)$ is isomorphic to the abelian group with
presentation:\vspace{1mm}
\halign{#&\quad#\hfill&\quad#\hfill&\quad#\hfill\cr
generators:
&$\denstein{a,b}$& with $a\in I$ or $b\in I$\vspace{1mm}\cr
relations:
&$\denstein{a,b}=-\denstein{b,a}$& if $a\in I$ or $b\in I$\vspace{1mm}\cr
&$\denstein{a,b}+\denstein{a,c}=\denstein{a,b+c-abc}$& if
$a\in I$ or $b,c\in I$\vspace{1mm}\cr
&$\denstein{a,bc}=\denstein{ab,c}+\denstein{ac,b}$& if $a\in I$
or $b\in I$ or $c\in I.$\vspace{1mm}\cr}
\noindent The isomorphism maps $\denstein{a,b}$ 
to the Dennis-Stein element
$\denstein{a,b}_\circ\in K_2(R,I)$. 
\e{thm}
\be{proof}
See \cite{Maazen-Stienstra,Keune}.
\e{proof}
A little digression seems in order. We refer to \cite[\S9]{Milnor}
and \cite{Dennis-Stein} for more background.\nl
Let $n>2$.\nl 
For any unit $r\in R$ one has the elements   
$w_{ij}(r)\isdef x_{ij}(r)x_{ji}(-r^{-1})x_{ij}(r)$ 
and $h_{ij}(r)\isdef w_{ij}(r)w_{ij}(-1)$ in $\st_{n}(R)$,
where $i$ and $j$ are distinct integers between $1$ and $n$.\nl
Further, for every couple of units $r,s\in R$, 
$$h_{ij}(rs)h_{ij}^{-1}(r)h_{ij}^{-1}(s)\in \st_{n}(R)$$
determines an element $\steinberg{r,s}$ in $K_2(R)$, 
which does not depend on $i$ or $j$.\nl
And for all $a,b\in R$ such that $1-ab$ is a unit of $R$, 
$$x_{ji}(-b(1-ab)^{-1})x_{ij}(-a)x_{ji}(b)x_{ij}(a(1-ab)^{-1})
h_{ij}^{-1}(1-ab)\in \st_{n}(R)$$
determines the Dennis-Stein element
$\denstein{a,b}_\circ\in K_2(R)$ which does not depend on $i$ or $j$ 
either. Note the sign conventions.\nl
In $K_2(R)$  the following relations hold,
whenever the left-hand side is defined.
\be{eqnarray*}
\steinberg{r_1r_2,s}&=&\steinberg{r_1,s}\steinberg{r_2,s}\\
\steinberg{r,s}&=&\steinberg{s,r}^{-1}\\
\steinberg{r,-r}&=&1\\
\steinberg{r,1-r}&=&1\\
\denstc{a,b}&=&\denstc{b,a}^{-1} \\
\denstc{a,b}\denstc{a,c}&=&\denstc{a,b+c-abc}\\
\denstc{a,bc}&=&\denstc{ab,c}\denstc{ac,b}\\
\denstc{0,a}&=&1 \\
\steinberg{r,s}&=&\denstc{(1-r)s^{-1},s}.
\e{eqnarray*}
Note that we used an additive notation in dealing
with the symbols $\denstein{\;,\;}$ 
and a multiplicative notation for the corresponding
Dennis-Stein elements $\denstc{\;,\;}$.
Nevertheless we will often omit the ${\scriptstyle \circ}$ .
\be{prop}\label{propk2inv}
Let $R$ be a commutative ring and $\conj\colon R\ra R$ an involution.
The involution $t$ on $K_2(R)$ induced by $\conj$ satisfies 
$$t(\denstc{a,b})=\denstc{\ol b,\ol a}.$$
\e{prop}
\be{proof}
We will work in $\st_{2n}(R)$.
We drop the decorations of the anti-involution on the Steinberg group
and simply write $t$.
From definition~\ref{defantit} and \ref{involstek12} of the 
first chapter we deduce
$$t(x_{ij}(a))=x_{n+j\,n+i}(\ol a),$$
provided that $i$ and $j$ do not exceed $n$.\nl
Thus $t(w_{12}(r))=w_{n+2\,n+1}(\ol r)$ and
\be{eqnarray}
t(h_{12}^{-1}(r))&=&w_{n+2\,n+1}^{-1}(\ol r)w_{n+2\,n+1}^{-1}(-1)\nonumber\\
                 &=&w_{n+2\,n+1}(-\ol r)w_{n+2\,n+1}(1)\\
                 &=&w_{n+1\,n+2}(\ol r^{-1})w_{n+1\,n+2}(-1)\\
                 &=&h_{n+1\,n+2}(\ol r^{-1})\nonumber
\e{eqnarray}
In (1) we used the relation $w_{ij}(r)=w_{ij}^{-1}(-r)$ and
(2) follows from the relation $w_{ij}(r)=w_{ji}(-r^{-1})$.
See \cite[lemma 9.5]{Milnor}.
Hence 
\be{eqnarray*}
t(\denstc{a,b})
&=&h_{n+1\,n+2}((1-\ol{ab})^{-1})
x_{n+2\,n+1}(\ol a(1-\ol{ab})^{-1})x_{n+1\,n+2}(\ol b)\cdot\\
&&x_{n+2\,n+1}(-\ol a)x_{n+1\,n+2}(-\ol b(1-\ol{ab})^{-1})\\
&=&h_{n+1\,n+2}((1-\ol{ab})^{-1})
\denstc{-\ol b,-\ol a}h_{n+1\,n+2}(1-\ol{ab})\\
&=&\denstc{-\ol b,-\ol a}\steinberg{(1-\ol{ab})^{-1},1-\ol{ab}}\\
&=&\denstc{\ol b,\ol a}\denstc{-\ol{ab},-1}\steinberg{(1-\ol{ab})^{-1},-1}\\
&=&\denstc{\ol b,\ol a}
\e{eqnarray*}
which proves the assertion.
\e{proof}
To make life more congenial, 
we will assume $R$ to carry some additional structure.
In that way $H^1(K_2(R_n,\I_n))$ becomes fairly accessible for computations 
by the techniques of \cite{Clauwens;k}.
The following definition occurs implicitly in \cite{Joyal} and 
\cite{Joyal;vec}.
It describes a notion of what one could call `partial $\lambda$-ring'.
\be{defi}\label{deftheta}
Let $R$ be a commutative ring with identity and $k\in\N\cup\{\infty\}.$
A $k\lambda$-ring structure on $R$ consists of
operations $\theta^p\colon R\ra R$, 
for every prime number $p\leq k$,
which satisfy the following conditions\vspace{1mm}
\halign{#\hfill&\quad#\hfill&\quad#\hfill\cr
1)& $\theta^p(1)=0$&for all \, $p\leq k$\vspace{1mm}\cr
2)& $\theta^p(a+b)=\theta^p(a)+\theta^p(b)
 +\sum_{k=1}^{p-1}\frac{1}{p}{p\choose k}a^kb^{p-k}$
 &for all \, $p\leq k$\vspace{1mm}\cr
3)& $\theta^p(ab)=\theta^p(a)b^p+\theta^p(b)a^p-p\theta^p(a)\theta^p(b)$
 &for all \, $p\leq k$\vspace{1mm}\cr
4)& $\theta^p(\psi^q(a))=\psi^q(\theta^p(a))$
 &for all \, $p,q\leq k$\vspace{1mm}\cr}
\noindent here  
$\psi^q$ is defined by $\psi^q(a)\isdef a^q-q\theta^q(a).$\nl
We then call $R$ an $k\lambda$-ring.
\e{defi}
\be{remark}
It is easy to verify that multiplication by $p$ transforms the equations 1
to 4 into
\halign{#\hfill&\quad#\hfill&\quad#\hfill\cr
1')& $\psi^p(1)=1$&for all \, $p\leq k$\vspace{1mm}\cr
2')& $\psi^p(a+b)=\psi^p(a)+\psi^p(b)$
 &for all \, $p\leq k$\vspace{1mm}\cr
3')& $\psi^p(ab)=\psi^p(a)\psi^p(b)$
 &for all \, $p\leq k$\vspace{1mm}\cr
4')& $\psi^p(\psi^q(a))=\psi^q(\psi^p(a))$
 &for all \, $p,q\leq k.$\vspace{1mm}\cr}
Thus the so called Adams operations $\psi_p$ are 
ringhomomorphisms, which satisfy the compatibility
conditions 4'.\nl
Conversely, if $R$ is a torsion-free commutative ring equipped with
$\psi_p$ satisfying 1' to 4' such that
$\psi_p(a)\equiv a^p\pmod{pR}$ for all $p\leq k$,
then $R$ becomes a $k\lambda$-ring in the 
obvious way and the $\psi_p$ are the associated Adams operations.\nl
As far as the references to \cite{Joyal} and \cite{Joyal;vec} are concerned,
a few remarks are in order.
\be{itemize}
\item[$\cdot$]
We point out the differences in sign conventions between 
the definition in \cite{Joyal;vec} and the one above.
\item[$\cdot$]
Condition 4 in our list is equivalent to what is called
the permutability of $\theta_p$ and $\theta_q$ in \cite{Joyal}.
\e{itemize}
\e{remark}
The terminology is explained by the following theorem.
\be{thm} {\rm \cite[theorem 3]{Joyal}. }
The notions $\lambda$-ring and $\infty\lambda$-ring coincide.
\e{thm}
\be{lemma}\label{lringext}
Any structure of $k\lambda$-ring on a ring $R$ admits 
a unique extension to the rings $R[T]$ and $R_n$ 
for all $n\in\N\cup\{\infty\}$, under the condition
that $\theta_p(T)=0$ for all $p\leq k$.
\e{lemma}
\be{proof}
There exists a unique $k\lambda$-ring structure on the ring of integers $\Z$
defined by $\psi_p\isdef 1$.\nl
Since the polynomial ring $\Z[T]$ has no torsion and the condition 
$\theta_p(T)=0$ implies $\psi_p(T)=T^p$, the formula
$\psi_p(\sum a_iT^i)=\sum a_iT^{ip}$, determines a unique
structure of $k\lambda$-ring on $\Z[T]$.\nl
We now call upon \cite[theorem 3]{Joyal;vec}, which reads as follows.
If $R_1$ and $R_2$ are $k\lambda$-rings, then $R_1\te R_2$ can be provided
with a unique structure of $k\lambda$-ring, such that the canonical
maps $R_1\ra R_1\te R_2$ and $R_2\ra R_1\te R_2$ preserve every $\theta_p$.
Applying this theorem in our situation, proves the assertion for the ring
$R[T]=R\te \Z[T]$. \nl
From condition 2 in definition~\ref{deftheta} we deduce that
$f\equiv g\pmod{T^lR[T]}$ implies 
$\theta_p(f)\equiv\theta_p(g)\pmod{T^lR[T]}$.
Consequently the $k\lambda$-ring structure on $R[T]$ extends
uniquely to the rings $R_n$ for all $n\in \N\cup\{\infty\}$.
\e{proof}
\be{defi}
Let $\delta\colon R\ra\Omega_R$ be the universal derivation on $R$ and 
define $\Omega_{R_n,\I_n}\isdef\ker(\Omega_{R_n}\ra\Omega_R).$\nl
Define recursively 
\be{eqnarray*}
\Omega(R,1)&\isdef&\Omega_R,\\
\Omega(R,n+1)&\isdef&\Omega(R,n)\oplus
\frac{R\oplus\Omega_R}{\span\{((n+1)a,\delta a)\mid a\in R\}}.
\e{eqnarray*}
Define 
\[\widetilde{\Omega}(R,n)\isdef\cases{
\Omega(R,n)\oplus
\frac{{\displaystyle R}}{{\displaystyle(n+1)R}}& if $n$ is odd\vspace{1mm}\cr
\Omega(R,n)& if $n$ is even.\cr}\]
Define
\[\widetilde{K_2}(R_n,\I_n)\isdef\cases{
K_2(R_n,\I_n)& if $n$ is odd \vspace{1mm}\cr
\frac{{\displaystyle K_2(R_n,\I_n)}}{{\displaystyle\span\{\denstein{aT^n,T}\mid a\in R\}}}& 
if $n$ is even.\cr}\]
\e{defi}
\be{lemma}\label{beromeg}
As $R$-modules 
$$\frac{\Omega_{R_n,\I_n}}{\delta \I_n}=\Omega(R,n)\oplus\frac{R}{(n+1)R}$$
\e{lemma}
\be{proof}
Write $J$ for the ideal of $R[T]$ generated by $T^{n+1}$.
We have
$$\Omega_{R_n}=\frac{\Omega_{R[T]}}{J\Omega_{R[T]}+\delta J}$$
and as $R$-modules
$$\Omega_{R[T]}=(R\te_\Z\Omega_{{\displaystyle \Z[T]}})\oplus(R[T]\te_R\Omega_R),$$
So
$$\Omega_{R_n,\I_n}=\underbrace{(R\oplus\Omega_R)\oplus\cdots
\oplus(R\oplus\Omega_R)}_{n \quad{\rm copies}}\oplus\frac{R}{(n+1)R}.$$
Dividing out $\delta \I_n$ yields the desired result.
\e{proof}
We are now in the position to apply the machinery of \cite{Clauwens;k}
to our situation.
As a matter of fact, the construction in {\em loc. cit.}
yields a homomorphism
\[\nu_n\colon K_2(R_n,\I_n)\ra\frac{\Omega_{R_n,\I_n}}{\delta \I_n},\]
even when $R$ possesses a $(n+1)\lambda$-ring structure.
In view of lemma~\ref{beromeg} we obtain a homomorphism
\[\nu_n\colon K_2(R_n,\I_n)\ra\Omega(R,n)\oplus\frac{R}{(n+1)R}.\] 
Furthermore we obtain a homomorphism 
\[\widetilde{\nu_n}\colon\widetilde{K_2}(R_n,\I_n)\ra
\widetilde{\Omega}(R,n)\] 
whenever $R$ is a $n\lambda$-ring ($n>1$).
\be{thm}\label{thmnuniso}
$\nu_n$ and $\widetilde{\nu_n}$ are isomorphisms.
\e{thm}
\be{proof}
We refer to {\em loc. cit.} for the definitions of the $\nu_n$.
We proceed by applying induction on $n$.\nl
$n=1$:
$R$ is a $2\lambda$-ring and
$\nu_1\colon K_2(R_1,\I_1)\ra\Omega_R\oplus\frac{R}{2R}$
is determined by
$$\nu_1\denstein{aT,b}=(a\delta b,[a^2\theta^2(b)]),\quad
\nu_1\denstein{cT,T}=(0,[c]).$$
It is straightforward to check that 
$\nu_1^{-1}\colon\Omega_R\oplus\frac{R}{2R}\ra K_2(R_1,\I_1)$
is well defined by
$$\nu_1^{-1}(a\delta b,[c])=\denstein{aT,b}+\denstein{a^2\theta^2(b)T,T}+
\denstein{cT,T}$$
$n>1$:
Consider the diagram 
$$\diagram{
&&&&K_2(R_n,\I_n^n)&&&&\cr
&&&&\mapdown{\tau}&&&&\cr
0&\ra&\frac{{\displaystyle R\oplus\Omega_R}}%
{{\displaystyle (na,\delta a)}}&\mapright{\iota}&
K_2(R_n,\I_n)&\mapright{\kappa}&
\Omega(R,n-1)\oplus\frac{{R}}{{(n+1)R}}&\ra&0\cr
&&\mapdown{\pi}&&\mapdown{\chi}&&\mapdown{\pi}&&\cr
0&\ra&\frac{{\displaystyle R}}{{\displaystyle nR}}&\mapright{\iota}&
K_2(R_{n-1},\I_{n-1})&\mapright{\kappa}&\Omega(R,n-1)&\ra&0\cr}$$
Here $\chi$ and $\tau$ are the obvious maps
and $\ker(\chi)=\im(\tau)$.
In the top row $\kappa$ is the obvious direct summand of $\nu_n$ and
$\iota([a,b\delta c])=\denstein{aT^{n-1},T}+\denstein{bT^n,c}$.
In the bottom row $\kappa$ is the obvious direct summand of $\nu_{n-1}$ and
$\iota([a])=\denstein{aT^{n-1},T}$.
The maps denoted by $\pi$ are the cononical projections.
We compute
\begin{eqnarray*}
 &&\nu_n\denstein{aT^n,b}=[0,a\delta b]\in
  \frac{R\oplus\Omega_R}{(na,\delta a)}\\
 &&\nu_n\denstein{aT^{n-1},T}=[a,0]\in
  \frac{R\oplus\Omega_R}{(na,\delta a)}\\
 &&\nu_n\denstein{aT^n,T}=[a]\in\frac{R}{(n+1)R}
\end{eqnarray*}
Therefore the map $\iota$ in the top row is split
by the remaining summand of $\nu_n$;
and since $\pi\nu_n=\nu_{n-1}\chi$ this implies that
the map $\iota$ in the bottom row is split
by the remaining summand of $\nu_{n-1}$.
The bottom row is exact by the induction hypothesis.\nl
Suppose that $x\in K_2(R_n,\I_n)$ and 
$$\kappa(x)=0\in\Omega(R,n-1)\oplus\frac{R}{(n+1)R}\quad (\star).$$
Then there exists a $y\in\frac{R\oplus\Omega_R}{(na,\delta a)}$
such that $\iota(\pi(y))=\chi(x).$
The exactness of the column guarantees the existence of an element
$z\in K_2(R_n,\I_n^n)$ satisfying $x-\iota(y)=\tau(z)$.
Thus there exists an $r\in R$ such that
$x+\denstein{rT^n,T}\in \im(\iota)$.
But $[r]=0\in\frac{R}{(n+1)R}$ because of $(\star)$.
So $x\in\im(\iota)$. This proves that $\nu_n$
is an isomorphism. 

If $n$ is odd and $n>1$, then the notions $n\lambda$-ring and 
$(n+1)\lambda$-ring coincide and the preceding proves that 
$\widetilde{\nu_n}$
is an isomorphism.
For $n$ even consider the following diagram.
$$\diagram{
0&\ra&{{\displaystyle R\oplus\Omega_R}\over{\displaystyle(na,\delta a)}}&
\mapright{\widetilde\iota}&
\widetilde{K_2}(R_n,\I_n)&\mapright{\widetilde\kappa}&\Omega(R,n-1)&\ra&0\cr
&&\mapdown{\pi}&&\mapdown{\widetilde\chi}&&\mapdown{ 1}&&\cr
0&\ra&\frac{{\displaystyle R}}{{\displaystyle nR}}&\mapright{\iota}&
K_2(R_{n-1},\I_{n-1})&\mapright{\kappa}&\Omega(R,n-1)&\ra&0\cr}$$
and proceed as before.
\e{proof}
\be{cor}{}
If $R$ possesses a structure of $n\lambda$-ring, 
then $H^1(K_2(R_n,\I_n);t)$ is isomorphic to
$H^1(\widetilde\Omega(R,n);\widetilde{\nu_n}t\widetilde{\nu_n}^{-1})$
\e{cor}
\be{proof}
Note that
$t\denstein{aT^n,T}=\denstein{\ol{a}T^n,T}$ \,in $K_2(R_n,\I_n)$
and $\denstein{aT^n,T}$ is an odd torsion element of $K_2(R_n,\I_n)$
if $n$ is even. But odd torsion elements vanish in $H^1(K_2(R_n,\I_n))$,
so $H^1(K_2(R_n,\I_n))\cong H^1(\widetilde{K_2}(R_n,\I_n))$.
In view of the preceding theorem this yields the desired result.
\e{proof}

This enables us to compute these cohomology groups in the cases where
$\widetilde{\nu_n}$ is manageable.
The next theorem for instance, shows what these groups look
like for  $n=1$ and $n=2$.
For all abelian groups $A$ and numbers $k$ we write ${}_kA$ to denote  
$\{a\in A\mid ka=0\}.$ 
\be{thm}
Let $R$ be a $2\lambda$-ring and $\conj\colon R\ra R$ the identity.Then
\begin{eqnarray*}
 H^1(K_2(R_1,\I_1))
 &\cong&\frac{R}{2R}\oplus{}_2(\Omega_R)\\
 H^1(K_2(R_2,\I_2))
 &\cong&\{\alpha\in{}_2(\Omega_R)\mid
 (1+\phi^2)\alpha\in\delta({}_2R)\}\\
 &&\oplus\frac{R}{2R}\oplus
 {\Omega_R\over2\Omega_R+\delta R+\im(1+\phi^2)}
\end{eqnarray*}
where $\phi^2\colon\Omega_R\ra\Omega_R$
is given by $\phi^2(a\delta b)=\psi^2(a)(b\delta b-\delta\theta^2(b)).$
\e{thm}
\be{proof}
Again we refer to \cite{Clauwens;k}, for more details on the operations
$\phi^2$.
According to proposition~\ref{propk2inv}
$$t(\denstein{aT,b})=\denstein{b,-aT}=\denstein{aT,b}$$ and
$$t(\denstein{aT,T})=\denstein{-T,aT}=\denstein{aT,T}$$ 
in $K_2(R_1,\I_1)$.
So in view of the corollary to theorem~\ref{thmnuniso} we have 
$$H^1(K_2(R_1,\I_1))\cong H^1(\frac{R}{2R}\oplus\Omega_R;1).$$
The isomorphism
$$\widetilde{\nu_2}\colon K_2(R_2,\I_2)\ra\Omega_R\oplus
{R\oplus\Omega_R\over(2a,\delta a)}$$ is given by
\be{eqnarray*}
\widetilde{\nu_2}(\denstein{aT,b})&=&(a\delta b,[a^2\theta^2(b),
(a^2-\theta^2(a))\delta\theta^2(b)+
\theta^2(a)b\delta b+\theta^2(b)a\delta a]),\\
\widetilde{\nu_2}(\denstein{aT,T})&=&(0,[a,0]),\\
\widetilde{\nu_2}(\denstein{aT^2,b})&=&(0,[0,a\delta b]).
\e{eqnarray*}
Using proposition~\ref{propk2inv} we compute 
$$\widetilde{\nu_2}t\widetilde{\nu_2}^{-1}(\alpha,[b,\gamma])=
(\alpha,[-b,-(1+\phi^2)(\alpha)-\gamma]).$$
Hence
$$\ker(1+\widetilde{\nu_2}t\widetilde{\nu_2}^{-1})=
\{(\alpha,[b,\gamma])\mid2\alpha=0\hbox{ and }
[0,(1+\phi^2)(\alpha)]=[0,0]\},$$
$$\im(1-\widetilde{\nu_2}t\widetilde{\nu_2}^{-1})=
\{(0,[2b,2\gamma+(1+\phi^2)(\alpha)])\}$$
and the quotient of these groups equals the right-hand-side
of the second isomorphism.
\e{proof}
As far as stability is concerned we have:
\be{prop}
Let $n\neq 0$ be even.
If $R$ is a $(n+2)\lambda$-ring and $\conj=1$, then
$$H^1(\ker(\widetilde{K_2}(R_{n+2},\I_{n+2})\ra
\widetilde{K_2}(R_n,\I_n)))\;\cong\;
{{}_{n+2}\ker(2\delta)\over{}_{n+2}\ker(\delta)}\oplus\frac{R}{2R}.$$
\e{prop}
\be{proof}
Consider the exact sequence
$$0\ra{R\oplus\Omega_R\over((n+1)a,\delta a)}\oplus
{R\oplus\Omega_R\over((n+2)a,\delta a)}
\stackrel{\widetilde\iota}{\longrightarrow}
\widetilde{K_2}(R_{n+2},\I_{n+2})\ra
\widetilde{K_2}(R_n,\I_n)\ra 0,$$
where $\widetilde\iota$ is defined by
$$\widetilde\iota([a,b\delta c],[x,y\delta z])=
\denstein{aT^n,T}+\denstein{bT^{n+1},c}
+\denstein{xT^{n+1},T}+\denstein{yT^{n+2},z}.$$
A splitting $\sigma$ of $\widetilde\iota$ is given by
the appropriate direct summand of $\widetilde\nu_{n+2}$.
The involution $t$ on both $\widetilde{K_2}$-groups 
induces the involution $\sigma t\widetilde\iota$ on 
$${R\oplus\Omega_R\over((n+1)a,\delta a)}\oplus
{R\oplus\Omega_R\over((n+2)a,\delta a)}.$$
A little computation shows that
$$\sigma t\widetilde\iota([a,\alpha],[b,\beta])=
([a,\alpha],[-b,\delta a-(n+1)\alpha-\beta]).$$
Now $([a,\alpha],[b,\beta])\in \ker(1+\sigma t\widetilde\iota)$,
if and only if $([2a,2\alpha],[0,\delta a-(n+1)\alpha])=0$.\nl
Thus putting $n=2m$,
there exist $r,s\in R$ satisfying the relations:
\be{eqnarray*}
2a&=&(2m+1)r,\\
2\alpha&=&\delta r,\\
(2m+2)s&=&0 \mbox{ \ and}\\
\delta s&=&\delta a-(2m+1)\alpha.
\e{eqnarray*}
Hence 
$[a,\alpha]=[a,\delta a-\delta s-m\delta r]%
=[a+(2m+1)mr,\delta (a-s)]=[0,-\delta s]=[-s,0]$
and $2\delta s=(2m+2)s=0$.\nl
Conversely, if $[a,\alpha]=[s,0]$ for some $s\in R$ satisfying 
$2\delta s=(2m+2)s=0$,
then $([a,\alpha],[b,\beta])\in \ker(1+\sigma t\widetilde\iota)$.\nl
The observation that 
$\im(1-\sigma t\widetilde\iota)=\{([0,0],[2b,\beta'])\}$
completes the proof.
\e{proof}
The final contribution to the comprehension of the value group of $\omega_2$
comes from the following proposition.
\be{prop} \label{propdber}
{\rm Compare \cite[theorem 4.1.]{Giffen;k2}}.
Let $(R,\conj,u)$ be a commutative ring with antistructure.
If $n$ is even,
$$d\colon H^1(K_1(R_n,\I_n))\lra H^1(K_2(R_n,\I_n))$$ assigns to the class
$[x]$ of the element $x\in 1+\I_n$ the class $[\{x,-u\}]$. 
Recall that we identified $H^1(K_1(R_n,\I_n))$ and $H^1(1+\I_n)$.
\e{prop}
\be{proof}
We will work in $\gl_{2k}(R_n)$ and $\st_{2k}(R_n)$.\nl
Suppose $x\in 1+\I_n$ and $\ol x=x^{-1}$.
Let $X$ be the image of $x$ under the map
$1+\I_n\lra \gl_1(R_n)\injarrow \gl_k(R_n)$.
By definition $t_{\conj,u_n}(X)X=\pmatrix{X&0\cr 0&X^{-1}}$ and
$h_{1\,k+1}(x)$ is a lift of this element in $\st_{2k}(R_n)$.
According to lemma~\ref{expld} 
$$d([x])=d([X])=
[h_{1\,k+1}^{-1}(x)\,t_{\conj,u_n}(h_{1\,k+1}(x))].$$
But from the definition of $t_{\conj,u_n}$ we compute
\be{eqnarray*}
t_{\conj,u_n}(h_{1\,k+1}(x))
&=&t_{\conj,u_n}(w_{1\,k+1}(x)w_{1\,k+1}(-1))\\
&=&w_{1\,k+1}(-u_n^{-1})w_{1\,k+1}(u_n^{-1}\ol x)\\
&=&w_{1\,k+1}(-u_n^{-1})w_{1\,k+1}(-1)
w_{1\,k+1}(1)w_{1\,k+1}(u_n^{-1}x^{-1})\\
&=&h_{1\,k+1}(-u_n^{-1})h_{1\,k+1}^{-1}(-u_n^{-1}x^{-1}).
\e{eqnarray*}
Thus
\be{eqnarray*}
d([x])&=&
[h_{1\,k+1}^{-1}(x)h_{1\,k+1}(-u_n^{-1})h_{1\,k+1}^{-1}(-u_n^{-1}x^{-1})]\\
&=&[\{x,u_n\}]\\
&=&[\{x,-u\}\{x,-(1+T)\}].
\e{eqnarray*}
It remains to show that $\steinberg{x,-(1+T)}$ vanishes 
in $H^1(K_2(R_n,\I_n))$.
First note that $\steinberg{x,-u}$ is a cycle:
\be{eqnarray*}
t(\steinberg{x,-u})&=&t(\denstein{-u^{-1}(1-x),-u})\\
&=&\denstein{-u^{-1},-u(1-x^{-1})}\\
&=&\steinberg{-u^{-1},x^{-1}}\\
&=&\steinberg{x,-u}^{-1}.
\e{eqnarray*}
Now choose $y\in R_n$ such that $1-x^{-1}=yT$.  
So $1-x=-\ol{y}T(1+T)^{-1}$.
We compute
\be{eqnarray*}
(1-t)(\denstein{T,y})&=&\denstein{T,y}\denstein{-T(1+T)^{-1},\ol y}\\
&=&\denstein{T,y}\denstein{T,-(1+T)^{-1}\ol y}
\denstein{-(1+T)^{-1},\ol{y}T}\\
&=&\denstein{T,y-(1+T)^{-1}\ol y+y\ol{y}T(1+T)^{-1}}\cdot\\
& &\denstein{-(1+T)^{-1},(1+T)(x-1)}\\
&=&\denstein{T,y-(1+T)^{-1}\ol y+y\ol{y}T(1+T)^{-1}}\steinberg{x,-(1+T)}.
\e{eqnarray*}
But since 
$$(y-(1+T)^{-1}\ol y+y\ol{y}T(1+T)^{-1})T=1-x^{-1}+1-x+(1-x^{-1})(x-1)=0,$$
we have
$$y-(1+T)^{-1}\ol y+y\ol{y}T(1+T)^{-1}=zT^n \mbox{ \ for some \ } z\in R.$$
For $\steinberg{x,-u}$ is a cycle, so is $\denstein{T,zT^n}$.
What's more $\denstein{T,zT^n}$ is an odd torsion element in $K_2(R_n,\I_n)$,
because $0=\denstein{T^{n+1},z}=(n+1)\denstein{T,zT^n}$ and $n$ is even.
This finishes the proof. 
\e{proof}
\be{cor}{}
If $\; u=-1$ in the situation of proposition~\ref{propdber},
$ d $ is the zero map.
\e{cor}
\be{punt}
The composition of homomorphisms
$$\Arf^s(R,1,-1)\injarrow
L_0^s(R,1,-1)\stackrel{\lambda\omega_1^s}{\longrightarrow}
C(R)=\frac{R}{\span\{x+x^2\mid x\in R\}}$$
maps $\plane{a,b}$ to $[ab]$. 
This surjection splits by the homomorphism $[r]\mapsto\plane{r,1}$.\nl
Writing $\widetilde{\Arf}(R)$ for the kernel,
we obtain a splitting
$$\Arf^s(R,1,-1)\cong \widetilde{\Arf}(R)\oplus 
\frac{R}{\span\{x+x^2\mid x\in R\}}.$$
$\widetilde{\Arf}(R)$ is 
generated by $\arfred{a,b}\isdef\plane{a,b}+\plane{ab,1}$, where $a,b\in R$.
The following relations hold in $\widetilde{\Arf}(R)$:\vspace{1mm}
 \halign{#&\quad#\hfil&\quad#\hfill\cr
 &$\arfred{a,b_1+b_2}=\arfred{a,b_1}+\arfred{a,b_2}$&\vspace{1mm}\cr
 &$\arfred{a,b}=\arfred{b,a}$&\vspace{1mm}\cr
 &$\arfred{a,b}=0$&for $a\in 2R$\vspace{1mm}\cr
 &$\arfred{ax^2,b}=\arfred{a,bx^2}$&for every $x\in R$\vspace{1mm}\cr
 &$\arfred{a,b}=\arfred{a,ab^2}$&\vspace{1mm}\cr
 &$\arfred{a,1}=0$&\cr}
\e{punt}
The secondary Arf invariant is by definition the 
the restriction of $\omega_2$ to the
$\widetilde{\Arf}$-part of $\ker(\omega_1^s)$:
$$\widetilde{\Arf}(R)\injarrow
\ker(\omega_1^s)\stackrel{\omega_2}{\lra}\coker(d)=H^1(K_2(R_2,\I_2)).$$
The next theorem tells us what this invariant looks like for $n=2$.
\be{thm} 
$\omega_2(\arfred{a,b})=[\denstein{aT^2,b}]\in H^1(K_2(R_2,\I_2))$.
\e{thm}
\be{proof}
Let $\arfred{a,b}=\plane{a,b}+\plane{ab,1}$ be represented by
$$\left[\pmatrix{a&0&1&0\cr 0&ab&0&1\cr 0&0&b&0\cr 0&0&0&-1\cr}\right]
-\left[\pmatrix{0&0&1&0\cr 0&0&0&1\cr 0&0&0&0\cr 0&0&0&0\cr}\right]\in 
L_0^s(R,1,-1).$$
A lift of this element in $L_0^s(R_2,\alpha,-(1+T))$ is given by
$$l\isdef\left[\pmatrix{a&0&1&0\cr 0&ab&0&1\cr 0&0&b&0\cr 0&0&0&-1\cr}\right]
-\left[\pmatrix{0&0&1&0\cr 0&0&0&1\cr 0&0&0&0\cr 0&0&0&0\cr}
\right].$$ 
To apply the map $G$ of definition~\ref{defomegaend} we choose
 $$\gamma\isdef x_{24}(T^2-T)x_{13}(b(T-T^2))
 h_{12}(1+abT^2)x_{31}(-aT)x_{42}(-abT)\in \st_4(R_2)$$
as a lift of
$$\left(\pmatrix{a&0&1&0\cr 0&ab&0&1\cr 0&0&b&0\cr 0&0&0&-1\cr}
+u_2\pmatrix{a&0&0&0\cr 0&ab&0&0\cr  1&0&b&0\cr 0&1&0&-1\cr}\right)
\pmatrix{0&0&u_2^{-1}&0\cr 0&0&0&u_2^{-1}\cr 1&0&0&0\cr 0&1&0&0\cr}$$
$$=\pmatrix{1&0&b(T-T^2)&0\cr 0&1&0&T^2-T\cr 
-aT&0&1&0\cr 0&-abT&0&1\cr}\in E_4(R_2).$$
Using the definition of $t$ and the calculations in the proof of 
proposition~\ref{propk2inv} we find
$$t\gamma^{-1}=x_{24}(T-T^2)x_{13}(b(T^2-T))h_{34}(1-abT^2)
x_{31}(aT)x_{42}(abT).$$
A little computation shows that
\begin{eqnarray*}
G(l)&=&[\gamma^{-1}(t\gamma)]\\
&=&[\denstein{abT,T-T^2}\denstein{aT,b(T^2-T)}\\
&&h_{12}(1+abT^2)h_{34}(1-abT^2)
h_{42}(1-abT^2)h_{31}(1+abT^2)]\\
&=&[\denstein{aT^2,b}]\in H^1(K_2(R_2,\I_2)).
\end{eqnarray*}
But since $\omega_2(\arfred{a,b})=G(l)$ this finishes the proof.
\e{proof}
Taking the (primary) Arf invariant into account we have the following result.
\be{thm}\label{thmtotw}
Let $R$ be a $2\lambda$-ring.
The invariant
$$\Arf^s(R,1,-1)\ra\frac{R}{\{x+x^2\}}\oplus
\frac{\Omega_R}{2\Omega_R+\delta R+\{x\delta y+x^2y\delta y\mid x,y\in R\}}$$
maps $\plane{a,b}$ to $([ab],[a\delta b])$.
\e{thm}
\be{proof}
We compute $\phi^2(a\delta b)$ modulo $2\Omega_R+\delta R$:
\be{eqnarray*}
\phi^2(a\delta b)&\equiv&\psi^2(a)(b\delta b-\delta\theta^2(b))\\
                 &\equiv&(a^2-2\theta^2(a))(b\delta b-\delta\theta^2(b))\\
                 &\equiv&a^2b\delta b-a^2\delta\theta^2(b)\\
                 &\equiv&a^2b\delta b.
\e{eqnarray*}
Thus
$$2\Omega_R+\delta R+\im(1+\phi^2)=
2\Omega_R+\delta R+\{x\delta y+x^2y\delta y\mid x,y\in R\}.$$
In view of the preceding the rest is obvious.
\e{proof} 
\noindent
Let $R$ be an arbitrary commutative ring.
We recognize $$\frac{\Omega_R}{2\Omega_R+\delta R}$$
as an instance of a cyclic homology group 
{\em viz.} $HC_1(R/2R)$.

The assignment $a\mapsto a\delta a$ 
determines a well-defined homomorphism
$$q'\colon R\ra \frac{\Omega_R}{\delta R}.$$
Under the assumption that $2R=0$ 
$$\theta\colon R\ra R \qquad x\mapsto x^2$$
$$\theta'\colon \frac{\Omega_R}{\delta R}\ra \coker\,q' \qquad
[a\delta b]\mapsto [a^2b\delta b]$$
are well-defined homomorphisms.
From this point of view
$$\frac{R}{\{x+x^2\}}=\coker(1+\theta)$$
and
$$\frac{\Omega_R}{2\Omega_R+\delta R+\{x\delta y+x^2y\delta y\mid x,y\in R\}}
=\coker(1+\theta').$$
We are a bit sloppy here in \vspace{1mm} denoting the projection 
$\frac{{\displaystyle\Omega_R}}{{\displaystyle\delta R}}\ra \coker\,q'$ by 1.

These observations are the motivation for investigating
(operations on) cyclic homology groups.
In the next chapter we will construct the homomorphism
\[\Arf^s(R,1,-1)\ra \frac{R}{\{x+x^2\mid x\in R\}}\oplus
{\Omega_R\over 2\Omega_R+\delta R+\{(r+r^2\delta s)\delta s\mid
r,s\in R\}}\]
without the assumption that $R$ carries some extra structure.

It turns out that the right generalization in the non-commutative case
involves  the notion of quaternionic homology groups.
We will enter into details in the next chapter.

\newpage
\section{Examples.} 
\setcounter{altel}{0}
\setcounter{equation}{0}

\be{nitel}{Example}
Let $R=\Z[X,Y]$ be the polynomial ring in two variables.
\be{thm}\label{thmzxy}
$$L_0^s(R,1,-1)\cong\frac{R}{\{f+f^2\}}\oplus
\frac{\Omega_R}{2\Omega_R+\delta R+
\{f\delta g+f^2g\delta g\mid f,g\in R\}}.$$
\e{thm}
\be{proof}
First we claim that $L_0^s(R,1,-1)=\Arf^s(R,1,-1)$.\nl
\noindent Let $(M,[\phi],e)\in BQ(R,1,-1)$ be given.
Then $b_{[\phi]}(m)(m)=0$ for every $m\in M$.
Choose a basis element $f$ in $M$.
There exists an element $g\in M$ such that $b_{[\phi]}(g)=f^*$.
Thus we obtain a decomposition 
$$(M,[\phi],e)\cong(N,[\phi_{\mid N}],[f,g])\perp
(N^\perp,[\phi_{\mid N^\perp}],h),$$
where $N\isdef\span(f,g)$,
$N^{\perp}\isdef\{m\in M\mid b_{[\phi]}(m)(N)=0\}$ 
and $h$ is some class of bases.
Given the fact that $K_1(R)\cong\Z/2$ it may be necessary
to interchange the roles of $f$ and $g$ to get the right
class of bases at the right hand side.
In this decomposition the first summand is isomorphic to
$$(R^2,\left[\left(\begin{array}{cc}a&1\\0&b\end{array}\right)\right],
[(1,0),(0,1)])$$ for some $a,b\in R$.
An induction argument proves the claim.\nl
Furthermore $R$ has a structure of $\lambda$-ring by lemma~\ref{lringext}.
Next we claim that 
$$\frac{R}{\{f+f^2\}}\oplus
\frac{\Omega_R}{2\Omega_R+\delta R+\{f\delta g+f^2g\delta g\mid f,g\in R\}}
\ra\Arf^s(R,1,-1)$$
defined by $$\left([x],\sum[a\delta b]\right)\longmapsto
\plane{x,1}+\sum\arfred{a,b}$$
is a well defined inverse of the homomorphism  
in theorem~\ref{thmtotw}.
The only non-trivial point on our checklist is: show that
this map respects the relation $$a\delta bc+ab\delta c+ac\delta b=0.$$
This amounts to showing that the relation
$$\arfred{ a,bc}=\arfred{ ab,c}+\arfred{ ac,b}$$
holds in $\widetilde{\Arf}(R)$.
But this follows immediately from the identity
$$\arfred{ f,g }=
\arfred{ f\frac{\partial g}{\partial x},x}+
\arfred{ f\frac{\partial g}{\partial y},y}
\quad\hbox{for every }f,g\in R.$$
It suffices to prove this for monomials by additivity.
By using the relations in $\widetilde{\Arf}(R)$ we see that
$$\arfred{X^iY^j,X^kY^l}=
\arfred{X^iY^jkX^{k-1}Y^l,X}+
\arfred{X^iY^jX^klY^{l-1},Y}$$
whenever $k$ or $l$ is even.
By symmetry this is also true when $i$ or $j$ is even.
In the remaining case $i$, $j$, $k$ and $l$ are all odd and
\be{eqnarray*}
\arfred{X^iY^j,X^kY^l}&=&\arfred{XY,X^{i+k-1}Y^{j+l-1}}\\
&=&\arfred{XY,XYX^{i+k-2}Y^{j+l-2}}\\
&=&\arfred{XY,X^{(i+k-2)/2}Y^{(j+l-2)/2}}.
\e{eqnarray*}
An induction argument finishes the proof.
\e{proof}
\e{nitel}
\be{nitel}{Example}
Let $G$ be the group with presentation
$$G\isdef\langle X,Y,S\mid S^2=(XS)^2=(YS)^2=1,\quad XY=YX\rangle.$$
We study $\Arf^s(G)$ and $\Arf^h(G)$.
Recall that we are working with the anti-involution determined by
$\ol{g}=g^{-1}$ for all $g\in G$.
Let $H$ be the subgroup of $G$ generated by $X$ and $Y$.
These groups fit into the split short exact sequence
$$1\lra H\lra G\lra C_2\lra 1,$$
where $C_2$ is the group of order two generated by $S$.
Elements of order two in $G$ have the form $X^iY^jS$ for some $i,j\in \Z$.
Every element $f\in\mbbf\,_2[G]$ can be decomposed in a unique way as 
$f=f_-+f_+S$
with $f_-,f_+\in\mbbf\,_2[H]$.
\be{prop}\label{propsuf1}
$\Arf^{s,h}(G)$ is generated by the elements
$$\cases{
\plane{1,1}&\vspace{0.4mm}\cr 
\plane{X^{2i}Y^{2j+1}S,S} & with $j\geq 0$\vspace{0.4mm}\cr
\plane{X^{2i+1}Y^{2j}S,S}&  with $i\geq 0$\vspace{0.4mm}\cr
\plane{X^{2i+1}Y^{2j+1}S,S}& with $i\geq 0$\vspace{0.4mm}\cr
\plane{X^{2i}Y^{2j+1}S,XS} & with $j\geq 0$\vspace{0.4mm}\cr
\plane{X^{2i+1}Y^{2j+1}S,XS}& with $j\geq 0$\vspace{0.4mm}\cr
\plane{X^{2i+1}Y^{2j+1}S,YS} & with $i\geq 0.$\cr
}$$
\e{prop}
\be{remark}\label{remsuf1}
We say that an element $f\in\mbbf\,_2[H]$ fulfils condition 1 resp. 2 if
all terms $X^iY^j$ of $f$ satisfy $i\geq 0$ resp. $j\geq 0$.
Using the fact that
for each $h\in \mbbf\,_2[H]$ there exist unique $h_0,h_1,h_2,h_3\in 
\mbbf\,_2[H]$ such that $$h=h_0^2+h_1^2x+h_2^2y+h_3^2xy,$$
we can reformulate
proposition~\ref{propsuf1} as follows.
Every element of $\Arf^{s,h}(G)$
is of the form
$$\plane{fS,S}+\plane{gS,XS}+\plane{hS,YS},$$
with
\be{enumerate}
\item[$\cdot$]
$f_1,f_3$ satisfy condition 1, $f_2$ satisfies 
condition 2 and $f_0\in\mbbf\,_2$
\item[$\cdot$]
$g_2,g_3$ satisfy condition 2 and $g_0=g_1=0$ 
\item[$\cdot$]
$h_3$ satisfies condition 1 and $h_0=h_1=h_2=0$.
\e{enumerate}
\e{remark}
\be{lemma}\label{lemsuf1}
Every element of $\Arf^{s,h}(G)$ is a sum of elements of the form
$$\plane{X^mY^nS,S},\quad
\plane{X^mY^nS,XS} ,\quad
\plane{X^mY^nS,YS} .$$
\e{lemma} 
\be{proof}
It suffices to prove this for generators
$\plane{X^iY^jS,X^kY^lS}$.
Conjugation by $X$ and $Y$ yields
$$\plane{X^iY^jS,X^kY^lS}=\cases{\plane{X^{i\pm2}Y^jS,X^{k\pm2}Y^lS}&\cr
                             \plane{X^{i}Y^{j\pm2}S,X^{k}Y^{l\pm2}S}.&\cr}$$
This proves that our generator has the desired form whenever one 
of the exponents $i$, $j$, $k$ or $l$ is even.\nl
If all exponents are odd, we have
$$\plane{X^iY^jS,X^kY^lS}=\plane{XYS,X^{k-i+1}Y^{l-j+1}S}$$ 
where both $k-i+1$
and $l-j+1$ are odd. But since
\be{eqnarray*}
\plane{XYS,X^{2i+1}Y^{2j+1}S}&=&
\plane{XYS,X^{i+1}Y^{j+1}SXYSX^{i+1}Y^{j+1}S}\\
&=&\plane{XYS,X^{i+1}Y^{j+1}S}
\e{eqnarray*}
and $$\plane{XYS,XYS}=\plane{XYS,1}=\plane{1,1}=\plane{S,S},$$
we can use an induction argument to prove the assertion in this case.
\e{proof}
We turn to the proof of the proposition.
\be{proof}
By lemma~\ref{lemsuf1} it suffices to prove the claim for the elements
$$\plane{X^mY^nS,S},\quad
\plane{X^mY^nS,XS} ,\quad
\plane{X^mY^nS,YS} .$$
\be{enumerate}
\item[$\diamond$]
$\plane{X^mY^nS,S}$\nl
We may assume that $m$ or $n$ is odd by using the relations
$$\plane{S,S}=\plane{1,1}$$
$$\plane{X^{2m}Y^{2n}S,S}=\plane{X^mY^nSSX^mY^nS,S}=\plane{X^mY^nS,S}.$$
Further we may assume that the odd exponent is positive since
$$\plane{X^{m}Y^{n}S,S}=\plane{SX^mY^nSS,S}=\plane{X^{-m}Y^{-n}S,S}.$$
\item[$\diamond$]
$\plane{X^mY^nS,XS}$\nl
We may assume that $n$ is odd by 
$$\plane{X^{2m}Y^{2n}S,XS}=
\plane{S,X^{-2m+1}Y^{-2n}S}=\plane{X^{2m-1}Y^{2n}S,S}$$
\be{eqnarray*}
\plane{X^{2m+1}Y^{2n}S,XS}&=&
\plane{X^{m+1}Y^{n}SXSX^{m+1}Y^{n}S,XS}\\
&=&\plane{X^{m+1}Y^{n}S,XS}.
\e{eqnarray*}
And we may assume that $n$ is positive since
$$\plane{X^{m}Y^{n}S,XS}=\plane{XSX^mY^nSXS,XS}=
\plane{X^{-m+2}Y^{-n}S,XS}.$$
\item[$\diamond$]
$\plane{X^mY^nS,YS}$\nl
We may assume that $n$ is odd by 
$$\plane{X^{2m}Y^{2n}S,YS}=
\plane{X^{2m}Y^{2n-1}S,S}$$
$$\plane{X^{2m+1}Y^{2n}S,YS}=
\plane{XS,X^{-2m}Y^{-2n+1}S}=\plane{X^{-2m}Y^{-2n+1}S,XS}.$$
We may assume that $m$ is odd by the relation
\be{eqnarray*}
\plane{X^{2m}Y^{2n+1}S,YS}&=&
\plane{X^{m}Y^{n+1}SYSX^{m}Y^{n+1}S,YS}\\
&=&\plane{X^{m}Y^{n+1}S,YS}.
\e{eqnarray*}
And we may assume that $m$ is positive since
$$\plane{X^{m}Y^{n}S,YS}=\plane{YSX^mY^nSYS,YS}=
\plane{X^{-m}Y^{-n+2}S,YS}.$$
\e{enumerate}
This completes the proof.
\e{proof}
The Arf invariant
$$\Arf^s(G)\lra
K(G)=\frac{\mbbf\,_2[G]}{\span\{a+\ol{a},b+b^2\mid a,b\in \mbbf\,_2[G]\}}$$
which maps
$$\cases{
\plane{X^iY^jS,X^kY^lS}\vspace{.5mm} & to \  $[X^{i-k}Y^{j-l}]$\cr
\plane{X^iY^jS,1}=\plane{1,1} & to \  $[1]$\cr}$$
splits by
$$\cases{
[X^iY^j]   &$\mapsto \plane{X^iY^jS,S}$ \cr
[X^iY^jS]  &$\mapsto \plane{1,1}.$ \cr}$$
We write $\widetilde{\Arf}(G)$ for the remaining summand.
Thus $$\Arf^s(G)\cong\widetilde{\Arf}(G)\oplus K(G).$$
Observe that the inclusion 
$\mbbf\,_2[H]\injarrow\mbbf\,_2[G]$ induces an isomorphism
$$K(H)\bijarrow K(G),$$
with inverse $[a]\longmapsto [a_-+a_+\ol{a_+}].$\nl
$\widetilde{\Arf}(G)$ is 
generated by 
$$\arfred{a,b}\isdef\plane{a_+S,b_+S}+\plane{a_+\ol{b_+}S,S},$$
where $a=\ol{a},b=\ol{b}$ in $\mbbf\,_2[G]$.
The following relations hold in $\widetilde{\Arf}(G)$:\vspace*{1mm}
 \halign{#&\quad#\hfil&\quad#\hfill\cr
 &$\arfred{a,b}=\arfred{a_+S,b_+S}$&\vspace{1mm}\cr
 &$\arfred{a,1}=\arfred{1,a}=0$&\vspace{1mm}\cr
 &$\arfred{a,S}=\arfred{S,a}=0$&\vspace{1mm}\cr
 &$\arfred{a,b_1+b_2}=\arfred{a,b_1}+\arfred{a,b_2}$&\vspace{1mm}\cr
 &$\arfred{a,b}=\arfred{b,a}$&\vspace{1mm}\cr
 &$\arfred{\ol{c}ac,b}=\arfred{a,cb\ol{c}}$&for every 
 $c\in \mbbf\,_2[G]$\vspace{1mm}\cr
 &$\arfred{a,b}=\arfred{a,\ol{b}ab}$&\vspace*{1mm}\cr}
\noindent 
Now we consider the representation $\rho\colon\mbbf\,_2[G]\ra M_2(R)$
of $G$ over the ring $R\isdef\mbbf\,_2[H]$
determined by
\be{eqnarray*}
X&\longmapsto&\pmatrix{X&0\cr0&X^{-1}\cr}\\
Y&\longmapsto&\pmatrix{Y&0\cr0&Y^{-1}\cr}\\
S&\longmapsto&\pmatrix{0&1\cr1&0\cr}
\e{eqnarray*}
and the diagram
$$\diagram{\Arf^s(G)&
{\buildrel \psi \over {\hbox to 100pt{\rightarrowfill}}}
&\Arf^s(R,1,1)\cr
\mapdown{\imath}&&\mapdown{\jmath}\cr
L^s(G)&\mapright{\widetilde{\rho}}\quad 
L_0^s(M_2(R),\alpha,1)\quad\mapright{\gamma}&
L_0^s(R,1,1).\cr}$$
Here $\imath$ and $\jmath$ are inclusion maps,\nl
$\widetilde{\rho}$ is induced by $\rho$,\nl
$U\isdef\pmatrix{0&1\cr1&0\cr}$,\nl
$\alpha(A)\isdef UA^tU$ for all $A\in M_2(R)$,\nl
$\gamma$ is the composition of the `scaling-isomorphism'
$$L_0^s(M_2(R),\alpha,1)\mapright{\cong}L_0^s(M_2(R),{\sf transpose},1)$$
and the `Morita-isomorphism'
$$L_0^s(M_2(R),{\sf transpose},1)\mapright{\cong}L_0^s(R,1,1).$$
\be{lemma}\label{lempsiar}
$\plane{X^iY^jS,X^kY^lS} \stackrel{\psi}{\longmapsto}
\plane{X^{-i}Y^{-j},X^kY^l}+\plane{X^iY^j,X^{-k}Y^{-l}}$.
\e{lemma}
\be{proof}
$\imath$ maps
$\plane{X^iY^j,X^kY^l}$
to
$$\left[\pmatrix{X^iY^jS&1\cr 0&X^kY^lS\cr}\right]-
\left[\pmatrix{0&1\cr 0&0\cr}\right],$$
$\widetilde{\rho}$ maps this element to
$$\left[\pmatrix{0&X^iY^j&1&0\cr X^{-i}Y^{-j}&0&0&1\cr
0&0&0&X^kY^l\cr0&0&X^{-k}Y^{-l}&0\cr}\right]-
\left[\pmatrix{0&0&1&0\cr0&0&0&1\cr0&0&0&0\cr0&0&0&0\cr}\right],$$
$\gamma$ maps this element to
$$\left[\pmatrix{X^iY^j&0&0&1\cr 0&X^{-i}Y^{-j}&1&0\cr
0&0&X^kY^l&0\cr 0&0&0&X^{-k}Y^{-l}\cr}\right]-
\left[\pmatrix{0&0&0&1\cr 0&0&1&0\cr 0&0&0&0\cr 0&0&0&0\cr}\right].$$
Now we apply the isometry $\pmatrix{U&0\cr 0&I\cr}$.
Note that this isometry is admissible since its 
class in $K_1(R)$ is trivial.
This yields
$$\left[\pmatrix{X^{-i}Y^{-j}&0&1&0\cr 0&X^{i}Y^{j}&0&1\cr
0&0&X^kY^l&0\cr 0&0&0&X^{-k}Y^{-l}\cr}\right]-
\left[\pmatrix{0&0&1&0\cr 0&0&0&1\cr 0&0&0&0\cr 0&0&0&0\cr}\right].$$
This element is equal to
$\jmath\left(\plane{X^{-i}Y^{-j},X^kY^l}+\plane{X^iY^j,X^{-k}Y^{-l}}\right).$
\e{proof}
Consequently
$$\psi(\arfred{fS,gS})=
\arfred{\ol{f},g}+\arfred{f,\ol{g}}\in\widetilde{\Arf}(R)$$
for all $f,g\in\mbbf\,_2[H]$.
We are now in the position to apply the machinery of the previous section
and in particular the secondary Arf invariant 
$$\widetilde{\Arf}(R)\ra
\frac{\Omega_R}{\delta R+\{(a+a^2b)\delta b\mid a,b\in R\}}$$
of theorem~\ref{thmtotw}.
\be{thm}
The invariant
\be{eqnarray*}
\Arf^s(G)&\lra &
\frac{R}{\span\{a+\ol{a},b+b^2\mid a,b\in R\}}\oplus
\frac{\Omega_R}{\delta R+\{(a+a^2b)\delta b\mid a,b\in R\}}\\
\plane{fS,gS}&\longmapsto&([f\ol{g}],[\ol{f}\delta g+f\delta\ol{g}]),
\e{eqnarray*}
is injective and the elements mentioned in proposition~\ref{propsuf1}
constitute a basis for $\Arf^s(G)$.
\e{thm}
\be{proof}
By the reformulation of proposition~\ref{propsuf1} in remark~\ref{remsuf1}
it suffices to prove that 
$$\plane{f'S,S}+\plane{gS,XS}+\plane{hS,YS}=0\quad \Longrightarrow\quad
 f'=g=h=0,$$
whenever $f',g,h\in\mbbf\,_2[H]$ satisfy the conditions mentioned in
remark~\ref{remsuf1}.
Suppose $\xi\isdef\plane{f'S,S}+\plane{gS,XS}+\plane{hS,YS}=0$.
Define $f\isdef f'+gX^{-1}+hY^{-1}$.
Then
$$\xi=\plane{fS,S}+\arfred{gS,XS}+\arfred{hS,YS}$$
and $f$ still fulfils the condition of
remark~\ref{remsuf1}.
The image 
$$([f],[(gX^{-1}+\ol{g}X)X^{-1}\delta X +(hY^{-1}+\ol{h}Y)Y^{-1}\delta Y])$$
of $\xi$ vanishes in
$$\frac{R}{\span\{a+\ol{a},b+b^2\mid a,b\in R\}}\oplus
\frac{\Omega_R}{\delta R+\{(a+a^2b)\delta b\mid a,b\in R\}}.$$
We will exploit the following facts to show that $f=g=h=0$.
\be{enumerate}
\item[$\cdot$]
For each $h\in R$ there are unique $h_0,h_1,h_2,h_3\in R$ such that\nl
$h=h_0^2+h_1^2X+h_2^2Y+h_3^2XY$.
\item[$\cdot$]
If $h\in R$ is symmetric, i.e. $\ol h=h$ and the constant term of $h$ 
is zero, then $h=p+\ol p$ for some $p\in R$.
\item[$\cdot$]
If $h\in R$ is symmetric,
then $h_0^2$, $h_1^2X$, $h_2^2Y$ and $h_3^2XY$ are symmetric.
\e{enumerate}
The fact that $[f]=0$ guarantees the existence of $a,b\in R$ such that
$$f=a+a^2+b+\ol{b}.$$
This implies: $f_0=0$ and $a_0^2+a^2$ is symmetric.
So $a_0+a=a_0+a_0^2+a_1^2X+a_2^2Y+a_3^2XY$ is symmetric as well.
By applying induction on 
$$\max\{|i|+|j|\, \mid X^iY^j \mbox{ is a term of } a+a^2\}$$ 
we conclude that $a+a^2$ is symmetric.
Hence $f_1^2X+f_2^2Y+f_3^2XY$ is symmetric, but the conditions on
$f_1,f_2,f_3$ make this impossible unless $f=0$.\nl
Since
$[(gX^{-1}+\ol{g}X)X^{-1}\delta X +(hY^{-1}+\ol{h}Y)Y^{-1}\delta Y]=0$
there exist $a,b,c\in R$ such that
$$(gX^{-1}+\ol{g}X)X^{-1}\delta X +(hY^{-1}+\ol{h}Y)Y^{-1}\delta Y
=(a+a^2)X^{-1}\delta X +(b+b^2)Y^{-1}\delta Y +\delta c$$
Since
\be{eqnarray*}
\delta c&=&\delta(c_0^2+c_1^2X+c_2^2Y+c_3^2XY)\\
        &=&c_1^2XX^{-1}\delta X +c_2^2YY^{-1}\delta Y 
           +c_3^2XYX^{-1}\delta X +c_3^2XYY^{-1}\delta Y ,
\e{eqnarray*}
we may assume that $c_0=0$ and
it follows that
$$gX^{-1}+\ol{g}X=a+a^2+c_1^2X+c_3^2XY,$$
$$hY^{-1}+\ol{h}Y=b+b^2+c_2^2Y+c_3^2XY.$$
Substituting $g=g_2^2Y+g_3^2XY$ and $h=h_3^2XY$ gives us the identities
$$g_2^2X^{-1}Y+g_3^2Y+\ol{g_2^2X^{-1}Y+g_3^2Y}=a+a^2+c_1^2X+c_3^2XY,$$
$$h_3^2X+\ol{h_3^2X}=b+b^2+c_2^2Y+c_3^2XY.$$
From these equations we deduce that
$a_0=a$ and $b_0=b$, thus $a+a^2=b+b^2=0$.
Hence $c_1=c_2=c_3=0$.
But then the restrictions on $g_2$, $g_3$ and $h_3$ 
imply $g_2=g_3=h_3=0$.
This finishes the proof.
\e{proof}
\e{nitel}

\newpage
{\Large {\bf \be{center}
Chapter III \vspace{4mm}\\
Hochschild, cyclic and quaternionic homology.
\e{center}}}
\vspace{6mm}

\setcounter{section}{0}
\section{Definitions and notations.}\label{defhomolo}
\setcounter{altel}{0}
\setcounter{equation}{0}

In the fourth section of the previous chapter we explained why we are 
interested in constructing certain operations on cyclic homology groups.
We start by summing up
the definitions of the various homologies we need.
We refer to \cite{LQ,Loday} for more details.

Let $k$ denote a commutative ring with identity.
\be{defi}
A simplicial $k$-module is a series of $k$-modules $\{M_n\mid n\in \N\},$ 
endowed with $k$-module homomorphisms
\[d_i\colon M_n\ra M_{n-1}\quad\mbox{  for all }\quad i\in\{0,1,\ldots,n\}\]
\[s_i\colon M_n\ra M_{n+1}\quad\mbox{  for all }\quad i\in\{0,1,\ldots,n\},\]
satisfying
\be{eqnarray*}
d_id_j&=&d_{j-1}d_i \quad\mbox{ if } i<j\\
d_is_j&=&\cases{
s_{j-1}d_i  &  if  $i<j$ \cr
1           &  if  $j\leq i\leq j+1$\cr
s_jd_{i-1}  &  if  $i>j+1$ \cr}\\
s_is_j&=&s_{j+1}s_i \quad\mbox{ if } i\leq j.
\e{eqnarray*}
\e{defi}
\be{defi}
A cyclic $k$-module is a simplicial $k$-module
$\{M_n \mid n\in\N\}$ equipped with homomorphisms
\[x\colon M_n\ra M_n\] satisfying
\begin{eqnarray*}
x^{n+1}&=&1\\
d_ix&=&-xd_{i-1}\quad\mbox{ for all }\quad i\in\{1,\ldots,n\}\\
d_0x&=&(-1)^nd_n\\
s_ix&=&-xs_{i-1}\quad\mbox{ for all }\quad i\in\{1,\ldots,n\}\\
s_0x&=&(-1)^{n+1}x^2s_n.
\end{eqnarray*}
\e{defi}
\be{defi}
A quaternionic $k$-module consists of a simplicial
$k$-module $\{M_n\mid n\in\N\}$ and homomorphisms 
\[\left\{\begin{array}{l}x\colon M_n\ra M_n\\
                         y\colon M_n\ra M_n\end{array}\right.\]
satisfying
\[\begin{array}{lcll}
x^{n+1}&=&y^2&\\
xyx&=&y&\\
d_ix&=&-xd_{i-1}&\mbox{ for all }\quad i\in\{1,\ldots,n\}\\
s_ix&=&-xs_{i-1}&\mbox{ for all }\quad i\in\{1,\ldots,n\}\\
d_iy&=&(-1)^nyd_{n-i}&\mbox{ for all }\quad i\in\{0,\ldots,n\}\\
s_iy&=&(-1)^{n+1}ys_{n-i}&\mbox{ for all }\quad i\in\{0,\ldots,n\}.
\end{array}\]
\e{defi}
\be{defi}
A quaternionic $k$-module is called a dihedral $k$-module when
$y^2=1$.
\e{defi}
\be{exa}\label{kanex1} 
Let $R$ be a $k$-algebra. 
We write $R^{n+1}$ as an abbreviation for the (n+1)-fold tensor product 
$R\te_kR\te_k\cdots\te_kR.$
The $k$-modules\[M_n\isdef R^{n+1}\]
and the homomorphisms $d_i$ and $s_i$ determined by
\be{eqnarray*}
d_i(a_0\te\cdots\te a_n)&\isdef&\cases{
a_0\te\cdots\te a_ia_{i+1}\te\cdots\te a_n&for $0\leq i<n$\cr
a_na_0\te a_1\te\cdots\te a_{n-1}&for $i=n$\cr}\\
s_i(a_0\te\cdots\te a_n)&\isdef&
a_0\te\cdots\te a_i\te1\te a_{i+1}\te\cdots\te a_n
\quad\mbox{ for all }0\leq i\leq n
\end{eqnarray*}
constitute a simplicial $k$-module.
The homomorphisms $x\colon R^{n+1}\ra R^{n+1}$ determined by 
\[x(a_0\te\cdots\te a_n)\isdef
(-1)^na_n\te a_0\te\cdots\te a_{n-1}\]
make this simplicial module into a cyclic module.\nl
If in addition $R$ is equipped with an anti-involution of $k$-algebras
$\conj\colon R\ra R$, it even becomes a dihedral module by defining 
\[y(a_0\te\cdots\te a_n)\isdef
(-1)^{\frac{1}{2}n(n+1)}
(\overline{a_0}\te\overline{a_n}\te\cdots\te\overline{a_1}).\]
\e{exa}
\be{exa}\label{kanex2}
More general, given a $k$-algebra $R$ and a $R$-bimodule $P$ 
we can turn
\[M_n\isdef P\te_kR^n\]
into a simplicial $k$-module through the homomorphisms
\begin{eqnarray*}
d_i(p\te r_1\te\cdots\te r_n)&\isdef&\cases{
pr_1\te r_2\te\cdots\te r_n& for $i=0$\cr
p\te r_1\te\cdots\te r_ir_{i+1}\te\cdots\te r_n& for $0<i<n$\cr
r_np\te r_1\te\cdots\te r_{n-1}& for $i=n$\cr}\\
s_i(p\te r_1\te\cdots\te r_n)&\isdef&
p\te r_1\te\cdots\te r_i\te1\te r_{i+1}\te\cdots\te r_n \\
&&\mbox{ for } 0\leq i\leq n
\end{eqnarray*}
\e{exa}
\be{defi}\label{defh}
For every simplicial $k$-module $M_*$ one constructs the chain complex
${\cal B}(M_*)$ called Hochschild complex as follows:
$$\diagram{\cdots\mapright{b}M_{n+1}\mapright{b}M_n\mapright{b}M_{n-1}
\mapright{b}\cdots\mapright{b}M_0}$$ 
where $$ b\isdef\sum_{i=0}^n(-1)^id_i.$$
The Hochschild-homology of $M_*$ is by definition the homology
of this chain complex.\nl
In case $M_*$ is the simplicial $k$-module of example~\ref{kanex1} 
we denote this chain complex by $(R^*,b)$ and its homology by $H_*(R).$
\e{defi}
\be{defi}\label{defhc}
If $M_*$ is a cyclic $k$-module one can build a double complex
${\cal C}(M_*)$:
$$\diagram{\vdots&&\vdots&&\vdots&&\vdots&&\cr
\da&&\da&&\da&&\da&&\cr
M_n&\mapleft{1-x}&M_n&\mapleft{L}&M_n&
\mapleft{1-x}&M_n&\mapleft{}&\cdots\cr
\mapdown{b}&&\mapdown{-b'}&&\mapdown{b}&&\mapdown{-b'}&&\cr
M_{n-1}&\mapleft{1-x}&M_{n-1}&\mapleft{L}&M_{n-1}&
\mapleft{1-x}&M_{n-1}&\lla&\cdots\cr
\da&&\da&&\da&&\da&&\cr
\vdots&&\vdots&&\vdots&&\vdots&& \cr}$$
where \begin{eqnarray*}b&\isdef&\sum_{i=0}^n(-1)^id_i\\
                       b'&\isdef&\sum_{i=0}^{n-1}(-1)^id_i\\
                       L&\isdef&\sum_{i=0}^nx^i \end{eqnarray*}
The cyclic homology $HC_n(M_*)$ of $M_*$
is by definition the n-th homology of the total
complex $\tot{\cal C}(M_*)$ associated to ${\cal C}(M_*)$, i.e.
\[HC_n(M_*)\isdef H_n(\tot{\cal C}(M_*)).\]
In the case that $M_*$ is the cyclic module 
of example~\ref{kanex1} we 
denote this cyclic homology by \[HC_n(R).\]
\e{defi}
\be{defi}\label{defhq}
If $M_*$ is a quaternionic module one can build a double complex 
${\cal D}(M_*)$ as follows:
\halign{\hfil$#$\hfil&\quad\hfil$#$\hfil&\quad\hfil$#$\hfil&\quad
\hfil$#$\hfil&\quad\hfil$#$\hfil&\quad\hfil$#$\hfil&\quad
\hfil$#$\hfil&\quad\hfil$#$\hfil&\quad\hfill$#$\hfil&\quad\hfil$#$\hfil\cr
\vdots&&\vdots&&\vdots&&\vdots&&\vdots&\cr
\da&&\da&&\da&&\da&&\da&\cr
M_n&\stackrel{\alpha}{\leftarrow}&M_n\oplus M_n&\stackrel{\beta}{\leftarrow}
&M_n\oplus M_n&
\stackrel{\gamma}{\leftarrow}&M_n&\stackrel{N}{\leftarrow}
&M_n&\leftarrow\cdots\cr
\mapdown{b}&&\mapdown{ -\widetilde{B}}&
&\mapdown{ \widehat{B}}&&\mapdown{-b'}&&\mapdown{b}&\cr
M_{n-1}&\stackrel{\alpha}{\leftarrow}&M_{n-1}
\oplus M_{n-1}&\stackrel{\beta}{\leftarrow}
&M_{n-1}\oplus M_{n-1}&
\stackrel{\gamma}{\leftarrow}&M_{n-1}&\stackrel{N}{\leftarrow}
&M_{n-1}&\leftarrow\cdots\cr
\da&&\da&&\da&&\da&&\da&\cr
\vdots&&\vdots&&\vdots&&\vdots&&\vdots&\cr}
where \begin{eqnarray*}b&\isdef&\sum_{i=0}^n(-1)^id_i\\
                       b'&\isdef&\sum_{i=0}^{n-1}(-1)^id_i\\
                       \widetilde{B}&\isdef& \pmatrix{b'&0\cr0&b\cr}\\
                       \widehat{B}&\isdef& \pmatrix{b&0\cr0&b'\cr}\\
                       L&\isdef&\sum_{i=0}^nx^i\\
                       N&\isdef&\sum_{i=0}^3Ly^i\\
                  \alpha&\isdef&\pmatrix{1-x&1-y\cr}\\
                  \beta&\isdef& \pmatrix{L&1+yx\cr-1-y&x-1\cr}\\
                  \gamma&\isdef&\pmatrix{1-x\cr yx-1\cr}
      \end{eqnarray*}
The quaternionic homology $HQ_n(M_*)$ of $M_*$
is by definition the n-th homology of the total
complex $\tot{\cal D}(M_*)$ associated to ${\cal D}(M_*)$ i.e.
\[HQ_n(M_*)\isdef H_n(\tot{\cal D}(M_*)).\]
In the case that $M_*$ is the quaternionic module of example~\ref{kanex1} we
denote this quaternionic homology by \[HQ_n(R).\]
\e{defi}

\newpage
\section{Reduced power operations.}\label{sechomoperaties}
\setcounter{altel}{0}
\setcounter{equation}{0}

In this section we will construct operations on various 
low dimensional homology groups. These operations will be used
later on to define new Arf invariants.
We feel that the material in this section is interesting 
in its own right.

\be{nota}
Let $p$ be a fixed prime number for the rest of this section.
For every $n\in \N$, $I_n$ denotes the set $\{1,2,\ldots,n\}$. 
$I_n$ will act as a set of indices.
The symmetric group of degree $p$, $S_p$ acts on the
$p$-fold cartesian product $I_n^p$ of $I_n$ by
$$\tau(i_1,\ldots,i_p)\isdef
(i_{\tau(1)},\ldots,i_{\tau(p)}) \mbox{ \ for all \ }(i_1,\ldots,i_n)\in 
I_n^p,\tau\in S_p.$$
Consider the permutation  $\sigma\isdef(1\,2\cdots p)^{-1}$.
Define $\Delta_n\isdef\{\gamma\in I_n^p\mid \sigma\gamma=\gamma\}$.
Let $\Gamma_n$ denote a set of representatives for the $\sigma$-orbits 
of the free action of $\sigma$ on $I_n^p-\Delta_n$.
\e{nota}

Let $R$ be an associative ring with identity.
Now recall the definitions of the Hochschild homology group $H_0(R)$ and
the cyclic homology group $HC_0(R)$.
Observe that both groups are equal to $\coker(b)$, where 
$b\colon R\te R\ra R$
is defined by $$b(r_1\te r_2)=r_1r_2-r_2r_1.$$
For all $r\in R$ we denote by $[r]$ the class of $r$ in $H_0(R)$.
\be{prop}\label{propthh0}
$\theta_p\colon H_0(R)\ra H_0(R/pR)$ defined by 
$$\theta_p([r])\isdef[r^p],$$
is a well-defined homomorphism.
\e{prop}
\be{proof}
For all maps $\alpha\colon I_n\ra R$ and elements
$\gamma=(i_1,\ldots,i_p)\in I_n^p$, we will write $\gamma(\alpha)$ instead of
$\alpha_{i_1}\alpha_{i_2}\cdots \alpha_{i_p}$.
We assert that
$$
\sum_{k=1}^p\sigma^k\gamma(\alpha)=p\gamma(\alpha)-
b\left(\sum_{l=1}^{p-1}\alpha_{i_1}\cdots \alpha_{i_l}\te 
\alpha_{i_{l+1}}\cdots \alpha_{i_p}\right). 
$$
This is easily verified by writing everything out.
For all $\alpha\colon I_2\ra R$, 
the following identity holds in $H_0(R/pR)$: 
\be{eqnarray*}
[(\alpha_{1}+\alpha_{2})^p]
&=&[\alpha_{1}^p+\alpha_{2}^p+
  \sum_{\gamma\in I_2^p-\Delta_2}\gamma(\alpha)]\\
&=&[\alpha_{1}^p+\alpha_{2}^p+
  \sum_{\gamma\in \Gamma_2}\sum_{k=\;1}^p\sigma^k\gamma(\alpha)]\\
&=&[\alpha_{1}^p+\alpha_{2}^p] \\ 
&=&[\alpha_{1}^p]+[\alpha_{2}^p].
\e{eqnarray*}
So it suffices to show that 
$[(b(\alpha_{1}\te \alpha_{2}))^p]=0$ in $H_0(R/pR)$.
Now then:
\be{eqnarray*}
[(b(\alpha_{1}\te \alpha_{2}))^p]&=&
[(\alpha_{1}\alpha_{2}-\alpha_{2}\alpha_{1})^p]\\
&=&[(\alpha_{1}\alpha_{2})^p+(-1)^p(\alpha_{2}\alpha_{1})^p]\\
&=&[\alpha_{1}\alpha_{2}(\alpha_{1}\alpha_{2})^{p-1}-
\alpha_{2}(\alpha_{1}\alpha_{2})^{p-1}\alpha_{1}]\\
&=&[b(\alpha_{1}\te \alpha_{2}(\alpha_{1}\alpha_{2})^{p-1})]\\
&=&0
\e{eqnarray*}
This proves the proposition.
\e{proof}

Recall the definitions of the Hochschild homology group $H_1(R)$ and
the cyclic homology group $HC_1(R)$: 
$$H_1(R)\isdef
\frac{\ker(b\colon R\te R\ra R)}{\im(b\colon R\te R\te R\ra R)}$$
$$HC_1(R)\isdef
\frac{\ker(b\colon R\te R\ra R)}{\im(b\colon R\te R\te R\ra R)+\im(1-x)}\,,$$
where 
$$b(r_1\te r_2)=r_1r_2-r_2r_1,$$
$$b(r_1\te r_2\te r_3)=r_1r_2\te r_3-r_1\te r_2r_3+r_3r_1\te r_2,$$
$$x\colon R\te R\ra R\te R
\mbox{ \ is defined by \ } x(r_1\te r_2)=-r_2\te r_1.$$
For all $\xi\in\ker(b\colon R\te R\ra R)$, we denote 
by $[\xi]$ the class of
$\xi$ in $H_1(R)$ as well as in $HC_1(R)$.\nl
Let $\alpha,\beta\colon I_n\ra R$ be set-theoretic maps.
For every $p$-tuple $\gamma=(i_1,\ldots,i_p)\in I_n^p$ we write 
$$\gamma(\alpha,\beta)$$ instead of 
$$\alpha_{i_1}\beta_{i_1}%
\alpha_{i_2}\beta_{i_2}\cdots\alpha_{i_{p-1}}\beta_{i_{p-1}}\te
\alpha_{i_p}\beta_{i_p}\in R\te R.$$
\be{thm}\label{thmoperaties}
The map $\theta_p\colon H_1(R)\ra HC_1(R/pR)$ determined by
\[
\left[\sum_{i\in I_n}\alpha_i\te\beta_i\right]\mapsto
\left[\sum_{i\in I_n}(\alpha_i\beta_i)^{p-1}\alpha_i\te\beta_i+
\sum_{\gamma\in\Gamma_n}\sum_{t=1}^{p-1}\left(t\sigma^t\gamma(\alpha,\beta)-
t\sigma^t\gamma(\beta,\alpha)\right)\right]
\]
is a well-defined homomorphism.
\e{thm} 
\be{remark}
In the case that $p=2$ this reads
$\theta_2\colon H_1(R)\ra HC_1(R/2R)$
$$\left[\sum_{i=1}^n\alpha_i\te\beta_i\right]\mapsto
\left[\sum_{i=1}^n\alpha_i\beta_i\alpha_i\te\beta_i+
\sum_{i<j}\left(\alpha_i\beta_i\te \alpha_j\beta_j+
\beta_i\alpha_i\te \beta_j\alpha_j\right)\right].$$
\e{remark}
We will prove this theorem with the help of a series of lemmas.
\be{lemma}\label{lemmacor}
Let $m>1$. For all $r_1,r_2,\ldots,r_m\in R$:
\be{eqnarray*}
\sum_{i=1}^mr_{i+1}r_{i+2}\cdots r_mr_1r_2\cdots r_{i-1}\te r_i&=&
(1-x)(r_1\cdots r_m\te1)\\
&+&b\left(\sum_{i=1}^{m-2}r_{i+2}\cdots r_m\te
r_1\cdots r_i\te r_{i+1}\right)\\
&+&b(1\te r_1\cdots r_{m-1}\te r_m)\\
&-&b(1\te r_1\cdots r_m\te1)
\e{eqnarray*}
\e{lemma}
\be{proof}
Simply a matter of writing everything out.
\e{proof}
\be{cor}{}
For all $\alpha,\beta\colon I_n\ra R$ and $\gamma\in I_n^p$
$$\left[\sum_{t=1}^{p-1}(t\sigma^t\gamma(\alpha,\beta)-
                          t\sigma^{t+1}\gamma(\alpha,\beta))\right]=
\left[\sum_{t=1}^p\sigma^t\gamma(\alpha,\beta)\right]=0.$$
\e{cor}
\be{cor}{}
$\theta_p$ does not depend on the choice of $\Gamma_n$.
\e{cor}
\be{lemma}\label{lemmatens}
Let $\F(R\times R)$ be the free abelian monoid on the set $R\times R$ and
$\te\colon \F(R\times R)\ra R\te R$ be the canonical morphism.
There is a bijective correspondence
between homomorphisms on $\ker(b\colon R\te R\ra R)$ and morphisms
on $\ker(b\te\colon \F(R\times R)\ra R)$ which kill all elements of the
form \vspace{1mm}\nl
$\begin{array}{ll}(u,0)&u\in R\\
                 (0,u)&u\in R\\
                 (u,v+w)+(u,-v)+(u,-w)&u,v,w\in R\\
                 (u+v,w)+(-u,w)+(-v,w)&u,v,w\in R.
\end{array}$
\e{lemma}
\be{proof}
To a homomorphism $f$ on $\ker(b)$, we associate the morphism $f\te$ on 
$\ker(b\te)$. It is clear that this morphism meets all requirements.

Conversely suppose $f$ is a morphism on $\ker(b\te)$ as in the statement 
above.
Define the homomorphism $g$ on $\ker(b)$ as follows:
If $\xi=\sum_i\alpha_i\te\beta_i$ belongs to 
$\ker(b)$, we choose
$\eta=\sum_i(\alpha_i,\beta_i)$ as a lift of $\xi$ in $\ker(b\te)$,
and define $g(\xi)\isdef f(\eta)$. 
Let us verify that this is well-defined.

Suppose $\tilde{\eta}=\sum_i(\tilde{\alpha}_i,\tilde{\beta}_i)$ 
is another lift of $\xi$ in $\ker(b\te)$.
Consider the difference $\eta-\tilde{\eta}$ in the free abelian group
$\F\!g(R\times R)$.
By definition of the tensor-product, this takes the form:
\be{eqnarray*}
&&\sum_{k_1}\left\{(u_{k_1}+v_{k_1},w_{k_1})-
(u_{k_1},w_{k_1})-(v_{k_1},w_{k_1})\right\}+\\
&&\sum_{k_2}\left\{(u_{k_2},w_{k_2})+(v_{k_2},w_{k_2})-
(u_{k_2}+v_{k_2},w_{k_2})\right\}+\\
&&\sum_{k_3}\left\{(u_{k_3},v_{k_3}+w_{k_3})-
(u_{k_3},v_{k_3})-(u_{k_3},w_{k_3})\right\}+\\
&&\sum_{k_4}\left\{(u_{k_4},v_{k_4})+(u_{k_4},w_{k_4})-
(u_{k_4},v_{k_4}+w_{k_4})\right\}
\e{eqnarray*}
for certain $u_{k_i},v_{k_i},w_{k_i}\in R.$
As a consequence we have in $\F(R\times R)$:
\be{eqnarray*}
\eta&+&\sum_{k_1}\left\{(u_{k_1},w_{k_1})+
       (-u_{k_1},w_{k_1})+(0,w_{k_1})\right\}+\\
    & &\sum_{k_1}\left\{(v_{k_1},w_{k_1})+
       (-v_{k_1},w_{k_1})+(0,w_{k_1})\right\}+\\
    & &\sum_{k_2}\left\{(u_{k_2}+v_{k_2},w_{k_2})+
       (-u_{k_2},w_{k_2})+(-v_{k_2},w_{k_2})\right\}+\\
    & &\sum_{k_2}\left\{2(0,w_{k_2})\right\}+\\
    & &\sum_{k_3}\left\{(u_{k_3},v_{k_3})+
       (u_{k_3},-v_{k_3})+(u_{k_3},0)\right\}+\\
    & &\sum_{k_3}\left\{(u_{k_3},w_{k_3})+(u_{k_3},-w_{k_3})+
       (u_{k_3},0)\right\}+\\
    & &\sum_{k_4}\left\{(u_{k_4},v_{k_4}+w_{k_4})+
       (u_{k_4},-v_{k_4})+(u_{k_4},-w_{k_4})\right\}+\\
    & &\sum_{k_4}\left\{2(u_{k_4},0)\right\}=\\
\tilde{\eta}&+&\sum_{k_1}\left\{(u_{k_1}+v_{k_1},w_{k_1})+
               (-u_{k_1},w_{k_1})+(-v_{k_1},w_{k_1})\right\}+\\
            & &\sum_{k_1}\left\{2(0,w_{k_1})\right\}+\\
            & &\sum_{k_2}\left\{(u_{k_2},w_{k_2})+
               (-u_{k_2},w_{k_2})+(0,w_{k_2})\right\}+\\
            & &\sum_{k_2}\left\{(v_{k_2},w_{k_2})+
               (-v_{k_2},w_{k_2})+(0,w_{k_2})\right\}+\\
            & &\sum_{k_3}\left\{(u_{k_3},v_{k_3}+w_{k_3})+
               (u_{k_3},-v_{k_3})+(u_{k_3},-w_{k_3})\right\}+\\
            & &\sum_{k_3}\left\{2(u_{k_3},0)\right\}+\\
            & &\sum_{k_4}\left\{(u_{k_4},v_{k_4})+
               (u_{k_4},-v_{k_4})+(u_{k_4},0)\right\}+\\
            & &\sum_{k_4}\left\{(u_{k_4},w_{k_4})+(u_{k_4},-w_{k_4})+
               (u_{k_4},0)\right\}.
\e{eqnarray*}
This implies $f(\eta)=f(\tilde{\eta})$. Hence $g$ is well-defined.
The rest is obvious.
\e{proof}
We want to apply this lemma to the map 
$$\tilde{\theta}_p\colon\ker(b\te)\ra HC_1(R/pR)$$
defined by
$$
\sum_{i\in I_n}(\alpha_i,\beta_i)\mapsto
\left[\sum_{i\in I_n}(\alpha_i\beta_i)^{p-1}\alpha_i\te\beta_i+
\sum_{\gamma\in\Gamma_n}\sum_{t=1}^{p-1}\left(t\sigma^t\gamma(\alpha,\beta)-
t\sigma^t\gamma(\beta,\alpha)\right)\right]
$$
But first we need another lemma to show that $\tilde{\theta}_p$  
is well-defined in the sense that the formula on the right-hand side
defines a cycle in $HC_1(R/pR)$.
\be{lemma}
For all $\alpha,\beta\colon I_n\ra R$ with 
$\sum_{i\in I_n}(\alpha_i,\beta_i)\in\ker(b\te)$
$$b\left(\sum_{i\in I_n}(\alpha_i\beta_i)^{p-1}\alpha_i\te\beta_i+
\sum_{\gamma\in\Gamma_n}\sum_{t=1}^p\left(t\sigma^t\gamma(\alpha,\beta)-
t\sigma^t\gamma(\beta,\alpha)\right)\right)=0.
$$
\e{lemma}
\be{proof}
Writing $\ol{\gamma}(\alpha,\beta)$ instead of 
$\alpha_{i_1}\beta_{i_1}\cdots \alpha_{i_p}\beta_{i_p}$,
for every $\gamma=(i_1,\ldots,i_p)\in I_n^p$, the expression becomes
\be{eqnarray*}
\lefteqn{\sum_{\gamma\in\Delta_n}(\ol{\gamma}(\alpha,\beta)-
               \ol{\gamma}(\beta,\alpha))+}\hspace{2ex}\\
& &\hspace*{-9ex}\sum_{\gamma\in\Gamma_n}\sum_{t=1}^{p-1}
(t\ol{\sigma^t\gamma}(\alpha,\beta)-
t\ol{\sigma^{t+1}\gamma}(\alpha,\beta)-
t\ol{\sigma^t\gamma}(\beta,\alpha)+
t\ol{\sigma^{t+1}\gamma}(\beta.\alpha))\\
&=&\sum_{\gamma\in\Delta_n}(\ol{\gamma}(\alpha,\beta)-\ol{\gamma}(\beta,\alpha))+
\sum_{\gamma\in\Gamma_n}\sum_{t=1}^p(\ol{\sigma^t\gamma}(\alpha,\beta)-
\ol{\sigma^t\gamma}(\beta,\alpha))\\
&=&\sum_{\gamma\in\Delta_n}(\ol{\gamma}(\alpha,\beta)-\ol{\gamma}(\beta,\alpha))+
\sum_{\gamma\in I_n^p-\Delta_n}(\ol{\gamma}(\alpha,\beta)-\ol{\gamma}(\beta,\alpha))\\
&=&\sum_{\gamma\in I_n^p}(\ol{\gamma}(\alpha,\beta)-\ol{\gamma}(\beta,\alpha))\\
&=&\left(\sum_{i\in I_n}\alpha_i\beta_i\right)^p-
   \left(\sum_{i\in I_n}\beta_i\alpha_i\right)^p\\
&=&0.
\e{eqnarray*}
This proves the assertion.
\e{proof}
We proceed by showing that $\tilde{\theta}_p$ is a morphism on 
$\ker(b\te)$.
\be{punt}
Suppose we are given $\alpha,\beta\colon I_n\ra R$ and 
$\alpha',\beta'\colon I_{n'}\ra R$, such that
$$\eta=\sum_{i\in I_n}(\alpha_i,\beta_i) \quad\mbox{ and }\quad
\eta'=\sum_{i\in I_{n'}}(\alpha_i',\beta_i')$$ 
are in $\ker(b\te)$.
Let's say 
$$r\isdef\sum_{i\in I_n}\alpha_i\beta_i=\sum_{i\in I_n}\beta_i\alpha_i
\quad\mbox{ and }\quad
r'\isdef\sum_{i\in I_{n'}}\alpha_i'\beta_i'=
\sum_{i\in I_{n'}}\beta_i'\alpha_i'.$$
We identify the disjoint union $I_n\vee I_{n'}$ and $I_{n+n'}$.
Define $\tilde{\alpha}\colon I_{n+n'}\ra R$ by
$$\tilde{\alpha}(i)\isdef\cases{\alpha(i)& if $i\in I_n$\cr
\alpha'(i)&if $i\in I_{n'}$\cr}$$
and define $\tilde{\beta}$ in the same way.
The map $I_{n+n'}\ra I_2$ defined by  
$$i\mapsto\cases{1& if $i\in I_n$\cr 2& if $i\in I_{n'}$\cr}$$
induces a map $\pi\colon I_{n+n'}^p\ra I_2^p$ which preserves
the $\sigma$-action.
Therefore
$$\Gamma_{n+n'}=\Gamma_n\cup \Gamma_{n'}\cup
\bigcup_{\lambda\in\Gamma_2}\pi^{-1}(\lambda).$$
Using this terminology we equate
\be{eqnarray*}
\lefteqn{\tilde{\theta}_p(\eta+\eta')-\tilde{\theta}_p(\eta)
-\tilde{\theta}_p(\eta')}\\
&=&\left[\sum_{t=1}^{p-1}\sum_{\lambda\in\Gamma_2}%
\sum_{\gamma\in\pi^{-1}(\lambda)}%
(t\sigma^t\gamma(\tilde{\alpha},\tilde{\beta})-
t\sigma^t\gamma(\tilde{\beta},\tilde{\alpha}))\right]\\
&=&\left[\sum_{t=1}^{p-1}\sum_{\lambda\in\Gamma_2}%
(t\sigma^t\lambda(\rho)-t\sigma^t\lambda(\rho))\right]\\
&=&0,
\e{eqnarray*}
where $\rho\colon I_2\ra R$ is defined by $\rho(1)=r$ and $\rho(2)=r'$.\nl
And $\lambda(\rho)=\rho_{i_1}\cdots\rho_{i_{p-1}}\te \rho_{i_p}$ if
$\lambda=(i_1,\ldots,i_p)\in I_2^p$.
\e{punt}
\be{punt}
Now it is time to apply lemma~\ref{lemmatens} and show that
$\tilde{\theta}_p$ induces a homomorphism $\theta_p'$ on
$\ker(b\colon R\te R\ra R):$
\be{itemize}
\item[$\diamond$]
It is clear that  
$\tilde{\theta}_p(u,0)=\tilde{\theta}_p(0,u)=0$, for all $u\in R$.
\item[$\diamond$]
$\tilde{\theta}_p((u,v+w)+(u,-v)+(u,-w))=0,$
for all $u,v,w\in R$:\nl
Define 
$\alpha,\beta\colon I_3\ra R$ by 
$$\be{array}{lll}
\alpha(1)\isdef u &\alpha(2)\isdef u& \alpha(3)\isdef u\\
\beta(1)\isdef v+w &\beta(2)\isdef -v& \beta(3)\isdef -w.
\e{array}$$
The map $I_3\ra I_2$ defined by $1\mapsto1,\;\;2\mapsto2,\;\;3\mapsto2$
induces a map $\pi\colon I_3^p\ra I_2^p$ which preserves the 
$\sigma$-action.
Define
$\alpha',\beta'\colon I_2\ra R$ by
$$\be{array}{ll}
\alpha'(1)\isdef u&\alpha'(2)\isdef u\\
\beta'(1)\isdef -v&\beta'(2)\isdef -w.
\e{array}$$
And finally we define 
$$\gamma_1\isdef u\beta_1'u\beta_2'\cdots u\te \beta_{i_p}'\qquad
\gamma_2\isdef \beta_1'u\beta_2'u\cdots \beta_{i_p}'\te u$$
for all $\gamma=(i_1,\ldots,i_p)\in I_2^p.$
$$\tilde{\theta}_p((u,v+w)+(u,-v)+(u,-w))=
\tilde{\theta}_p\left(\sum_{i\in I_3}(\alpha_i,\beta_i)\right).$$
\be{eqnarray*}
\lefteqn{(u(v+w))^{p-1}u\te(v+w)-(uv)^{p-1}u\te v-(uw)^{p-1}u\te w}
\hspace{10ex}\\
&=&-\sum_{\gamma\in I_2}\gamma_1+\sum_{\gamma\in\Delta_2}\gamma_1\\
&=&-\sum_{\gamma\in I_2-\Delta_2}\gamma_1
\phantom{xxxxxxxxxxxxxxxxxxxxxxxxx}\\
&=&-\sum_{\gamma\in\Gamma_2}\sum_{t=1}^p(\sigma^t\gamma)_1
\e{eqnarray*}
\be{eqnarray*}
\lefteqn{\sum_{\gamma\in\Gamma_3}\sum_{t=1}^{p-1}
(t\sigma^t\gamma(\alpha,\beta)-t\sigma^t\gamma(\beta,\alpha))}\hspace{10ex}\\
&=&\sum_{\lambda\in\Gamma_2}\sum_{t=1}^{p-1}
\sum_{\gamma\in\pi^{-1}(\lambda)}
(t\sigma^t\gamma(\alpha,\beta)-t\sigma^t\gamma(\beta,\alpha))\\
&&+\sum_{\gamma\in\Gamma_2}\sum_{t=1}^{p-1}
(t\sigma^t\gamma(\alpha',\beta')-t\sigma^t\gamma(\beta',\alpha'))
\e{eqnarray*}
But for all $\lambda\in \Gamma_2$ we have 
\be{eqnarray*}
\lefteqn{\left[\sum_{\gamma\in\pi^{-1}(\lambda)}
(t\sigma^t\gamma(\alpha,\beta)-t\sigma^t\gamma(\beta,\alpha))\right]}\\
&=&\pm\left[t((u(v+w))^{p-1}\te u(v+w)-
((v+w)u)^{p-1}\te (v+w)u)\right]\\
&=&0,
\e{eqnarray*}
since $[(ab)^{k}\te ab-(ba)^{k}\te ba]=0$ in $HC_1(R)$. Further
\be{eqnarray*}
\lefteqn{\left[\sum_{\gamma\in\Gamma_2}\sum_{t=1}^{p-1}
(t\sigma^t\gamma(\alpha',\beta')-t\sigma^t\gamma(\beta',\alpha'))\right]}\\
&=&\left[\sum_{\gamma\in\Gamma_2}\sum_{t=1}^{p-1}
(t(\sigma^{t+1}\gamma)_2-t(\sigma^t\gamma)_2)\right]\\
&=&\left[-\sum_{\gamma\in\Gamma_2}\sum_{t=1}^{p}
(\sigma^{t}\gamma)_2\right]
\e{eqnarray*}
Conclusion:
\be{eqnarray*}
\tilde{\theta}_p((u,v+w)+(u,-v)+(u,-w))&=&
\left[-\sum_{\gamma\in\Gamma_2}\sum_{t=1}^p
((\sigma^t\gamma)_1+(\sigma^t\gamma)_2)\right]\\
&=&\left[-\sum_{\gamma\in\Gamma_2}\sum_{t=1}^p
\sigma^t\gamma(\beta',\alpha')\right]\\
&=&0 
\e{eqnarray*}
according to the corollary following lemma~\ref{lemmacor}. 
\item[$\diamond$]
In a similar way one can prove that 
$\tilde{\theta}_p((u+v,w)+(-u,w)+(-v,w))=0,$
for all $u,v,w\in R$.
\e{itemize}
Thus we obtain a homomorphism $\theta_p'\colon\ker(b)\lra HC_1(R/pR)$.
\e{punt}
\be{prop}\label{prophc}
$\theta_p'(u\te v+v\te u)=[(uv)^{p-1}\te uv]\mbox{ for all }
u,v\in R.$
\e{prop}
\be{proof}
Define $\alpha,\beta\colon I_2\ra R$ by 
$\alpha(1)\isdef u,\;\;\alpha(2)\isdef v,$
$\beta(1)\isdef v,\;\;\beta(2)\isdef u.$
We equate
\be{eqnarray*}
\lefteqn{\theta_p'(u\te v+v\te u)}\\
&=&\left[(uv)^{p-1}u\te v+(vu)^{p-1}v\te u+\sum_{\gamma\in\Gamma_2}
\sum_{t=1}^{p-1}t\sigma^t\gamma(\alpha,\beta)-
t\sigma^t\gamma(\beta,\alpha)\right]
\e{eqnarray*}
The permutation $I_2\ra I_2$ determined by 
$1\mapsto2,\;\;2\mapsto1,$ induces a 
permutation $\pi\colon I_2^p\ra I_2^p$ which preserves the $\sigma$-action.
Since $\pi\gamma(\alpha,\beta)=\gamma(\beta,\alpha)$ for every 
$\gamma\in \Gamma_2,$ the term involving the double sum in the equation
above vanishes.\nl
Adding this to the fact that
$[(uv)^{p-1}u\te v+(vu)^{p-1}v\te u]=[(uv)^{p-1}\te uv]$
proves the proposition.
\e{proof}
\be{punt}
To finish the proof of theorem~\ref{thmoperaties} 
it only remains to show that
$$\theta_p'(uv\te w-u\te vw+wu\te v)=0 \mbox{ \ for all \ } u,v,w\in R.$$
For this purpose we define 
$\alpha,\beta\colon I_3\ra R$ by 
$$\be{array}{lll}
\alpha(1)\isdef uv&\alpha(2)\isdef vw&\alpha(3)\isdef wu\\
\beta(1)\isdef w&\beta(2)\isdef u&\beta(3)\isdef v
\e{array}$$
We use proposition~\ref{prophc} to equate
\be{eqnarray*}
\lefteqn{\theta_p'(uv\te w-u\te vw+wu\te v)}\\
&=&\theta_p'((uv\te w+vw\te u+wu\te v)-(u\te vw+vw\te u))\\
&=&[(uvw)^{p-1}uv\te w+(vwu)^{p-1}vw\te u+(wuv)^{p-1}wu\te v\\
& &+\sum_{\gamma\in \Gamma_3}\sum_{t=1}^{p-1}
(t\sigma^t\gamma(\alpha,\beta)-t\sigma^t\gamma(\beta,\alpha))-
(uvw)^{p-1}\te uvw].
\e{eqnarray*}
The permutation $I_3\ra I_3$ defined by 
$1\mapsto3,\;\;2\mapsto1,\;\;3\mapsto2,$ induces a permutation $\pi$
of $I_3^p$ which respects the $\sigma$-action.
Since $\pi\gamma(\alpha,\beta)=\gamma(\beta,\alpha)$ for every 
$\gamma\in \Gamma_3,$ the term involving the double sum in the equation
above vanishes.
And because
$$[(uvw)^{p-1}uv\te w+(vwu)^{p-1}vw\te u+(wuv)^{p-1}wu\te v
-(uvw)^{p-1}\te uvw]=0,$$
we are done.
\e{punt}
This completes the proof of theorem~\ref{thmoperaties}. 
\be{punt} 
Let $B\colon HC_0(R)\ra H_1(R)$ denote the homomorphism determined by 
$[r]\mapsto [r\te1+1\te r]=[1\te r]$.
The composition of $B$ and $\theta_p\colon H_1(R)\ra HC_1(R/pR)$ yields 
a homomorphism $q\colon HC_0(R)\ra HC_1(R/pR)$, which, 
as a consequence of proposition~\ref{prophc},
maps $[r]$ to $[r^{p-1}\te r].$
\e{punt}
\be{thm}
The homomorphism
$\theta_p\colon H_1(R)\ra HC_1(R/pR)$ induces a homomorphism 
$HC_1(R)\ra \coker(q).$
\e{thm}
\be{proof}
This is an immediate consequence of proposition~\ref{prophc}.
\e{proof}
\be{punt}\label{altdefhq}
Now recall the definition~\ref{defhq} 
of the quaternionic 
homology group $HQ_1(R)$.
There is an isomorphism
\[HQ_1(R)\bijarrow
{\ker((b\;\;1-y)\colon(R\te R)\oplus R\ra R)\over 
\span\left\{\begin{array}{c}(r\te s+s\te r,-rs-\ol{rs}),\\
(u\te v+\ol{u}\te\ol{v},vu-uv),\\
(0,2(w+\ol{w})),\\
(xy\te z-x\te yz+zx\te y,0)\e{array}\right\}}\]
defined by
$$[\varpi,a,b]\mapsto[\varpi,0,b+a+\ol{a}].$$
Here 
$$b\colon R\te R\ra R 
\mbox{ \ is defined by \ }
b(r_1\te r_2)=r_1r_2-r_2r_1,$$
$$y\colon R\ra R \mbox{ \ is defined by \ }y(r)=\ol{r}.$$
\e{punt}
\be{punt}\label{punteentensx}
The correspondence $x\mapsto[1\te x,0]$ obviously defines a
homomorphism $\nu_R\colon R\ra HQ_1(R)$.
\e{punt}
\be{lemma}\label{lemmaxtensx}
The map $x\mapsto[x\te x,0]$ 
determines a well-defined homomorphism $\mu_R\colon R\ra\coker\,\nu_R$.
\e{lemma}
\be{proof}
For all $x,y\in R$,\nl
$\mu_R(x+y)-\mu_R(x)-\mu_R(y)=[x\te y+y\te x,0]=[\nu_R(xy)]=0$.
The rest is obvious.
\e{proof}
\be{thm}\label{thmhqoperatie}
There exists a well-defined homomorphism
$$\vartheta=\vartheta_R\colon
HQ_1(R)\ra\coker\,(\mu_{R/2R})={HQ_1(R/2R)\over\span\{[x\te x,0]\mid
x\in R\}}$$
defined by:
\be{eqnarray*}
\left[
\sum_{i\in I_n}\alpha_i\te \beta_i, c\right]&\mapsto&
\left[\sum_{i\in I_n}\alpha_i\beta_i\alpha_i\te \beta_i
+\sum_{i\in I_n}\alpha_i\beta_i\te \beta_i\alpha_i\right.\\
& &\hspace*{.5ex}+\left.
\sum_{i<j}(\alpha_i\beta_i+\beta_i\alpha_i)\te
(\alpha_j\beta_j+\beta_j\alpha_j)+c\te\ol{c},c^2\right]
\end{eqnarray*}
\e{thm}
The proof will come about in a few steps.
\be{remark}
\be{eqnarray*}
\vartheta\left(\left[\sum_{i\in I_n}\alpha_i\te \beta_i, c\right]\right)
&=&\left[\sum_{i\in I_n}\alpha_i\beta_i\alpha_i\te \beta_i
+\left(\sum_{i\in I_n}\alpha_i\beta_i\right)\te
\left(\sum_{i\in I_n}\beta_i\alpha_i\right)\right.\\
& &+\left.\sum_{\gamma\in\Gamma_2}(\gamma(\alpha,\beta)-
\gamma(\beta,\alpha))+c\te\ol{c},c^2\right].
\e{eqnarray*}
$\sum_{\gamma\in\Gamma_2}[\gamma(\alpha,\beta)-
\sigma\gamma(\alpha,\beta)]=0$,
since $[r\te s+s\te r,0]=0$ in $\coker(\mu_{R/2R})$.
Thus it is clear that $\vartheta$ does not depend on the way the sum
$\sum_{i\in I_n}\alpha_i\te \beta_i$, is ordered.
\e{remark}
\be{lemma}
If $(\sum_{i\in I_n}\alpha_i\te \beta_i,c)\in\ker(b\;\;1-y)$, then
\be{eqnarray*}
&&\left(\sum_{i\in I_n}\alpha_i\beta_i\alpha_i\te \beta_i+
\alpha_i\beta_i\te \beta_i\alpha_i+\right.\\
&&\qquad
\left.\sum_{i<j}(\alpha_i\beta_i+\beta_i\alpha_i)\te(\alpha_j\beta_j+
\beta_j\alpha_j)+c\te\ol{c},c^2\right)
\e{eqnarray*}
is a cycle in $HQ_1(R/2R)$.
\e{lemma}
\be{proof}
The image of this expression under the homomorphism $(b\;\;1-y)$ equals
\be{eqnarray*}
&&\sum_{i\in I_n}(\alpha_i\beta_i)^2+(\beta_i\alpha_i)^2+
\alpha_i\beta_i\beta_i\alpha_i+\beta_i\alpha_i\alpha_i\beta_i+\\
&&\sum_{i<j}(\alpha_i\beta_i+\beta_i\alpha_i)
(\alpha_j\beta_j+\beta_j\alpha_j)+
(\alpha_j\beta_j+\beta_j\alpha_j)(\alpha_i\beta_i+\beta_i\alpha_i)+\\
&&c\ol{c}+\ol{c}c+c^2+\ol{c}^2=\\
&&\sum_{i,j\in I_n}(\alpha_i\beta_i+\beta_i\alpha_i)(\alpha_j\beta_j+
\beta_j\alpha_j)+c\ol{c}+\ol{c}c+c^2+\ol{c}^2=\\
&&(c+\ol{c})^2+c\ol{c}+\ol{c}c+c^2+\ol{c}^2=0
\e{eqnarray*}
\e{proof}
\be{lemma}\label{lemmaextens}
As before $\F(R\times R)$ denotes the free abelian monoid on the set 
$R\times R$ and
$\te\colon \F(R\times R)\ra R\te R$ is the canonical mapping.
Compare lemma~\ref{lemmatens}.
There is a bijective correspondence between homomorphisms on 
$$\ker((b\;\;1-y)\colon(R\te R)\oplus R\ra R)$$ 
and morphisms on
$$\ker\left((b\;\;1-y)\pmatrix{\te&0\cr0&1\cr}\right)=
\ker((b\!\te\;\;1-y)\colon \F(R\times R)\oplus R\ra R),$$ 
which kill all elements of the form\vspace{1mm}\nl
$\begin{array}{ll}((u,0),0)&u\in R\\
                   ((0,u),0)&u\in R\\
                   ((u,v+w)+(u,-v)+(u,-w),0)&u,v,w\in R\\
                   ((u+v,w)+(-u,w)+(-v,w),0)&u,v,w\in R.
\end{array}$
\e{lemma}
\be{proof}
Modulo a few minor adjustments the proof of lemma~\ref{lemmatens} will do.
\e{proof} 
We apply this lemma to the map
$\tilde{\vartheta}\colon\ker(b\!\te\,\;1-y)\ra \coker(\mu_{R/2R})$
defined by
\begin{eqnarray*}
\tilde{\vartheta}\left(\sum_{i\in I_n}(\alpha_i,\beta_i),c\right)&=&
\left[\sum_{i\in I_n}\alpha_i\beta_i\alpha_i\te \beta_i+
\alpha_i\beta_i\te \beta_i\alpha_i+\right.\\
& &\left.
\;\sum_{i<j}(\alpha_i\beta_i+\beta_i\alpha_i)\te
(\alpha_j\beta_j+\beta_j\alpha_j)+c\te\ol{c},c^2\right].
\end{eqnarray*}
\be{lemma}
$\tilde{\vartheta}$ is a morphism on $\ker(b\!\te\;\;1-y)$.
\e{lemma}
\be{proof}
If $$\eta\isdef\left(\sum_{i\in I_n}(\alpha_i,\beta_i),c\right)
\quad\mbox{ and }\quad 
\eta'\isdef\left(\sum_{i\in I_{n'}}(\alpha_i',\beta_i'),c'\right)$$ 
are in $\ker(b\!\te\;\;1-y)$, then
\be{eqnarray*}
\lefteqn{\tilde{\vartheta}(\eta+\eta')-\tilde{\vartheta}(\eta)-
\tilde{\vartheta}(\eta')}\hspace{2ex}\\
&=&\left[\left(\sum_{i\in I_n}\alpha_i\beta_i+\beta_i\alpha_i\right)\te
\left(\sum_{i\in I_{n'}}\alpha_i'\beta_i'+\beta_i'\alpha_i'\right)+\right.\\
& &\left.(c+c')\te(\ol{c}+\ol{c'})+c\te\ol{c}+c'\te\ol{c'},
(c+c')^2+c^2+{c'}^2\right]\\
&=&\left[(c+\ol{c})\te(c'+\ol{c'})+c\te\ol{c'}+c'\te\ol{c},cc'+c'c\right]\\
&=&\left[c\te c'+\ol{c}\te c'+c'\te\ol{c}+\ol{c}\te\ol{c'},cc'+c'c\right]\\
&=&0
\e{eqnarray*}
which proves the lemma.
\e{proof}
The next step is to prove that $\tilde{\vartheta}$ induces a homomorphism 
$\vartheta'$ on 
$\ker(b\;\;1-y)$. 
\be{punt}
We use lemma~\ref{lemmaextens}:
\be{enumerate}
\item[$\diamond$]
It is clear that 
$\tilde{\vartheta}((u,0),0)=\tilde{\vartheta}((0,u),0)=0$
for all $u\in R$.
\item[$\diamond$]
For all $u,v,w\in R$, 
\be{eqnarray*}
\lefteqn{\tilde{\vartheta}((u,v+w)+(u,-v)+(u,-w),0)}\hspace{4ex}\\
&=&[u(v+w)u\te(v+w)+uvu\te v+uwu\te w+\\
& &\;u(v+w)\te(v+w)u+uv\te vu+uw\te wu+\\
& &\;(u(v+w)+(v+w)u)\te(uv+vu)+\\
& &\;(u(v+w)+(v+w)u)\te(uw+wu)+\\
& &\;(uv+vu)\te(uw+wu),0]\\
&=&[uvu\te w+uwu\te v+uv\te wu+uw\te vu+\\
& &\;(uv+vu)\te(uw+wu)]\\
&=&0.
\e{eqnarray*}
\item[$\diamond$]
In the same fashion one proves that
$\tilde{\vartheta}((u+v,w)+(-u,w)+(-v,w),0)=0$
for all $u,v,w\in R$. 
\e{enumerate}
\e{punt}
Finally we use \ref{altdefhq} to verify that $\vartheta'$ induces 
the promised homomorphism $\vartheta$.
\be{punt}
For all $r,s,u,v,w,x,y,z\in R$ we have 
\be{eqnarray*}
\lefteqn{\vartheta'(r\te s+s\te r,rs+\ol{rs})}\hspace{5ex}\\
&=&[rsr\te s+rs\te sr+srs\te r+sr\te rs+\\
& &\;(rs+sr)\te(rs+sr)+(rs+\ol{rs})\te(rs+\ol{rs}),(rs+\ol{rs})^2]\\
&=&[0,rsrs+rs\ol{rs}+\ol{rs}rs+\ol{rs}\ol{rs}]\\
&=&0,\\
\lefteqn{\vartheta'(u\te v+\ol{u}\te\ol{v},vu-uv)}\hspace{5ex}\\
&=&[uvu\te v+\ol{u}\ol{v}\ol{u}\te\ol{v}+uv\te vu+
\ol{u}\ol{v}\te\ol{v}\ol{u}+\\
& &\;(uv+vu)\te(\ol{u}\ol{v}+\ol{v}\ol{u})+(uv+vu)\te(\ol{uv+vu}),
(uv+vu)^2]\\
&=&[uvu\te v+\ol{u}\ol{v}\ol{u}\te\ol{v}+uv\te vu+
\ol{u}\ol{v}\te\ol{v}\ol{u},(uv+vu)^2]\\
&=&0,\\
\lefteqn{\vartheta'(xy\te z+x\te yz+zx\te y,0)}\hspace{5ex}\\
&=&[xyzxy\te z+xyzx\te yz+zxyzx\te y+\\
& &\;xyz\te zxy+xyz\te yzx+zxy\te yzx+\\
& &\;(xyz+zxy)\te(xyz+yzx)+\\
& &\;(xyz+zxy)\te(zxy+yzx)+\\
& &\;(xyz+yzx)\te(zxy+yzx),0]\\
&=&0,\\
\lefteqn{\vartheta'(0,2(w+w'))=0.}\hspace*{5ex}
\e{eqnarray*}
\e{punt}
This step completes the proof of theorem~\ref{thmhqoperatie}.

\newpage
\section{Morita invariance.}
\setcounter{altel}{0}
\setcounter{equation}{0}

\be{thm}\label{thmmorita}
Let $A$ be the ring of $m\times m$-matrices over  the $k$-algebra $R$.
The trace-maps $\tr\colon A^n\ra R^n$ determined by
\[\tr(X_1\te X_2\te\cdots\te X_n)\isdef
\sum_{i_1,\ldots,i_n}\left(X_1\right)_{i_1i_2}
\te\left(X_2\right)_{i_2i_3}\te\cdots\te\left(X_n\right)_{i_ni_1}\]
yield a chain equivalence between the Hochschild complexes
$(A^*,b)$ and $(R^*,b)$.
A chain inverse is given by the maps $\iota\colon R^n\ra A^n$ defined by
\[\iota(r_1\te r_2\te\cdots\te r_n)\isdef E_{11}(r_1)\te
\cdots \te E_{11}(r_n)\]  
Where $E_{ij}(r)$
denotes the $m\times m$-matrix with $r$ in the
$(i,j)$-entry and zeros in all other entries.
\e{thm}
\be{proof}
It's easy to check that $\tr$ and $\iota$ are chain maps.
We immediately see that $\tr\comp\iota=1$.
We will show that $\iota\comp\tr\simeq1$ simply by giving a chain homotopy.
For that purpose we proceed to introduce the following definitions:\nl
Define
$$\gamma\colon A^{n+1}\ra A^{n+1}$$ by 
$$\gamma(X_0\te X_1\te\cdots\te X_n)\isdef
(-1)^{n+1}\sum_{i=1}^m E_{i1}(1)\te E_{1i}(1)X_nX_0\te X_1\te\cdots\te 
X_{n-1},$$ 
$$s\colon A^n\ra A^{n+1}$$
by
$$s(X_1\te\cdots\te X_n)\isdef X_1\te\cdots\te X_n\te 1$$
and finally
$$\chi_n\colon A^n\ra A^{n+1}$$
by
$$\chi_n\isdef (-1)^{n+1}\sum_{k=1}^n\gamma^ks.$$
The following relations are valid:
\setcounter{equation}{0}
\be{eqnarray}
\sum_{i=1}^mE_{i1}(1)E_{1i}(1)&=&1\\
d_0\gamma&\stackrel{1}{=}&(-1)^{n+1}d_n \nonumber\\
d_0\gamma^k&=&(-1)^{n+1}d_n\gamma^{k-1} \mbox{ \ if \ } k>0\\
d_i\gamma&=&-\gamma d_{i-1} \quad\mbox{ \ if \ } 1\leq i<n\\
d_i\gamma^k&=&=(-1)^k\gamma^kd_{i-k}\mbox{ \ if \ }k\leq i<n\\
d_1\gamma^2&\stackrel{3}{=}&-\gamma d_0\gamma \\
           &\stackrel{2}{=}&(-1)^n\gamma d_n \nonumber\\
           &=&(-1)^n\gamma d_{n-1} \nonumber\\
d_i\gamma^k&\stackrel{4}{=}&(-1)^{i-1}\gamma^{i-1}d_1\gamma^{k-i+1}\\
           &\stackrel{5}{=}&(-1)^{n+i-1}\gamma^id_{n-1}\gamma^{k-i-1} \nonumber\\
           &=&(-1)^{n+k}\gamma^{k-1}d_{n+i-k}\mbox{ \ if \ }0<i<k\nonumber\\
\gamma sd_n&=&\gamma d_ns\\
E_{1i}(1)XE_{j1}(1)&=&E_{11}(X_{ij})\\
d_n\gamma^ns&\stackrel{8}{=}&\iota\comp\tr
\e{eqnarray}
Now we are in the position to prove that
$$b\chi_n+\chi_{n-1}b=1-\iota\tr:$$
\be{eqnarray*}
b\chi_n
&=&(-1)^{n+1}\sum_{k=1}^n\sum_{i=0}^n(-1)^id_i\gamma^ks\\
&=&(-1)^{n+1}(\sum_{k=1}^n(d_0\gamma^ks+(-1)^nd_n\gamma^ks)+
   \sum_{k=1}^n\sum_{i=1}^{n-1}(-1)^id_i\gamma^ks)\\
&\stackrel{2}{=}&1-d_n\gamma^ns+\\
& &(-1)^{n+1}(\sum_{k=1}^{n-1}\sum_{i=k}^{n-1}(-1)^id_i\gamma^ks+
   \sum_{k=2}^n\sum_{i=1}^{k-1}(-1)^id_i\gamma^ks)\\
&\stackrel{4\,6}{=}&1-\iota\tr+\\
& &(-1)^{n+1}(\sum_{k=1}^{n-1}\sum_{i=k}^{n-1}(-1)^{i+k}\gamma^kd_{i-k}s+
   \sum_{k=2}^{n}\sum_{i=1}^{k-1}(-1)^{n+i+k}\gamma^{k-1}d_{n+i-k}s)\\
&=&1-\iota\tr+\\
& &(-1)^{n+1}(\sum_{k=1}^{n-1}\sum_{m=0}^{n-k-1}(-1)^m\gamma^kd_ms+
   \sum_{k=2}^{n}\sum_{m=n-k+1}^{n-1}(-1)^m\gamma^{k-1}d_ms)\\
&=&1-\iota\tr+
   (-1)^{n+1}(\sum_{k=1}^{n-1}\sum_{m=0}^{n-1}(-1)^m\gamma^kd_ms)\\
&\stackrel{7}{=}&1-\iota\tr+
   (-1)^{n+1}(\sum_{k=1}^{n-1}\sum_{m=0}^{n-1}(-1)^m\gamma^ksd_m)\\
&=&1-\iota\tr-\chi_{n-1}b.
\e{eqnarray*}
\e{proof}

Let $\conj\colon R\ra R$ be an anti-involution of $k$-algebras.
We extend this anti-involution to an anti-involution $\conj\colon A\ra A$
by defining $(\overline{X})_{ij}=\overline{X_{ji}}$ for every $X\in A$.
According to example~\ref{kanex1} 
we may regard
both $R^*$ and $A^*$ as quaternionic modules.
\be{thm}
The map $\tr$ induces isomorphisms
$$H_*(A)\mapright{\tr} H_*(R)$$
$$HC_*(A)\mapright{\tr} HC_*(R)$$
$$HQ_*(A)\mapright{\tr} HQ_*(R)$$
\e{thm}
\be{proof}
It is clear from the definitions that both $\iota$ and $\tr$ preserve
$x$ and $y$.
\e{proof}
\be{thm}\label{thmmorhq}
The following diagrams commute
\be{itemize}
\item[$\diamond$]
$$\diagram{HC_0(A)&\mapright{B}&H_1(A)\cr
\mapdown{\tr}&&\mapdown{\tr}\cr
HC_0(R)&\mapright{B}&H_1(R)\cr}$$
\item[$\diamond$]
$$\diagram{H_0(A)&\mapright{\theta_p}&H_0(A/pA)\cr
\mapdown{\tr}&&\mapdown{\tr}\cr
H_0(R)&\mapright{\theta_p}&H_0(R/pR)\cr}$$
\item[$\diamond$]
$$\diagram{
HC_1(A)&\mapright{\theta_p}&HC_1(A/pA)/\im(q)\cr
\mapdown{\tr}&&\mapdown{\tr}\cr
HC_1(R)&\mapright{\theta_p}&HC_1(R/pR)/\im(q)\cr}$$
\item[$\diamond$]
$$\diagram{
HQ_1(A)&\mapright{\vartheta_A}&\coker(\mu_{(A/2A)})\cr
\mapdown{\tr}&&\mapdown{\tr}\cr
HQ_1(R)&\mapright{\vartheta_R}&\coker(\mu_{(R/2R)})\cr}$$
\e{itemize}
\e{thm}
\be{proof}
A little examination of the definitions shows that
$B$,
$\theta_p$
and $\vartheta$ commute with
$\iota$.
\e{proof}

\newpage
\section{Generalized Arf invariants.}
\setcounter{altel}{0}
\setcounter{equation}{0}

Let $(R,\alpha,u)$ be a ring with anti-structure with $u=\pm 1$.
Thus $u$ is central and $\alpha$ is an anti-involution.

\be{thm}
The map
\[\Upsilon\colon\Arf^h(R,\alpha,u)\ra\coker(1+\vartheta_R)\]
determined by
\[\plane{a,b}\mapsto [a\te b,ab]\]
is a well-defined homomorphism.
\e{thm}
We are a bit sloppy here in denoting the projection 
$HQ_1(R)\ra\coker(\mu_{R/2R})$ by $1$.
\be{proof}
Recall the presentation of $\Arf^h(R,\alpha,u)$ from 
theorem~\ref{thmarfgr}.\nl
For all $a,b\in \Lambda_1(R)$ the element $(a\te b,ab)$
is a cycle in $HQ_1(R/2R)$: 
$$(b\;\;1-y)(a\te b,ab)=
ab+ba+ab+ba=0.$$
Next we will check that $\Upsilon$ respects all the relations of the 
aforementioned presentation.
\begin{enumerate}
\item
obvious
\item
obvious
\item
$[a\te b+b\te a,ab+ba]=0$
\item
$[a\te(x+\alpha(x)),a(x+\alpha(x)]=[a\te x+\alpha(x)\te a,ax+xa]=0$ 
\item
$[a\te \alpha(x) bx+xa\alpha(x)\te b,a\alpha(x) bx+xa\alpha(x) b]=$\nl
$[a\alpha(x) b\te x+xa\te \alpha(x) b,a\alpha(x) bx+xa\alpha(x) b]=$\nl
$[a\alpha(x) b\te x+bxa\te \alpha(x),a\alpha(x) bx+xa\alpha(x) b]=0$
\item
$\vartheta([a\te b,ab])=[aba\te b+ab\te ba+ab\te ba,abab]=[aba\te b,abab]$
\item
Suppose
$$\pmatrix{X&Y\cr Z&T\cr}\in \gl_{2n}(R)
\mbox{ \ satisfies \ }
t_{\alpha,u}\left(\pmatrix{X&Y\cr Z&T\cr}\right)=
\pmatrix{X&Y \cr Z&T \cr}^{-1}.$$
Then using the relations for $X$, $Y$, $Z$ and $T$, we equate
\be{eqnarray*}
\lefteqn{(1+\vartheta)
[X^\alpha\te T+Z\te Y^\alpha,X^\alpha T]}\hspace{2ex}\\
&=&[X^\alpha\te T+Z\te Y^\alpha+\\
& &X^\alpha TX^\alpha\te T+
ZY^\alpha Z\te Y^\alpha+X^\alpha T\te TX^\alpha+ZY^\alpha\te Y^\alpha Z+\\
& &(X^\alpha T+TX^\alpha)\te(ZY^\alpha+Y^\alpha Z)+X^\alpha T\te T^\alpha X,
X^\alpha T+(X^\alpha T)^2]\\
&=&[X^\alpha ZY^\alpha\te T+TX^\alpha Z\te Y^\alpha,X^\alpha ZY^\alpha T]\\
&=&[X^\alpha Z\te Y^\alpha T,X^\alpha ZY^\alpha T].
\e{eqnarray*}
Now theorem~\ref{thmmorhq} finishes the job.
\e{enumerate}
This finishes the proof.
\e{proof}
\be{remark}
In the case that $R$ is commutative and $\alpha$ is the identity we have
\be{eqnarray*}
\coker(1+\vartheta_R)&=&\frac{R}{\span\{x+x^2\mid x\in R\}}\oplus\\
& &{\Omega_R\over 2\Omega_R+\delta R+\{(r+r^2\delta s)\delta s\mid
r,s\in R\}}\\
&=&\coker(1+\theta_2\colon H_0(R)\ra H_0(R/2R))\oplus\\
& &\coker(1+\theta_2\colon H_1(R)\ra HC_1(R/2R)).
\e{eqnarray*}
This can be verified by a little examination of \ref{altdefhq} and
the definitions of $\theta_2$ and $\vartheta$ in
proposition~\ref{propthh0},
theorem~\ref{thmoperaties} and
theorem~\ref{thmhqoperatie}.
The projection of 
$$\Upsilon\colon\Arf^h(R,1,-1)\ra\coker(1+\theta_2)\oplus\coker(1+\theta_2)$$
on the first summand is just the old primary Arf invariant.
The secondary Arf invariant 
$$\Arf^s(R,1,-1)\lra {\Omega_R\over 2\Omega_R+\delta R+
\{(r+r^2\delta s)\delta s\mid r,s\in R\}}$$
factors through the projection of $\Upsilon$ on the second summand.
\e{remark}

\newpage
{\Large {\bf \be{center}
Chapter IV \vspace{4mm}\\
Applications to group rings.
\e{center}}}
\vspace{6mm}
\setcounter{section}{0}
\section{Quaternionic homology of group rings.}\label{secquahom}
\setcounter{altel}{0}
\setcounter{equation}{0}

The following exposition is based upon the work of J.-L. Loday in
\cite{Loday}.\nl
Let $k$ be a commutative ring with identity, $G$ a group and $k[G]$ the
group algebra of $G$ over $k$.
By providing $k[G]$ with the anti-involution $\conj$ determined by
$\ol{g}=g^{-1}$ for all $g\in G$, \hspace{1ex} 
$k[G]\te_k k[G]^n$ becomes a quaternionic module by means of 
example~\ref{kanex1} of chapter III.
\be{nota}
Denote by $\Gamma$ the set of conjugacy classes of $G$ and by $C\colon G\ra
\Gamma$ the map which assigns to $g\in G$ its conjugacy class $C(g)$.
Further we choose a section $S\colon\Gamma\ra G$ of $C$ such that
$S(C(g^{-1}))=(S(C(g))^{-1}$ for every $g\in G$ with $C(g)\neq C(g^{-1})$.
Finally,
for every set $V$ endowed with a right $G$-action 
we supply the free $k$-module $k[V]$ 
with a $k[G]$-bimodule structure by letting $G$ act trivially from the
left-hand side on $V$. 
\e{nota}
\be{punt}
For every $z\in G$, the right action $C(z)\times G\ra C(z)$ of $G$ on $C(z)$
defined by 
$(x,g)\mapsto g^{-1}xg$ for all
$x\in C(z)$ and $g\in G$,
makes $k[C(z)]$ into a $k[G]$-bimodule.
\e{punt}                           
\be{lemma} 
The map $$\phi\colon k[G]\te_k k[G]^n\ra\bigoplus_{z\in\im S}
k[C(z)]\te_k k[G]^n$$
determined by
$$\phi(g\te g_1\te\cdots\te g_n)\isdef
\cases{g_1\cdots g_ng\te g_1\te\cdots\te
g_n & if $gg_1\cdots g_n\in C(z)$\cr
0&otherwise,\cr}$$
is an isomorphism of simplicial modules with inverse determined by
$$h\te g_1\te\cdots\te g_n\mapsto(g_1\cdots g_n)^{-1}h\te g_1\te\cdots\te
g_n\mbox{  for all  }h\in C(z).$$
\e{lemma}
\be{proof}
See \cite{Loday}.
\e{proof}
\be{defi}
We say that $z\in G$ is of type \newline 
$\be{array}{lll}
\mbox{ \ \ }  &1&\mbox{   if  }\;z=\zm,\\
&2&\mbox{   if  }\;\zm\in C(z)\mbox{ \ and \ } z\neq \zm,\\
&3&\mbox{   if  }\;\zm\not\in C(z).
\e{array}$\nl
For each $i\in\{1,2,3\}$ let $S_i$ denote the subset of $\im S$
consisting of all elements of type $i$. 
Notice that $\im S$ is the disjoint union of the $S_i$.
Now $S$ was chosen in such a way that $z\in S_3\Leftrightarrow \zm\in 
S_3$, This allows us to write $S_3$ as a disjoint union of sets $S_3^+$
and $S_3^-$ such that $z\in S_3^+\Leftrightarrow \zm\in S_3^-.$
\e{defi}
\be{defi}\label{defxeny}
The simplicial module  $k[C(z)\cup C(\zm)]\te_k k[G]^n$
becomes a quaternionic module
by defining
$$x(g\te g_1\te\cdots\te g_n)\isdef
(-1)^n(g_1\cdots g_n)^{-1}gg_1\cdots g_n
\te(g_1\cdots g_n)^{-1}g\te g_1\te\cdots\te g_{n-1}$$
and 
$$y(g\te g_1\te\cdots\te g_n)\isdef
(-1)^{\frac{n(n+1)}{2}}(g_1\cdots g_n)^{-1}g^{-1}g_1\cdots g_n
\te g_n^{-1}\te\cdots\te g_1^{-1}$$
\e{defi}
\be{thm}\label{thmophak}
$$\phi\colon k[G]\te_k k[G]^n\ra\bigoplus_{z\in S_1\cup S_2\cup S_3^+}
k[C(z)\cup C(\zm)]\te_k k[G]^n$$
is an isomorphism of quaternionic modules.
\e{thm}
\be{proof} 
This is easy to check. The maps
$x$ and $y$ were defined so as to make $\phi$ respect the quaternionic
structure.
\e{proof}
\be{punt}\label{koleq}
For every group $G$ we define
$$d_i\colon k[G]^{n+1}\ra k[G]^n$$
by
$$\be{array}{rcl}
d_i(g_0\te g_1\te\cdots\te g_n)&\isdef&
g_0\te\cdots\te g_ig_{i+1}\te\cdots\te g_n \quad\mbox{ if  }\quad 
0\leq i<n\\
d_n(g_0\te g_1\te\cdots\te g_n)&\isdef&
             g_0\te\cdots\te g_{n-1},
\e{array}
$$
$$d\colon k[G]^{n+1}\ra k[G]^n \quad\mbox{ \ by \ }\quad
d\isdef\sum_{i=0}^n(-1)^id_i$$ 
and 
$$d'\colon k[G]^{n+1}\ra k[G]^n\quad\mbox{ \ by \ }\quad
d'\isdef\sum_{i=0}^{n-1}(-1)^id_i.$$ 
Now the map
$$s\colon k[G]^{n+1}\ra k[G]^{n+2}\mbox{ \ determined by \ }
s(g_0\te\cdots\te g_n)\isdef1\te g_0\te\cdots\te g_n$$
satisfies $sd+ds=sd'+d's=1$ and therefore
provides for a 
chain contraction of both the chain complexes
$$(k[G]^{*+1},d)\quad\mbox{ \ and \ }\quad (k[G]^{*+1},d').$$
Now let $G$ be a group and $H$ a subgroup of $G$.
Choose a set-theoretic section $\beta\colon H\bs G\ra G$, of the
canonical projection $\pi\colon G\ra H\bs G$,
satisfying $\beta(H)=1$ and define $\gamma\isdef\beta\comp\pi$.\nl
In what follows we will give homotopy-inverse maps of the inclusion-induced
maps
$$j_*\colon(k[H]^{*+1},d)\ra (k[G]^{*+1},d)$$
$$j_*'\colon(k[H]^{*+1},d')\ra (k[G]^{*+1},d')$$
and appropriate chain homotopies.\nl
The chain map $p_*$ determined by 
$$p_n\colon k[G]^{n+1}\ra k[H]^{n+1}$$ 
$$p_n(g_0\te\cdots\te g_n)\isdef$$
$$g_0\gamma(g_0)^{-1}\te\gamma(g_0)g_1\gamma(g_0g_1)^{-1}\te\cdots\te
\gamma(g_0g_1\cdots g_{n-1})g_n\gamma(g_0\cdots g_n)^{-1}$$
is a chain inverse to $j_*$, through the homotopies
$$h_n\colon k[H]^{n+1}\ra k[H]^{n+2} \quad\mbox{\  defined by \ }\quad
 h_n\isdef0$$ and
$$\ol{h}_n\colon k[G]^{n+1}\ra k[G]^{n+2}\quad\mbox{\  defined by \ }\quad
\ol{h}_n\isdef s(j_np_n-1).$$
Thus \be{eqnarray*}
p_nj_n-1&=&dh_n+h_{n-1}d\\
j_np_n-1&=&d\ol{h}_n+\ol{h}_{n-1}d.\e{eqnarray*}
Analogously we define
$$p_n'\colon k[G]^{n+1}\ra k[H]^{n+1}\quad 
\mbox{ \ by \ }\quad p_n'\isdef 0,$$
$$h_n'\colon k[H]^{n+1}\ra k[H]^{n+2}
\quad\mbox{ \ by \ }\quad h_n'\isdef -s$$
and
$$\ol{h'}_n\colon k[G]^{n+1}\ra k[G]^{n+2}\quad \mbox{ \ by \ }\quad
\ol{h'}_n\isdef -s.$$
Then again $p_n'$ determines a chain map and
\be{eqnarray*}
p_n'j_n'-1&=&d'h_n'+h_{n-1}'d'\\
j_n'p_n'-1&=&d'\ol{h'}_n+\ol{h'}_{n-1}d'.
\e{eqnarray*}
\e{punt}
\be{defi}\label{defgzstreep}
For all $z\in G$ one defines the subgroups $\gz$ and $\gzs$ of $G$ by
$$\gz\isdef\{g\in G\mid gz=zg\}$$
and
$$\gzs\isdef\left\{g\in G\mid g^{-1}zg\in\{z,\zm\}\right\}.$$
\e{defi}
Notice that
\be{itemize}
\item[$\cdot$]
$\gz=\gzm$.
\item[$\cdot$]
the correspondences $\gz\bs G\ra C(z)$ and 
$\gz\bs\gzs\ra\{z,\zm\}$ determined by
$\gz a\mapsto a^{-1}za$ are bijective.
\item[$\cdot$]
$\gzs$ acts from the right on $\{z,\zm\}$ by conjugation 
and this makes $k[z,\zm]$ into a $k[\gzs]$-bimodule.
\e{itemize}
\be{thm}\label{thmdc}
For all $z\in G$ the inclusion $\gzs\subseteq G$ induces a morphism
$$k[z,\zm]\te_kk[\gzs]^n\lra k[C(z)\cup C(\zm)]\te_kk[G]^n$$
$$a\te g_1\te\cdots\te g_n\mapsto a\te g_1\te\cdots\te g_n$$
of quaternionic modules.
\e{thm}
\be{proof}
We distinguish between three cases and keep \ref{koleq} and 
definition~\ref{defgzstreep} in mind.
\be{enumerate}
\item 
For all $z$ of type 1 we have $\gz=\gzs$ and the inclusion $\gz\subseteq G$
induces a morphism
of quaternionic modules
$$\diagram{
k[z]\te_kk[\gz]^n\cr
\isodown{}\cr
k\te_{k[\gz]}k[\gz]^{n+1}\cr
\mapdown{}\cr
k\te_{k[\gz]}k[G]^{n+1}\cr
\isodown{}\cr
k[C(z)]\te_kk[G]^n\cr}$$
mapping
$$z\te g_1\te\cdots\te g_n\mbox{ \ to \ }z\te g_1\te\cdots\te g_n.$$
Formulas for $x$ and $y$ can be found in definition~\ref{defxeny}.
\item 
For all $z$ of type 2 we have $C(z)=C(\zm)$ and the 
inclusion $\gzs\subseteq G$
induces a morphism of quaternionic modules
$$\diagram{
k[z,\zm]\te_kk[\gzs]^n\cr
\isodown{}\cr
k\te_{k[\gz]}k[\gzs]^{n+1}\cr
\mapdown{}\cr
k\te_{k[\gz]}k[G]^{n+1}\cr
\isodown{}\cr
k[C(z)]\te_kk[G]^n\cr}$$
mapping
$$a\te g_1\te\cdots\te g_n\mbox{ \ to \ }a\te g_1\te\cdots\te g_n.$$
Formulas for $x$ and $y$ can be found in definition~\ref{defxeny}.
\item 
For all $z$ of type 3 we have $\gz=\gzm=\gzs$ and the inclusion 
$\gz\subseteq G$
induces a morphism of quaternionic modules
$$\diagram{
k[z,\zm]\te_kk[\gz]^n\cr
\isodown{}\cr
(k\te_{k[\gz]}k[\gz]^{n+1})\oplus (k\te_{k[\gzm]}k[\gzm]^{n+1})\cr
\mapdown{}\cr
(k\te_{k[\gz]}k[G]^{n+1})\oplus (k\te_{k[\gzm]}k[G]^{n+1})\cr
\isodown{}\cr
k[C(z)\cup C(\zm)]\te_kk[G]^n\cr}$$
mapping
$$a\te g_1\te\cdots\te g_n\mbox{ \ to \ }a\te g_1\te\cdots\te g_n.$$
Formulas for $x$ and $y$ can be derived from definition~\ref{defxeny}.
\e{enumerate}
In all cases this morphism induces a chain map of the associated
quaternionic double complexes.
\e{proof}
By applying $k\te_{k[\gz]}-$ in the various situations of \ref{koleq}
that occur here, we see that these maps are 
chain equivalences on the columns by Shapiro's lemma.
Further \ref{koleq} enables us to compute explicit
chain inverses and chain homotopies.
To obtain the inverse homomorphism on the level of quaternionic homology
we use the following lemma.
\be{lemma}\label{lemmadubbelcomplexiso}
Suppose $j\colon\cee\ra\cees$ is a chain map of double complexes
$$\diagram{C_{20}& & & & &          &\ces_{20}&&&&\cr
\mapdown{d_{20}^v}&&\vdots&&&     &\mapdown{\des_{20}^v}&&\vdots&&\cr
C_{10}&\mapleft{d_{11}^h}&C_{11}&\cdots&&\mapright{j}
&\ces_{10}&\mapleft{\des_{11}^h}&\ces_{11}&\cdots&\cr    
\mapdown{d_{10}^v}&&\mapdown{d_{11}^v}&&&
&\mapdown{\des_{10}^v}&&\mapdown{\des_{11}^v}&&\cr
C_{00}&\mapleft{d_{01}^h}&C_{01}&\mapleft{d_{02}^h}&C_{02}&&
\ces_{00}&\mapleft{\des_{01}^h}&\ces_{01}&\mapleft{\des_{02}^h}&\ces_{02}\cr}
$$
which is a chain equivalence on the columns.
Let
$p_{*\,k}$ be a chain inverse of $j_{*\,k}$ and
\be{eqnarray*}
p_{m\,k}j_{m\,k}-1&=&d_{m+1\,k}^vh_{m\,k}+h_{m-1\,k}d_{m\,k}^v\\
j_{m\,k}p_{m\,k}-1&=&\des_{m+1\,k}^v\hes_{m\,k}+\hes_{m-1\,k}\des_{m\,k}^v.
\e{eqnarray*}
Then  $$\tau\colon H_1(\tot\cees)\ra H_1(\tot\cee)$$
defined by $$[a,b]\mapsto 
[p_{10}a+p_{10}\des_{11}^h\hes_{01}b+h_{00}d_{01}^hp_{01}b,p_{01}b]$$
for all $(a,b)\in\ker(\des_{10}^v\;\; \des_{01}^h)$,
is the inverse of $$j_*\colon H_1(\tot\cee)\ra H_1(\tot\cees).$$
\e{lemma}
\be{proof}
The map $j_*$ is an isomorphism
since $j$ is an equivalence on the columns.
By definition of double complex:
\be{eqnarray*}
d_{m-1\,k}^hd_{m\,k}^v+d_{m\,k-1}^vd_{m\,k}^h&=&0 
\quad\mbox{ \ for all \ }m,k\in\N\\
\des_{m-1\,k}^h\des_{m\,k}^v+\des_{m\,k-1}^v\des_{m\,k}^h&=&0 
\quad\mbox{ \ for all \ }m,k\in\N
\e{eqnarray*}
Now suppose $a\in\ces_{10}$ and $b\in\ces_{01}$ satisfy
$$\des_{10}^va+\des_{01}^hb=0.$$
Then 
\be{eqnarray*}
\lefteqn{d_{10}^v(p_{10}a+p_{10}\des_{11}^h\hes_{01}b+h_{00}d_{01}^hp_{01}b)+
d_{01}^hp_{01}b }\hspace{4em} \\
&=&-p_{00}\des_{01}^hb-p_{00}\des_{01}^h\des_{11}^v\hes_{01}b+
p_{00}j_{00}d_{01}^hp_{01}b\\
&=&-p_{00}\des_{01}^hj_{01}p_{01}b+p_{00}\des_{01}^hj_{01}p_{01}b\\
&=&0
\e{eqnarray*}
proves that $\tau([a,b])\in H_1(\tot\cee).$
Further we equate
\be{eqnarray*}
\lefteqn{j_*\tau([a,b])-[a,b]}\hspace{2em}\\
&=&[j_{10}p_{10}(a+\des_{11}^h\hes_{01}b)+
j_{10}h_{00}d_{01}^hp_{01}b-a,(j_{01}p_{01}-1)b]\\
&=&[(j_{10}p_{10}-1)(a+\des_{11}^h\hes_{01}b)+
j_{10}h_{00}d_{01}^hp_{01}b,0]\\
&=&[j_{10}p_{10}(j_{10}p_{10}-1)(a+\des_{11}^h\hes_{01}b)+
j_{10}p_{10}j_{10}h_{00}d_{01}^hp_{01}b,0].
\e{eqnarray*}
To obtain this last identity we used the fact that
$$(j_{10}p_{10}-1)(a+\des_{11}^h\hes_{01}b)+
j_{10}h_{00}d_{01}^hp_{01}b\in\ker(\des_{10}^v)$$
and
$$(j_{10}p_{10}-1)\ker(\des_{10}^v)\subseteq\im(\des_{20}^v).$$
To continue the computation we define
$c\isdef j_{10}p_{10}(j_{10}p_{10}-1)(a+\des_{11}^h\hes_{01}b).$
\be{eqnarray*}
\lefteqn{[c+j_{10}p_{10}j_{10}h_{00}d_{01}^hp_{01}b,0]}\\
&=&[c+j_{10}h_{00}p_{00}j_{00}d_{01}^hp_{01}b,0]\\
&=&[c+j_{10}h_{00}p_{00}\des_{01}^hj_{01}p_{01}b,0]\\
&=&[c+j_{10}h_{00}p_{00}\des_{01}^h(\des_{11}^v\hes_{01}+1)b,0]\\
&=&[c-j_{10}h_{00}p_{00}\des_{10}^va-
j_{10}h_{00}p_{00}\des_{10}^v\des_{11}^h\hes_{01}b,0]\\
&=&[c-j_{10}h_{00}d_{10}^vp_{10}(a+\des_{11}^h\hes_{01}b),0]\\
&=&[c-j_{10}p_{10}(j_{10}p_{10}-1)(a+\des_{11}^h\hes_{01}b),0]\\
&=&0.
\e{eqnarray*}
Thus we find $j_*\tau=1$ and since $j_*$ is already
an isomorphism, this proves the lemma.
\e{proof}
\be{punt}\label{puntranden}
We apply lemma~\ref{lemmadubbelcomplexiso} to the situation of
theorem~\ref{thmdc}:
Write $\D_1$ for the double complex 
$$\D\left(k[z,\zm]\te k[\gzs]^*\right),$$ 
$\D_2$ for the double complex
$$\D\left(k[C(z)\cup C(\zm)]\te k[G]^*\right)$$ 
and $$j\colon \D_1\ra\D_2$$ for the
chain map induced by the morphism of theorem~\ref{thmdc}. 
See definition~\ref{defhq} of chapter III  for the definition of $\D$.
Picture $\D_1:$ 
$$\diagram{
k[z,\zm]\te k[\gzs]^2&&&\cr
\mapdown{d_{20}^v}&&&\cr
k[z,\zm]\te k[\gzs]&\mapleft{d_{11}^h}&
(k[z,\zm]\te k[\gzs])\oplus&\vspace{-1.0ex}\cr
&&\hspace{1em}(k[z,\zm]\te k[\gzs])&\cr
\mapdown{d_{10}^v}&&\mapdown{d_{11}^v}&\cr
k[z,\zm]&\mapleft{d_{01}^h}&k[z,\zm]\oplus k[z,\zm]&\mapleft{d_{02}^h}
k[z,\zm]\oplus k[z,\zm]\cr}$$
and 
$\D_2:$ 
$$\diagram{
k[C(z)\cup C(\zm)]\te k[G]^2&&\cr
\mapdown{\des_{20}^v}&&\cr
k[C(z)\cup C(\zm)]\te k[G]&\mapleft{\des_{11}^h}&
(k[C(z)\cup C(\zm)]\te k[G])\oplus\vspace{-1.0ex}\cr
&&\hspace{1em}(k[C(z)\cup C(\zm)]\te k[G])\cr
\mapdown{\des_{10}^v}&&\mapdown{\des_{11}^v}\cr
k[C(z)\cup C(\zm)]&\mapleft{\des_{01}^h}&k[C(z)\cup C(\zm)]\oplus 
k[C(z)\cup C(\zm)]\mapleft{\des_{02}^h}\cr}$$
We use \ref{koleq}, definition~\ref{defxeny}, theorem~\ref{thmdc} and 
lemma~\ref{lemmadubbelcomplexiso} to obtain the following formulas.
\be{eqnarray*}
d_{10}^v,\des_{10}^v&\colon&a\te g\mapsto g^{-1}ag-a\\
d_{20}^v,\des_{20}^v&\colon&a\te g_1\te g_2\mapsto g_1^{-1}ag_1\te g_2-a\te
g_1g_2+a\te g_1\\
d_{11}^v,\des_{11}^v&\colon&(a\te g_1,0)\mapsto (-g_1^{-1}ag_1,0)\\
&&(0,b\te g_2)\mapsto (0,b-g_2^{-1}bg_2)\\
d_{01}^h, \des_{01}^h&\colon&(a,0)\mapsto 0\\
&&(0,b)\mapsto b-b^{-1}\\
d_{11}^h,\des_{11}^h&\colon&(a\te g_1,0)\mapsto a\te g_1+
g_1^{-1}ag_1\te g_1^{-1}a\\
&&(0,b\te g_2)\mapsto b\te g_2+g_2^{-1}b^{-1}g_2\te g_2^{-1}\\
d_{02}^h,\des_{02}^h&\colon&(a,0)\mapsto(a,-a-a^{-1})\\
&&(0,b)\mapsto(b+b^{-1},0)\\
p_{10}&\colon&g_1^{-1}ag_1\te g_2\mapsto 
\gamma(g_1)g_1^{-1}ag_1\gamma(g_1)^{-1}\te\gamma(g_1)g_2\gamma(g_1g_2)^{-1}\\
p_{01}&\colon&(g_1^{-1}ag_1,0)\mapsto 0\\
&&(0,g_2^{-1}bg_2)\mapsto(0,\gamma(g_2)g_2^{-1}bg_2\gamma(g_2)^{-1})\\
h_{00}=0&&\\
\hes_{01}&\colon&(g_1^{-1}ag_1,0)\mapsto(a\te g_1,0)\\
&&(0,g_2^{-1}bg_2)\mapsto(0,b\te g_2-b\te g_2\gamma(g_2)^{-1})
\e{eqnarray*}
\e{punt}
\be{thm}\label{thminvketen}
The inverse
$$\tau\colon H_1(\tot(\D_2))\lra H_1(\tot(\D_1))$$
of $j_*$ is \underline{determined} by
$(x,y)\longmapsto$
$$(\gamma(g_1)g_1^{-1}ag_1\gamma(g_1)^{-1}\te\gamma(g_1)g_2\gamma(g_1g_2)^{-1}
+b\te b,(0,\gamma(g_4)g_4^{-1}cg_4\gamma(g_4)^{-1})),$$
where
\be{eqnarray*}
x&=&g_1^{-1}ag_1\te g_2\in k[C(z)\cup C(\zm)]\te k[G]\\  
y&=&(g_3^{-1}bg_3,g_4^{-1}cg_4)\in k[C(z)\cup C(\zm)]\oplus
k[C(z)\cup C(\zm)].
\e{eqnarray*}
\e{thm}
\be{proof}
Under the given conditions we have
\be{eqnarray*}
\lefteqn{p_{10}\des_{11}^h\hes_{01}(g_3^{-1}bg_3,g_4^{-1}cg_4)}\\
&=&p_{10}\des_{11}^h(b\te g_3,c\te g_4-c\te g_4\gamma(g_4)^{-1})\\
&=&p_{10}(b\te g_3+g_3^{-1}bg_3\te g_3^{-1}b+c\te g_4+
g_4^{-1}c^{-1}g_4\te g_4^{-1}\\
&&-c\te g_4\gamma(g_4)^{-1}-\gamma(g_4)g_4^{-1}c^{-1}g_4\gamma(g_4)^{-1}\te
\gamma(g_4)g_4^{-1})\\
&=&b\te g_3\gamma(g_3)^{-1}+\gamma(g_3)g_3^{-1}bg_3\gamma(g_3)^{-1}\te
\gamma(g_3)g_3^{-1}b.
\e{eqnarray*} 
Applying the first relation of the list of theorem~\ref{thmrelaties} yields
$$[b\te g_3\gamma(g_3)^{-1}+\gamma(g_3)g_3^{-1}bg_3\gamma(g_3)^{-1}\te
\gamma(g_3)g_3^{-1}b]=[b\te b].$$
Using the formula for $\tau$ in lemma~\ref{lemmadubbelcomplexiso} 
yields the desired result.
\e{proof}
\be{thm}\label{thmrelaties}
For every $g,g_1,g_2\in\gzs$ and $a\in\{z,\zm\}$,
the following relations are valid in $H_1(\tot(\D_1))$.
\be{enumerate}
\item
$[g_1^{-1}ag_1\te g_2+a\te(g_1-g_1g_2),0,0]=0,$
\item
$[0,z+\zm,0]=0,$
\item
$[0,a,a+a^{-1}]=0$ and $[0,0,2(z+\zm)]=0,$
\item
$[z\te z,0,z+\zm]=0,$
\item
$[z\te g-\zm\te g,0,z-g^{-1}zg]=0,$
\item
$[z\te (g_1+g_2-g_1g_2),0,\epsilon(g_1,g_2)]=0,$ where
$$\epsilon(g_1,g_2)\isdef
\cases{z-\zm&if $g_1,g_2\not\in\gz$\cr 0&otherwise\cr}.$$
\e{enumerate}
\e{thm}
\be{proof}
\be{enumerate}
\item[{\em 1}]
follows immediately from the definition of $d_{20}^v.$
\item[{\em 2}]
is clear since $d_{02}^h(0,z)=(z+\zm,0).$
\item[{\em 3}]
$d_{02}^h(a,0)=(a,-a-a^{-1})$ and {\em 2} imply that
$[0,0,2(z+\zm)]=0$. The rest is obvious.
\item[{\em 4}]
Using the definitions of $d_{11}^h$ and $d_{11}^v$ we find
\be{eqnarray*}
0&=&[z\te g+g^{-1}zg\te g^{-1}z,-g^{-1}zg,0]\\
&=&[z\te g+z\te(z-g),0,z+\zm]  \mbox{ \ by {\em 1} and {\em 3} \ }\\ 
&=&[z\te z,0,z+\zm].
\e{eqnarray*}
\item[{\em 5}]
Using the definitions of $d_{11}^h$ and $d_{11}^v$ we equate
\be{eqnarray*}
0&=&[z\te g+g^{-1}\zm g\te g^{-1},0,z-g^{-1}zg]\\
&=&[z\te g+\zm\te(1-g),0,z-g^{-1}zg] \mbox{ \ by {\em 1} \ }\\ 
&=&[z\te g-\zm\te g,0,z-g^{-1}zg].
\e{eqnarray*}
Note that $[z\te1,0,0]=0$ by taking $g_1=g_2=1$ in {\em 1}.
\item[{\em 6}]
If $g_1\in\gz$, then {\em 6} follows from {\em 1}. \nl
If $g_1\not\in\gz$, then {\em 6} follows from {\em 1} and {\em 5}. 
\e{enumerate}
This completes the list of relations.
\e{proof}
\be{nota}
For every group $J$ we denote by 
$J_{{\rm ab}}$ the commutator quotient of $J$, 
i.e. $J_{{\rm ab}}=J/[J,J]$, and by $J_\#$ the quotient group 
$J_{{\rm ab}}/(J_{{\rm ab}})^2.$
\e{nota}
\be{thm}\label{thmdelenuitrek}
Let $k=\mbbf\:_2$.
\be{enumerate}
\item 
For every $z$ of type {\rm 1} the map
$$\eta\colon 
H_1(\tot(\D_1))\mapright{}\left((\gz)_\#/<z>\right)\times C_2$$
determined by
$$\left[\sum_i z\te g_i,n_1z,n_2z\right]\mapsto 
\left(\left[\prod_ig_i\right],t^{n_2}\right)$$
for all $n_1,n_2\in k$ and $g_i\in\gz$, where
$C_2$ denotes the cyclic group of order two generated by $t$,
is an isomorphism.
\item 
For every $z$ of type {\rm 2} the map
$$\eta\colon 
H_1(\tot(\D_1))\mapright{}\frac{(\gzs\times_{C_2}C_4)_\#}{<[z,t^2]>}$$
determined by
$$\left[\sum_ia_i\te g_i,\rho(n_1)z+\rho(n_2)\zm,
\rho(n_3)z+\rho(n_4)\zm\right]\mapsto\left[\prod_ig_i,t^n\right]$$
for all $n_1,n_2,n_3,n_4\in \Z$, $a_i\in\{z,\zm\}$ and $g_i\in\gzs$,
satisfying the cycle condition $\sum\rho(w(g_i))=\rho(n_3-n_4)$,
is an isomorphism. 
Here
$$\rho \mbox{ \ is the canonical map \ } \Z\ra k,$$
$$n\isdef \sum w(g_i)+2\left(n_1+n_2+n_4+\sum w'(a_i)w(g_i)\right),$$
$$w'(z)\isdef0,\,\,w'(\zm)\isdef1,$$
$$w(g)\isdef w'(g^{-1}zg) \mbox{ \ for all \ } g\in\gzs. $$
And $\gzs\times_{C_2}C_4$ is the pull-back of the diagram
$$\diagram{&&C_4\cr&&\mapdown{\pi_1}\cr\gzs&\mapright{\pi_2}&C_2\cr}$$
Here $C_4$ denotes the cyclic group of order four generated by $t$,
$\pi_1$ is the non-trivial map and $\pi_2(g)\isdef t^{w(g)}$ for 
every $g\in \gzs$.
\item
For every $z$ of type {\rm 3} the map
$$\eta\colon H_1(\tot(\D_1))\mapright{}(\gz)_\#$$
determined by
$$\left[\sum_i a_i\te g_i,n_1z+n_2\zm,n_3z+n_4\zm\right]\mapsto 
\left[\prod_ig_iz^{n_1+n_2+n_3}\right]$$
for all $n_1,n_2,n_3,n_4\in k$, $a_i\in\{z,\zm\}$ and $g_i\in\gz$,
satisfying the cycle condition $n_3=n_4$,
is an isomorphism.
\e{enumerate}
\e{thm}
\be{proof}
We will not enter into all the details of the proof; it is not difficult
but rather tedious.
\be{enumerate}
\item[{\em 1}] 
The data in ~\ref{puntranden} make it is easy to verify that the map
on $$\ker(d_{10}^v\;\;d_{01}^h)=(k[z]\te k[\gz])\oplus k[z]\oplus k[z]$$ 
determined by the expression in the definition of $\eta$ 
is a homomorphism which
vanishes on $\im(d_{02}^h)$, $\im(d_{20}^v)$ 
and $\im(d_{11}^h\;\;d_{11}^v).$\nl
Theorem~\ref{thmrelaties} enables us to check that
the inverse of $\eta$ is determined by 
$$([g],t^n)\mapsto [z\te g,0,\ol{n}z]$$
for every $g\in\gz$, $n\in\Z$.
\item[{\em 2}] 
Again $\eta$ is a well-defined homomorphism.
The inverse homomorphism is determined by 
$$[g,t^n]\mapsto [z\te g,0,\rho(\ent((n+1)/2))z+\rho(\ent(n/2))\zm]$$ 
for all
$g\in\gzs$ and $n\in\Z$ satisfying $\rho(w(g))=\rho(n)$.
\item[{\em 3}] 
The homomorphism $\eta^{-1}$ maps $[g]$ to $[z\te g,0,0]$ for all $g\in\gz$.
\e{enumerate}
Here $\ent$ denotes the entier function.
\e{proof}
\be{nota}
Write $$\Sigma(G)=
\bigoplus_{z\in S_1}\left(\left((\gz)_\#/<z>\right)\times C_2\right)
\oplus
\bigoplus_{z\in S_2}\frac{(\gzs\times_{C_2}C_4)_\#}{<[z,t^2]>}
\oplus
\bigoplus_{z\in S_3^+} (G_z)_\#$$
\e{nota}
\be{thm}\label{thmhqiso}
We have an isomorphism $\Psi\colon HQ_1(\mbbf\;_2[G])\lra\Sigma(G)$.
\e{thm}
\be{proof}
By theorem~\ref{thmophak}, theorem~\ref{thminvketen} and 
theorem~\ref{thmdelenuitrek}.
\e{proof}

\newpage
\section{Managing Coker$(1+\vartheta).$} 
\setcounter{altel}{0}
\setcounter{equation}{0}

Before we start with our reflections on 
$\coker(1+\vartheta_{{\displaystyle \mbbf\;_2[G]}})$
recall the theorems~\ref{thmophak}, \ref{thminvketen},  
\ref{thmdelenuitrek} and \ref{thmhqiso} of the previous section.
\be{lemma}
The isomorphism $\Psi$ induces an isomorphism
$$\Psi_1\colon\coker(\nu)\bijarrow\coker(\Psi\comp\nu)$$
and
$$\diagram{\coker(\Psi\comp\nu)
&=&\bigoplus_{z\in S_1}\left(\left((\gz)_\#/<z>\right)\times C_2\right)
\oplus\cr
&&\bigoplus_{z\in S_2}\left(\gzs\right)_\#/<z>
\oplus\hfill\cr
&&\bigoplus_{z\in S_3^+} (G_z)_\#/<z>\hfill\cr}
$$
\e{lemma}
\be{proof}
To determine $\coker(\Psi\comp\nu)$ 
we compute $$\Psi(\nu(x))=\Psi([1\te x,0,0])$$ for $x\in\mbbf\;_2[G]$.
We may assume that $x\in G$.
There exist $z\in\im\,S$ and $g\in G$ such that $x=g^{-1}zg.$
Now 
$$\phi(1\te g^{-1}zg)=g^{-1}zg\te g^{-1}zg,$$
\be{eqnarray*}
\tau([g^{-1}zg\te g^{-1}zg,0,0])&=&
[\gamma(g)g^{-1}zg\gamma(g)^{-1}\te
\gamma(g)g^{-1}zg\gamma(zg)^{-1},0,0]\vspace{.5mm}\\
&=&[\gamma(g)g^{-1}zg\gamma(g)^{-1}\te\gamma(g)g^{-1}zg\gamma(g)^{-1},0,0]
\vspace{.5mm}\\
&=&\cases{[z\te z,0,0]& if $z$ is of \vspace{.5mm}type 1\cr
          [z^{\pm1}\te z^{\pm1},0,0]& if $z$ is of \vspace{.5mm}type 2\cr
          [z\te z,0,0]& if $z$ is of type 3\cr}
\e{eqnarray*}
applying the isomorphism $\eta$ we find
$$\Psi([1\te g^{-1}zg,0,0])=
\cases{([z],1)=1\in((\gz)_\#/<z>)\times C_2
&if $z$ is \vspace{1mm} of type 1\cr
[z,1]\in(\gzs\times_{C_2}C_4)_\#/<[z,t^2]>
&if $z$ \vspace{1mm}is of type 2\cr
[z]\in(\gz)_\#&if $z$ is of type 3\cr}$$
The rest is clear now.
\e{proof}
\be{defi}
Let $G$ be a group. Define \[\tilde{F}(z)\isdef\cases{
{\displaystyle {(G_z)_\#\over 
<\{x\in G\mid x=z \,\vee\, x^2=z\}>}}\times C_2,& 
if $z$ is of\vspace{1mm} type $1$ \cr
{\displaystyle {(\ol{G_z})_\# \over
<\{x\in G\mid x=z\,\vee\, x^2=z\}>}}, & if $z$ is of \vspace{1mm}type $2$\cr 
{\displaystyle {(G_z)_\#\over<\{x\in G\mid x=z\,\vee\, x^2=z\}>}},
&if $z$ is of type $3$.\cr
}\]
\e{defi}
\be{lemma}
The isomorphism $\Psi_1$ induces an isomorphism 
$$\Psi_2\colon\coker(\mu)\bijarrow\coker(\Psi_1\comp\mu)$$ and
$$
\coker(\Psi_1\comp\mu)=
\bigoplus_{z\in S_1\cup S_2\cup S_3^+}\tilde{F}(z)$$
\e{lemma}
\be{proof} 
To determine $\coker(\Psi_1\comp\mu)$
we compute $\Psi_1([x\te x,0,0])$.
Again we may assume that $x\in G$.
There exist $z\in\im\,S$ and $g\in G$ such that $x^2=g^{-1}zg.$
Observe that $(gxg^{-1})^2=z$.
Now $$\phi(x\te x)=x^2\te x=g^{-1}zg\te x.$$
Notice that $\gamma(gx)=\gamma(g)$ since $gxg^{-1}\in\gz$.
\be{eqnarray*}
\tau([g^{-1}zg\te x,0,0])&=&
[\gamma(g)g^{-1}zg\gamma(g)^{-1}\te 
\gamma(g)x\gamma(gx)^{-1},0,0]\vspace{.5mm}\\
&=&[\gamma(g)g^{-1}zg\gamma(g)^{-1}\te 
\gamma(g)x\gamma(g)^{-1},0,0]\vspace{.5mm}\\
&=&\!\cases{
[z\te \gamma(g)g^{-1}gxg^{-1}g\gamma(g)^{-1},0,0]
&if $z$ is of\vspace{.5mm} type 1\cr
[z^{\pm1}\te \gamma(g)g^{-1}gxg^{-1}g\gamma(g)^{-1},0,0]
&if $z$ is of\vspace{.5mm} type 2\cr
[z\te \gamma(g)g^{-1}gxg^{-1}g\gamma(g)^{-1},0,0]
&if $z$ is of type 3\cr}
\e{eqnarray*}
applying the isomorphism $\eta$ we find
$$\Psi_1([x\te x,0,0])=
\cases{([gxg^{-1}],1)\in((G_z)_\#/<[z]>)\times C_2
&if $z$ is of\vspace{.5mm} type 1\cr
        [gxg^{-1}]\in(\gzs)_\#/<[z]>              
&if $z$ is of\vspace{.5mm} type 2\cr
        [gxg^{-1}]\in(\gz)_\#/<[z]>              
&if $z$ is of type 3\cr}$$
This proves the claim.
\e{proof}

The isomorphism $\Psi_2$ induces an isomorphism 
$$\Psi_3\colon\coker(1+\vartheta)\bijarrow
\coker(\Psi_2(1+\vartheta)\Psi^{-1}):$$
$$\diagram{
HQ_1(\mbbf\:_2[G])&\mapright{1+\vartheta}&\coker(\mu)
&\lra&\coker(1+\vartheta)\cr
\mapdown{\Psi}&&\mapdown{\Psi_2}&&\mapdown{\Psi_3}\cr
{\cal S}&\lra&\coker(\Psi_1\comp\mu)&
\lra&\coker(\Psi_2(1+\vartheta)\Psi^{-1})\cr}$$
\be{lemma}\label{identlem}
$\coker(\Psi_2(1+\vartheta)\Psi^{-1})$ arises from $\coker(\Psi_1\comp\mu)$
by imposing the following identifications.
For all $z$ of type
\be{enumerate}
\item[1]
identify $$([g],t^i)\in\tilde{F}(z) \mbox{ \  and \ }
([g],t^i)\in\tilde{F}(1)$$
\item[2]
identify $$[g]\in\tilde{F}(z) \mbox{ \  and \ }
\cases{
([g],t^{w(g)})\in\tilde{F}(z^2)& if $z^2$ is\vspace{.5mm} of type 1\cr
[g]\in\tilde{F}(z^2)& if $z^2$ is of type 2\cr}$$
\item[3]
identify $$[g]\in\tilde{F}(z) \mbox{ \  and \ }
\cases{
([g],1)\in\tilde{F}(z^2)& if $z^2$ is of\vspace{.5mm} type 1\cr
[g]\in\tilde{F}(z^2)& if $z^2$ is of \vspace{.5mm}type 2\cr
[g]\in\tilde{F}(z^2)& if $z^2$ is of type 3.\cr}$$
\e{enumerate}
\e{lemma}
\be{proof}
\be{enumerate}
\item
Let $([g],t^i)\in((\gz)_\#/<z>)\times C_2.$
$$\Psi^{-1}([g],t^i)=[g^{-1}z\te g,0,iz]$$
since 
\be{eqnarray*}
\phi(g^{-1}z\te g)&=&(z\te g) \quad\mbox{and}\quad
\phi(iz)=iz,\\
\tau([z\te g,0,iz])&=&[z\te g\gamma(g)^{-1},0,iz]\;=\;[z\te g,0,iz],\\
\eta([z\te g,0,iz])&=&([g],t^i).
\e{eqnarray*}
$\vartheta([g^{-1}z\te g,0,iz])=[g^{-1}z^2\te g+z\te z+iz\te z,0,iz^2]=
[g^{-1}\te g,0,i]$.
$$\Psi_2(g^{-1}\te g,0,i])=([g],t^i)\in\tilde{F}(z^2)=\tilde{F}(1)$$ 
since 
\be{eqnarray*}
\phi(g^{-1}\te g)&=&(1\te g) \quad\mbox{and}\quad  \phi(i)\;=\;i,\\
\tau([1\te g,0,iz])&=&[1\te g\gamma(g)^{-1},0,i]\;=\;[1\te g,0,i],\\ 
\eta([1\te g,0,i])&=&([g],t^i).
\e{eqnarray*}
\item
Let $[g,t^i]\in(\gzs\times_{C_2}C_4)_\#/<[z,t^2]>$.
$$\Psi^{-1}([g,t^i])=[g^{-1}z\te g,0,y],$$ 
where $y=\ent((i+1)/2)z+\ent(i/2)\zm$,
since 
\be{eqnarray*}
\phi(g^{-1}z\te g)&=&(z\te g)\quad\mbox{ and}\quad \phi(y)\;=\;y,\\ 
\tau([z\te g,0,y])&=&
[z\te g\gamma(g)^{-1},0,y]\;=\;[z\te g,0,y],\\
\eta([z\te g,0,y])&=&[g,t^{w(g)+2\ent(i/2)}]\;=\;[g,t^i].
\e{eqnarray*}
$\vartheta([g^{-1}z\te g,0,y])=
[g^{-1}z^2\te g,0,\ent((i+1)/2)z^2+\ent(i/2)z^{-2}]$.\nl 
Note that 
$[z\te z^{\pm1},0,0]=[z^{\pm1}\te z,0,0]=0$ in $\coker(\mu)$.\nl
Define $y'\isdef\ent((i+1)/2)z^2+\ent(i/2)z^{-2}$.
$$\Psi_2([g^{-1}z^2\te g,0,y'])=
\cases{
([g],t^i)\in\tilde{F}(z^2)& if $z^2$ is of type 1\cr
[g]\in\tilde{F}(z^2)& if $z^2$ is of type 2\cr}$$
since
\be{eqnarray*}
\phi(g^{-1}z^2\te g)&=&(z^2\te g) \quad\mbox{ and }\quad \phi(y')\;=\;y',\\
\tau([z^2\te g,0,y'])&=&[z^2\te g\gamma(g)^{-1},0,y']\;=\;
[z^2\te g,0,y'],\\
\eta([z^2\te g,0,y'])&=&
\cases{
([g],t^i)\in\tilde{F}(z^2)& if $z^2$ is of type 1\cr
[g]\in\tilde{F}(z^2)& if $z^2$ is of type 2\cr}
\e{eqnarray*}
\item
Let $[g]\in(\gz)_\#$.
$$\Psi^{-1}([g])=[g^{-1}z\te g,0,0]$$ 
since 
\be{eqnarray*}
\phi(g^{-1}z\te g)&=&(z\te g),\\  
\tau([z\te g,0,0])&=&
[z\te g\gamma(g)^{-1},0,0]\;=\;[z\te g,0,0],\\
\eta([z\te g,0,0])&=&[g].
\e{eqnarray*}
$\vartheta([g^{-1}z\te g,0,0])=[g^{-1}z^2\te g,0,0]$.
$$\Psi_2([g^{-1}z^2\te g,0,0])=
\cases{
([g],1)\in\tilde{F}(z^2)& if $z^2$ is of type 1\cr
[g]\in\tilde{F}(z^2)& if $z^2$ is of type 2\cr
[g]\in\tilde{F}(z^2)& if $z^2$ is of type 3\cr}$$
since 
\be{eqnarray*}
\phi(g^{-1}z^2\te g)&=&(z^2\te g),\\
\tau([z^2\te g,0,0])&=&[z^2\te g\gamma(g)^{-1},0,0]\;=\;
[z^2\te g,0,0],\\
\eta([z^2\te g,0,0])&=&
\cases{
([g],1)\in\tilde{F}(z^2)& if $z^2$ is of type 1\cr
[g]\in\tilde{F}(z^2)& if $z^2$ is of type 2\cr
[g]\in\tilde{F}(z^2)& if $z^2$ is of type 3\cr}
\e{eqnarray*}
\e{enumerate}
This completes the proof.
\e{proof}
\be{defi}
Let $G$ be a group. For every $z\in G$ we define $\sqrt{z}$ as
the subgroup of $(G_z)_\#$ resp. $(\ol{G_z})_\#$ generated by the set
$$\{g\in G\mid g^{2^k}=z \mbox{ for some } k\in\N\}.$$
\e{defi}
\be{defi}
Define
\[\J(G)\isdef\lim_{\stackler{z}{\lra}} F(z),\]
where \[F(z)\isdef\cases{
{\displaystyle {(G_z)_\#\over\sqrt{z}}}\times C_2
& if $ z$ is of\vspace{1mm} type $1$\cr 
{\displaystyle {(\ol{G_z})_\# \over\sqrt{z}}} 
& if $z$ is of\vspace{1mm} type $2$\cr
{\displaystyle {(G_z)_\#\over\sqrt{z}}}
& if $z$ is of type $3$\cr}\]
and the limit is taken with respect to the homomorphisms 
\be{enumerate}
\item[$\cdot$]
$F(z)\lra F(x^{-1}zx)$ for every $x\in G$
defined by
$$\cases{
([g],t^i)\mapsto([x^{-1}gx],t^i)& for all $z$ of\vspace{1mm} type $1$\cr
[g]\mapsto [x^{-1}gx]& for all $z$ of type $2$ and $3$\cr}$$
\item[$\cdot$]
$F(z)\ra F(z^2)$ defined by 
$$\cases{
([g],t^i)\mapsto([g],t^i)& for all $z$ of\vspace{1mm} type $1$\cr
[g]\mapsto\cases{([g],t^{w(g)})& if $z^2$ is of type $1$\cr
                  [g]          & if $z^2$ is of type $2$\cr}&
for all $z$ of\vspace{1mm} type $2$\cr
[g]\mapsto\cases{([g],1)& if $z^2$ is of type $1$\cr
                   [g]   & if $z^2$ is of type $2$\cr
                   [g]   & if $z^2$ is of type $3$\cr}&
for all $z$ of type $3$\cr}$$
\item[$\cdot$]
$F(z)\ra F(z^{-1})$ defined by 
$[g]\mapsto[g]$
for all $z$ of type $3$.
\e{enumerate}
\e{defi}
\be{remark}\label{remjgstruct}
$$\J(G)\cong\bigoplus_{{\displaystyle c\in\kl(G)}} \L(c),$$
where 
$$\L(c)\isdef \lim_{\stackler{z\in c}{{\displaystyle \lra}}} F(z).$$
\e{remark}
\be{thm}
\[\coker(1+\vartheta)\cong\J(G)\]
\e{thm}
\be{proof}
Obvious in view of lemma~\ref{identlem}.
\e{proof}
\be{prop}\label{propupsarf}
Suppose $\plane{g,h}$ is an element of $\Arf^h(G)$.
The invariant $\Upsilon$ of {\rm chapter III} maps $\plane{g,h }$ to 
\[\cases{[1,t]\in \L([1]) & if $gh$ is of type $1$\cr
         [h]\in \L([gh]) & if $gh$ is of type $2$\cr}\]
Note that $gh$ is never of type $3$.
\e{prop}
\be{proof}
Define $z\isdef gh$.
Note that $g^2=h^2=1$ and $hzh=\zm$.
By definition $\Upsilon(\plane{g,h})=[g\te h,gh]\in\coker(1+\vartheta).$
By the definitions of $\phi$, $\tau$ and $\eta$:  
\be{eqnarray*}
\phi(gh)&=&gh, \\
\phi(g\te h)&=&hg\te h,\\
\tau([hg\te h,0,gh])&=&\tau([hzh\te h,0,z])\\
&=&[\gamma(h)hzh\gamma(h)^{-1}\te \gamma(h)h,0,z]\\
&=&[z^{-1}\te h,0,z]\\
\eta([z^{-1}\te h,0,z])&=&\cases{
([h],t)\in\tilde{F}(z)& if $z$ is of type 1\cr
[h]\in\tilde{F}(z)& if $z$ is of type 2\cr}
\e{eqnarray*}
Hence
$$\Psi_3([g\te h,gh])=\cases{
([h],t)=([1],t)& if $gh$ is of type 1\cr
[h]& if $gh$ is of type 2\cr}$$
\e{proof}
\be{lemma}
For all $z\in G$
$$\ker(\gzs\ra(\gzs)_\#/\sqrt{z})\subset \gz.$$
\e{lemma}
\be{proof}
Every commutator is a product of squares: 
$xyx^{-1}y^{-1}=x^2(x^{-1}y)^2y^{-2}$.
Every square of an element in $\gzs$ belongs to $\gz$.
If $y^{{2^k}}=z$, then $y\in\gz$.
\e{proof}
\be{lemma}\label{lemmantgh}
$\Upsilon(\plane{g,h})$ is never trivial in $\J(G)$.
\e{lemma}
\be{proof}
Define $z\isdef (gh)^{2^k}$, with $k$ large.
If $z$ is of type 1, the statement is true by 
proposition~\ref{propupsarf}.
If $z$ is of type 2, then $g\in\gzs\setminus\gz$. Therefore 
$[g]\in \L([z])$ 
cannot be trivial.
\e{proof}
Now we review one of the examples we encountered in section 5 of chapter II.
\be{nitel}{Example}
Let $G$ be the group with presentation
$$G\isdef\langle X,Y,S\mid S^2=(XS)^2=(YS)^2=1,\quad XY=YX\rangle.$$
To compute $\Arf^h(G)$ we determine
$$\J(G)=\bigoplus_{{\displaystyle c\in\kl(G)}} \L(c).$$
It is immediately clear from the presentation of $G$ that $\sqrt{1}=G$ and
for all $z\in H$ we have $\gzs=G$. 
(Recall that $H$ is the subgroup generated by $X$ and $Y$.)
Therefore
$$\left(\gzs\right)_\#=G/G^2=G/\langle X^2,Y^2\rangle.$$
A little examination shows that
$$\kl(G)=
\left\{[1]\right\}\cup
\left\{[X^{2i}Y^{2j+1}],\,[X^{2k+1}Y^{2l}],\,[X^{2m+1}Y^{2n+1}]
\quad\mid j,k,m\geq 0\right\}$$
and
$$\L(c)=\cases{C_2& if $c=[1]$\cr
G/\langle X^2,Y\rangle \cong C_2\times C_2 & if $c=[X^{2i}Y^{2j+1}]$\cr
G/\langle X,Y^2\rangle \cong C_2\times C_2 & if $c=[X^{2k+1}Y^{2l}]$\cr
G/\langle X^2,XY\rangle\cong C_2\times C_2 & if $c=[X^{2m+1}Y^{2n+1}]$.\cr}$$
Proposition~\ref{propsuf1} of chapter II says that the elements
$$\cases{
\plane{1,1}&\cr 
\plane{X^{2i}Y^{2j+1}S,S} & for $j\geq 0$\cr
\plane{X^{2i+1}Y^{2j}S,S}&  for $i\geq 0$\cr
\plane{X^{2i+1}Y^{2j+1}S,S}& for $i\geq 0$\cr
\plane{X^{2i+1}Y^{2j+1}S,XS}& for $j\geq 0$\cr
\plane{X^{2i+1}Y^{2j+1}S,YS} & for $i\geq 0.$\cr
\plane{X^{2i}Y^{2j+1}S,XS} & for $j\geq 0$\cr
}$$
generate $\Arf^h(G)$.
But since
\be{eqnarray*}
\Upsilon(\plane{1,1})                 &=&t   \in C_2\\
\Upsilon(\plane{X^{2i}Y^{2j+1}S,S})   &=&[S] \in\L([X^{2i}Y^{2j+1}])\\
\Upsilon(\plane{X^{2i+1}Y^{2j}S,S})   &=&[S] \in\L([X^{2i+1}Y^{2j}])\\
\Upsilon(\plane{X^{2i+1}Y^{2j+1}S,S}) &=&[S] \in\L([X^{2i+1}Y^{2j+1}])\\
\Upsilon(\plane{X^{2i+1}Y^{2j+1}S,XS})&=&[XS]\in\L([X^{2i}Y^{2j+1}])\\
\Upsilon(\plane{X^{2i+1}Y^{2j+1}S,YS})&=&[YS]\in\L([X^{2i+1}Y^{2j}])\\
\Upsilon(\plane{X^{2i}Y^{2j+1}S,XS})  &=&[XS]\in\L([X^{2i-1}Y^{2j+1}])
\e{eqnarray*}
we may conclude that these elements constitute a basis for $\Arf^h(G)$.
\e{nitel}
We revert to one of the examples of chapter I.
\be{nitel}{Example}
Let $G$ be the group with presentation
$$G\isdef\langle X,Y,S\mid S^2=(XS)^2=(YS)^4=(Y^2S)^2=1,\quad XY=YX\rangle.$$
\be{prop}
\be{eqnarray*}
\left\{\plane{1,1}\right\}&\cup&
\left\{\plane{X^{2i+1}Y^{2j}S,S}\mid i\geq 0\right\}\\
&\cup&\left\{\plane{X^{2i}Y^{4j+2}S,S}\mid j\geq 0\right\}\\
&\cup&\left\{\plane{X^{2i+1}Y^{4j+2}S,XS}\mid j\geq 0\right\}
\e{eqnarray*}
is a basis for $\Arf^{s,h}(G)$.
\e{prop}
\be{proof}
We know already that these elements generate $\Arf^h(G)$.
To prove independence we use our invariant $\Upsilon$.
We proceed to compute the summands $\L(c)$ of value group $\J(G)$.
It is not hard to verify that
$$\kl(G)=
\left\{[1]\right\}\cup
\left\{[X^{2i+1}Y^{2j}]\,\mid i\geq 0\right\}\cup
\left\{[X^{i}Y^{2j+1}]\,\mid j\geq 0\right\}.$$
We omit the proof.
$$\L(c)=\cases{
C_2
& if $c=[1]$\cr
G/\langle X,Y^2,(YS)^2\rangle \cong C_2\times C_2 
& if $c=[X^{2i+1}Y^{2j}]$\cr
G/\langle X^2,Y,(YS)^2\rangle \cong C_2\times C_2 
& if $c=[X^{2i}Y^{2j+1}]$ \cr
G/\langle X^2,XY,(YS)^2\rangle\cong C_2\times C_2 
& if $c=[X^{2i+1}Y^{2j+1}]$.\cr}$$
Note that the class of $S$ is non-trivial in any $\L(c)$.
Further, the classes of $X$, $S$ and $XS$ in $\L([X^{2i}Y^{2j+1}])$ as well
as in $\L([X^{2i+1}Y^{2j+1}])$ are distinct.
Now we can use the list of images
\be{eqnarray*}
\Upsilon(\plane{1,1})                 &=&t   \in C_2\\
\Upsilon(\plane{X^{2i+1}Y^{2j}S,S})   &=&[S] \in\L([X^{2i+1}Y^{2j}])\\
\Upsilon(\plane{X^{4i}Y^{4j+2}S,S})   &=&[S] \in\L([X^{2i}Y^{2j+1}])\\
\Upsilon(\plane{X^{4i+2}Y^{4j+2}S,S}) &=&[S] \in\L([X^{2i+1}Y^{2j+1}])\\
\Upsilon(\plane{X^{4i+1}Y^{4j+2}S,XS})&=&[XS]\in\L([X^{2i}Y^{2j+1}])\\
\Upsilon(\plane{X^{4i+3}Y^{4j+2}S,XS})&=&[XS]\in\L([X^{2i+1}Y^{2j+1}])
\e{eqnarray*}
to see that the assertion is true.
\e{proof}
\e{nitel}
\be{nitel}{Example}
Let $G$ be the group with presentation
$$\langle X,Y,Z,S\mid X,Y,Z \mbox{ commute },
S^2=(XS)^2=(YS)^2=(ZS)^2=1\rangle.$$
Let $c\in\kl(G)$ be the class of $XYZ$. 
The invariant $\Upsilon$ maps 
$$\xi\isdef\plane{XYS,SZ}+\plane{XZS,SY}+\plane{YZS,SX}+\plane{XYZS,S}
\in\Arf^h(G)$$
to the class $[SZSYSXS]=[1]\in \L(c)=G/\langle X^2,Y^2,Z^2,XYZ\rangle.$
But it is not clear at all whether $\xi$ is trivial in $\Arf^h(G)$.
\e{nitel}

\newpage
\section{Groups with two ends.}
\setcounter{altel}{0}
\setcounter{equation}{0}

We wish to prove that our invariant $\Upsilon$ is injective for all groups
having two ends.
For that purpose theorem~\ref{chargp2e} gives a suitable characterization
of these groups.
\be{nota}
Throughout this section
\[\be{array}{ll}
G     & \mbox{ denotes a group,}\\
E     & \mbox{ denotes a finite group,}\\
C     & \mbox{ denotes the infinite cyclic group,}\\
C_m   & \mbox{ denotes the cyclic group of order $m$,}\\
D     & \mbox{ denotes the infinite dihedral group }\\ 
      & \mbox{ with presentation }<S,T\mid S^2=(ST)^2=1>,\\
D_{m}& \mbox{ denotes the dihedral group of order $2m$}\\ 
      & \mbox{ with presentation }
        <\sigma,\tau\mid \sigma^2=(\sigma\tau)^2=\tau^m=1>.
\e{array}\]
\e{nota}

\be{thm} {\rm \cite{Wall;gpth} }
The following statements are equivalent;
\be{enumerate}
\item $G$ has two ends.
\item $G$ has an infinite cyclic subgroup of finite index.
\item $G$ has an infinite cyclic normal subgroup of finite index.
\e{enumerate}
\e{thm}
\be{proof}
We refer to {\em loc. cit.} for a proof.
\e{proof}
\be{defi}
A group extension of $C$ by $E$ is a short exact sequence of groups and
homomorphisms
\[1\ra C\ra G\ra E\ra1\]
The extension is called central if the image of $C$ is central in $G$.
\e{defi}
\be{thm} \label{chargp2e}
\begin{enumerate}
\item 
$1\ra C\ra G\ra E\ra1$ is a central extension if and only if $G$ fits into a
pull-back diagram
\be{eqnarray*}G&\ra&E\\\da&&\da\\C&\ra&C_m\e{eqnarray*}
\item 
$1\ra C\ra G\ra E\ra1$ is a non-central extension if and only if $G$
fits into a pull-back diagram 
\be{eqnarray*}G&\ra&E\\\da&&\da\\D&\ra&D_m\e{eqnarray*}   
\e{enumerate}
\e{thm}
\be{proof}
In the sequel we will regard $C$ as a subgroup of $G$.
\be{itemize}
\item[{\em 1.}] ``$\Rightarrow$''
Suppose $1\ra C\ra G\mapright{\pi}E\ra1$
is a central extension. 
Define the so-called transfer homomorphism $\phi\colon G\ra C$ as follows:
choose a set theoretic section $\alpha\colon E\ra G$ 
of the projection $\pi\colon G\ra E$ such that $\alpha(1)=1$ and define
$$\phi(g)\isdef\prod_{e\in E}\alpha(e)g\alpha(e\pi(g))^{-1} 
\mbox{ \ for all \ } g\in G.$$
Note that
\be{itemize}
\item[$\diamond$]
$\alpha(e)g\alpha(e\pi(g))^{-1}\in\ker\pi=C$ for all $e\in E$ and $g\in G$.
\item[$\diamond$]
$\phi$ does not depend on the choice of $\alpha$:\nl
If $\alpha'$ is another section of $\pi$ we have
$\alpha(e)\alpha'(e)^{-1}\in C$ for every $e\in E$.
Hence
\be{eqnarray*}
\prod_{e\in E}\alpha(e)g\alpha(e\pi(g))^{-1}&=&
\prod_{e\in E}\alpha'(e)g\alpha'(e\pi(g))^{-1}\cdot\\
&&\prod_{e\in E}\alpha(e)\alpha'(e)^{-1}\cdot\\
&&\prod_{e\in E}\alpha'(e\pi(g))\alpha(e\pi(g))^{-1}\\
&=&\prod_{e\in E}\alpha'(e)g\alpha'(e\pi(g))^{-1}
\e{eqnarray*}
\item[$\diamond$]
$\phi$ is a homomorphism:
\be{eqnarray*}
\phi(g_1g_2)&=&\prod_{e\in E}\alpha(e)g_1g_2\alpha(e\pi(g_1g_2))^{-1}\\
            &=&\prod_{e\in E}\alpha(e)g_1\alpha(e\pi(g_1))^{-1}\cdot
               \alpha(e\pi(g_1))g_2\alpha(e\pi(g_1)\pi(g_2))^{-1}\\
            &=&\phi(g_1)\phi(g_2)
\e{eqnarray*}
\item[$\diamond$]
$\phi(c)=c^{|E|}$ for every $c\in C$. 
Here $|E|$ denotes the cardinality of $E$.
\e{itemize}
Now it is easy to verify that $G$ fits into the pull-back diagram
$$\diagram{G&\mapright{\pi}&E\cr\mapdown{\phi'}
&&\mapdown{p\phi'\alpha}\cr C&\mapright{p}&C_m\cr}$$
where
$m\isdef|E|/[C:\im\phi]$,\nl
$[C:\im\phi]$ is the index of $\im\phi$ in $C$,\nl
$\phi'\isdef\epsilon\comp\phi$,\nl
$\epsilon$ is an isomorphism $\im\phi\ra C$ and \nl
$p\colon C\ra C_m$ is the canonical projection.\nl
Note that $p\phi'\alpha$ does not depend on $\alpha$.
\item[$\phantom{2.}$] ``$\Leftarrow$'' 
Suppose
$$\diagram{G&\mapright{\pi}&E\cr\mapdown{}
&&\mapdown{}\cr C&\mapright{p}&C_m\cr}$$
is a pull-back diagram, then
$$\diagram{1\lra&C\lra   & G     &\mapright{\pi}E\lra1\hfil\cr
                &c\mapsto&(c^m,1)&                    \hfil\cr
                &        &(c,e)  &\mapsto e           \hfil\cr}$$
is a central extension.
\item[{\em 2.}] ``$\Rightarrow$'' 
Suppose $1\ra C\ra G\mapright{\pi}E\ra1$ is a non-central extension.
Choose a set theoretic section $\alpha\colon E\ra G$ as before.
The homomorphism $$w\colon E\ra \aut(C)\cong C_2$$ defined by
$w(e)(c)\isdef\alpha(e)c\alpha(e)^{-1}$ for all $c\in C$ and $e\in E$,
does not depend on the choice of $\alpha$.
Let $\phi\colon\ker(w\pi)\ra C$ be the transfer homomorphism associated to
the central extension
$$1\ra C\ra \ker(w\pi)\mapright{\pi}\ker(w)\ra1.$$
Choose an element $u\in G\setminus\ker(w\pi)$ and
define  $$\psi\colon G\ra D$$ 
$$\psi(g)\isdef\cases{\phi(g)&if $g\in\ker(w\pi)$\cr
                 \phi(gu^{-1})S&if $g\in G\setminus\ker(w\pi)$\cr}.$$
For every $g\in\ker(w\pi)$ we equate
\be{eqnarray*}
\phi(ugu^{-1})^{-1}&=&u^{-1}\phi(ugu^{-1})u\\
&=&\prod_{e\in\ker(w)}u^{-1}\alpha(e)ugu^{-1}\alpha(e\pi(ugu^{-1}))^{-1}u\\
&=&\prod_{e\in\ker(w)}u^{-1}\alpha(e)u\alpha(\pi(u)^{-1}e\pi(u))^{-1}\cdot\\
&&\prod_{e\in\ker(w)}\alpha(\pi(u)^{-1}e\pi(u))g
\alpha(\pi(u)^{-1}e\pi(u)\pi(g))^{-1}\cdot\\
&&\prod_{e\in\ker(w)}\alpha(\pi(u)^{-1}e\pi(u)\pi(g))u^{-1}
\alpha(e\pi(ugu^{-1}))^{-1}u\\
&=&\phi(g).
\e{eqnarray*}
In particular $\phi(u^2)=1$ and $\psi$ is a homomorphism.\nl
Again it is easy to verify that $G$ fits into the pull-back diagram
$$\diagram{G&\mapright{\pi}&E\cr\mapdown{\psi'}
&&\mapdown{p\psi'\alpha}\cr D&\mapright{p}&D_{m}\cr}$$
where
$2m\isdef|E|/[D:\im\psi]$.\nl
Note that $m\cdot[D:\im\psi]=|\ker(w)|$ and
$|E|=2|\ker(w)|$.\nl
$\psi'\isdef\epsilon\comp\psi$,\nl
$\epsilon$ is an isomorphism $\im\psi\ra D$ and \nl
$p\colon D\ra D_m$ is the canonical projection.\nl
Note that $p\psi'\alpha$ does not depend on $\alpha$.
\item[$\phantom{2.}$] ``$\Leftarrow$'' 
If
$$\diagram{G&\mapright{\pi}&E\cr\mapdown{}
&&\mapdown{}\cr D&\mapright{p}&D_{m}\cr}$$
is a pull-back diagram, then
$$\diagram{1\lra&C\lra   & G     &\mapright{\pi}E\lra1\hfil\cr
                &c\mapsto&(c^m,1)&                    \hfil\cr
                &        &(d,e)  &\mapsto e           \hfil\cr}$$
is obviously a non-central extension.
\e{itemize}
This completes the proof.
\e{proof}

\newpage
\section{$\Upsilon$ for groups with two ends.}
\setcounter{altel}{0}
\setcounter{equation}{0}

This section is devoted to the following theorem.
\be{thm}\label{thmin2end}
The invariant $\Upsilon\colon\Arf^h(G)\ra\J(G)$ is injective for all 
groups $G$ having two ends.
\e{thm}
\be{lemma}\label{lem2powz}
For all $k\in\N$ the relation
$$\plane{a,b}=\plane{a,a(ab)^{2^k}}=\plane{b,b(ab)^{2^k}}$$
holds in $\Arf^h(G)$.
\e{lemma}
\be{proof}
By the relations mentioned in remark~\ref{remarfgrel} of chapter I,
we have
$\plane{a,b}=\plane{a,bab}=\plane{a,a(ab)^2}$ and $\plane{a,b}=\plane{b,a}$.
The rest is obvious.
\e{proof}
\be{lemma}\label{lemmaeindord}
If $\plane{a,az}$ and $\plane{b,bz}$ are elements of $\Arf^h(G)$ and 
$abz^i$ has finite order, for some $i\in\Z$,
then $\plane{a,az}=\plane{b,bz}$.
\e{lemma}
\be{proof}
First we consider the case that $i=0$.
Let us write $x=ab$ and let's say that the order 
of $x$ equals $2^km$ with $m$ odd. \nl
Note that $a^2=b^2=(az)^2=1 $ and $xz=zx$.
Using the relations of remark~\ref{remarfgrel} of chapter I we equate
\be{eqnarray*}
\plane{b,bz}&=&\plane{ax,axz}\\
&=&\plane{ax^m,ax^mz}\\
&=&\plane{ax^m,ax^mz^{2^k}}\\
&=&\plane{ax^m,ax^m(ax^maz)^{2^k}}\\
&=&\plane{ax^m,az}\\
&=&\plane{az,ax^m}\\
&=&\plane{az,az(azax^m)^{2^k}}\\
&=&\plane{az,az(aza)^{2^k}}\\
&=&\plane{az,a}\\
&=&\plane{a,az}
\e{eqnarray*}
The case $i\neq 0$ can be reduced to the previous case:
$$\plane{b,bz}=\plane{z^{-j}bz^j,z^{-j}bzz^j}=\plane{bz^{2j},bz^{2j+1}}$$
$$\plane{b,bz}=\plane{b,bbzb}=\plane{b,bz^{-1}}=\plane{bz^{2j},bz^{2j-1}}
=\plane{bz^{2j-1},bz^{2j}},$$
thus $\;\plane{b,bz}=\plane{bz^i,bz^{i+1}}$.
\e{proof}
\be{prop}
In the case where $G$ fits into a pull-back diagram 
$$\diagram{G&\lra&E\cr\mapdown{}&&\mapdown{}\cr C&\lra&C_m\cr}$$ 
$\Upsilon$ is injective.
\e{prop}
\be{proof}
Let $x\in \ker(\Upsilon)$.
The relations in $\Arf^h(G)$ listed in remark~\ref{remarfgrel} of chapter I
and remark~\ref{remjgstruct} on the structure of $\J(G)$ allow us to assume, 
without loss of generality,
that $$x=\sum\plane{a_i,a_iz}.$$
Every product of two elements of order two, is of finite order, 
since all elements of order two in $G$ 
take the form $(1,e)$.
So we may use lemma~\ref{lemmaeindord} to see that $x=0$ or $x=\plane{a,az}$.
But according to
lemma~\ref{lemmantgh} $\Upsilon(\plane{a,az})$ is non-trivial,
so the second case does not occur.
\e{proof}
It remains to show that $\Upsilon$ is injective for groups $G$ which fit
into a pull-back diagram 
$$\diagram{G&\mapright{\pi}&E\cr
\mapdown{\psi}&&\mapdown{\hat p}\cr
D&\mapright{p}&D_{2m}\cr}$$
\underline{Intermezzo}.

\be{defi}
Let $G$ be a group.
Suppose we have 2-primary elements $a,b\in G$ which satisfy
$[a,b^2]=[a^2,b]=1$. 
Here $[x,y]$ denotes the commutator $xyx^{-1}y^{-1}$.
Denote by $H$ the subgroup of $G$ generated by $a$ and $b$. The matrix
$$\pmatrix{a&1\cr 0&b\cr}+\pmatrix{a&1\cr 0&b\cr}^\alpha=
\pmatrix{a+a^{-1}&1\cr 1&b+b^{-1}\cr}=
\pmatrix{a(1+a^{-2})&1\cr 1&b(1+b^{-2})\cr}$$
is invertible, since $1+a^{-2}$ and $1+b^{-2}$ are nilpotent and central  
in $\mbbf\,_2[H]$.
It is therefore legitimate to define 
$$\wp(a,b)\isdef\left[\pmatrix{a&1\cr 0&b\cr}\right]-
\left[\pmatrix{0&1\cr0&0\cr}\right]\in L^h(H).$$
We call such elements of $L^h(H)$ as well as their images in $L^h(G)$
pseudo-arfian.
\e{defi}
Notice that $\wp(a,b)$ is not necessarily an element of 
$\Arf^h(H)$ or $\Arf^h(G)$.
However, applying theorem~\ref{iadiciso} of chapter I to the ring
$\;\mbbf\,_2[H]$ and its nilpotent ideal
$(a^2+1,b^2+1)$ yields an isomorphism  
$$L^h(H)\bijarrow L^h(H/<a^2,b^2>),$$  which maps $\wp(a,b)\in L^h(H)$ 
to the Arf-element
$\plane{a,b}\in \Arf^h(H/<a^2,b^2>):$
$$\diagram{&&\wp(a,b)&\longmapsto&\plane{a,b}\cr
\Arf^h(G)&\longleftarrow&\Arf^h(H)&\injarrow&\Arf^h(H/<a^2,b^2>)\cr
\bigcap&&\bigcap&&\bigcap\cr
L^h(G)&\longleftarrow&L^h(H)&\bijarrow&L^h(H/<a^2,b^2>)\cr}$$            
\be{punt}\label{defpseu}
Let $G$ be a group and $g,z\in G$.
Assume $g^{-1}zg=z^{-1}$ and $g$ is of finite order, say $2^rr_0$ 
with $r_0$ odd. Define  $H$ as the subgroup of $G$
generated by $z$ and $h\isdef g^{r_0}$.\nl
Since $h^{-1}zh=z^{-1}$, i.e. $h^2=(hz)^2$
we obtain a pseudo-arf element $\wp(h,hz)\in L^h(H).$
The question is whether this element depends on $h$.
\e{punt}
\be{thm}\label{thmwinjend}
Let $E$ be a finite group.
The invariant $\omega_1^h$ of chapter II induces an isomorphism
$$L^h(E)\lra\bigoplus k/\{x+x^2\mid x\in k\}$$
Here the summation runs through all representations $\rho\colon
\mbbf\,_2[E]\ra M_n(k)$ which take the form
$\rho(e)^{-1}=P^{-1}\rho(e)^tP$ for all $e\in E$,
for some invertible matrix $P$ and $t$ means matrix transpose.
What's more, the image of $\wp(h,hz)$ under this isomorphism
is the element which has $[\tr(\rho(z))]$ at the place with index $\rho$.
In particular $\wp(h,hz)$ does not depend on $h$.
\e{thm}
\be{proof}
Define $R\isdef\mbbf\,_2[E]/\rad(\mbbf\,_2[E])$ where $\rad$ means 
Jacobson radical.
For every ring $A$, we denote by $\widetilde{A}$ the truncated polynomial 
ring $A[T]/(T^3)$ and in this context
$(T)$ is the two-sided ideal of $\widetilde{A}$ generated by $T$.
Consider the following diagram.
$$\diagram{
L^h(E)&&\cr
\isodown{1}&&\cr
L^h(R)&{\buildrel \omega_1^h \over {\hbox to 50pt{\rightarrowfill}}}
&H^0(K_1(\widetilde{R}))\cr
\isodown{2}&&\isodown{6}\cr
L^h\left(\prod D_i\right)&&
H^0\left(\bigoplus K_1(\widetilde{D_i})\right)\cr
\isodown{3}&&\isodown{7}\cr
\bigoplus_j L^h(D_j)&&
\bigoplus_j H^0\left( K_1(\widetilde{D_j})\right)\cr
\isodown{4}&&\isodown{}\cr
\bigoplus_j L^h(k_j)&&
H^0\left(\bigoplus( k_j^*\oplus1+T\widetilde{k_j})\right)\cr
\isodown{5}&&\isodown{}\cr
\bigoplus_j\Arf^h(k_j)&&
H^0\left(\bigoplus_j 1+T\widetilde{k_j}\right)\cr
\hfill\searrow& &\swarrow\hfill\cr
&\bigoplus_j\coker(1+\sigma_j)&\cr
&\mapdown{\cong}&\cr
&\bigoplus_j \Z/2&\cr}$$
We elucidate the diagram.
\be{enumerate}
\item
It follows from theorem~\ref{iadiciso} of chapter I 
that $L^h(E)$ and $L^h(R)$ are isomorphic,
because $\rad(\mbbf\,_2[E])$ is a nilpotent ideal.
\item
The ring $R$ is artinian and $\rad(R)=0$, so we can apply the 
Wedderburn-Artin theorem. In our case this reads:
$R$ is isomorphic to a direct product of full matrix rings over finite 
fields of characteristic two.
Explicitly, $$R\cong \prod D_i;$$
here $D_i\isdef M_{n_i}(k_i)$ is the ring of $(n_i\times n_i)$-matrices
over the finite field $k_i$ and char$(k_i)=2$. \nl
Denote by $\rho_i$ the composition 
$\,\mbbf\,_2[E]\ra R\ra \bigoplus D_i\ra D_i.$
\item
Let $\,\conj\,$ denote the (anti-) involutions on 
$R$ and $\prod D_i$ induced
by the involution on $\mbbf\,_2[E]$.
Before we decompose $\prod D_i$ as a product of rings with involution,
we fix some notations concerning finite fields of characteristic two.
If $k$ is such a field, the group of automorphisms of $k$ is cyclic and 
generated by the Frobenius automorphism $\sigma\colon k\ra k$ which 
assigns to an element $x$ of $k$ its square.
The field trace $\tr\colon k\surarrow \mbbf\,_2$ induces an isomorphism
$\coker(1+\sigma)\ra\mbbf\,_2$.
If the degree of $k$ over $\mbbf\,_2$ is even there 
exists a unique automorphism of order two, which we denote by 
$\hat{\sigma}$.\nl
Now, for a factor $D=M_n(k)$ of $\prod D_i$ we have three possible cases:
\be{itemize}
\item
$D$ is invariant under the involution, i.e. $D=\ol{D}$, and the restriction
of $\,\conj\,$ to $k$ is $\hat{\sigma}$.
\item
$D=\ol{D}$ and the restriction of $\,\conj\,$ to $k$ is the identity.
Since the composition of the (anti-) involution $\,\conj\,$ with matrix 
transpose is a $k$-linear automorphism of $D$, this composition takes the 
form $X\mapsto PXP^{-1}$, for some invertible matrix $P$.
Further we may assume that $P$ is symmetric,
since this automorphism is of order two.
Thus for all $X\in D$ we have $\ol{X}=P^{-1}X^tP$.
\item
$\ol{D}\neq D$. So $D\times\ol{D}$ is a factor of $\prod D_i$.
If $D\times D^{\circ}$ is endowed with the involution $(x,y)\mapsto (y,x)$,
the map $D\times\ol{D}\ra D\times D^{\circ}$ defined by 
$(x,y)\mapsto(x,\ol y)$ 
is an isomorphism of rings with involution. 
Recall that ${\scriptstyle \circ}$ means opposite multiplication. 
\e{itemize}
Thus we obtain a decomposition of $\prod D_i$ in which
three different types of factors occur.
The $L$-groups split accordingly. See e.g. \cite{w2}.
We assert that only the groups $L^h(D)$, 
where $D$ is of the second type, survive.\nl
In the first case we have
$$L^h(D,\conj,1)\cong L^h(k,\hat{\sigma},1)$$
by Morita invariance.
But $L^h(k,\hat{\sigma},1)=0$ by \cite[\S6]{Wall;lfound}.\nl
For the third case we will show that quadratic 
modules $(M,\theta)$ over the ring $D\times D^{\circ}$,
with involution $\alpha(x,y)=(y,x)$, are in fact hyperbolic. 
Note that there is no need to worry about bases, because we are working in
$L^h$.
Define $\lambda\isdef(1,0)$,
$M_1\isdef \lambda M$ and $M_2\isdef(1+\lambda)M$. So $M=M_1\oplus M_2$.
Since $b_\theta\colon M\ra M^\alpha$ is an isomorphism and for all $m,n\in M$
$$b_\theta((1+\lambda)m)((1+\lambda)n)=\lambda(1+\lambda)b_\theta(m)(n)=0,$$
the restriction of $b_\theta$ to $M_2$ yields an isomorphism
$M_2\ra M_1^\alpha$.
Now it is easy to verify that the map
$$(M,\theta)\lra H(M_1)=(M_1\oplus M_1^\alpha,\upsilon)$$
defined by
$$ m\mapsto(\lambda m,b_\theta((1+\lambda)m))$$
is an isometry.
In Walls terminology \cite{Wall;lfound} this says 
that $(D\times\{0\})M$ is a subkernel of $M$.
This proves our assertion.\nl
We use the index $j$ to refer to summands of the second type.
\item 
$L^h(D_j,\conj,1)$ is isomorphic to $L^h(k_j,1,1)$ by 
Morita invariance.
\item
Since the field trace $\tr\colon k\lra\mbbf\,_2$ is surjective we can choose
an element $a$ in $k$ such that $\tr(a)=1$.
The Arf invariant $$\omega_1^h\colon \Arf^h(k,1,1)\ra \coker(1+\sigma)\cong \Z/2$$
maps $\plane{a,1}$ to the non-trivial element in $\Z/2$.
Combining this with the fact that $L^h(k,1,1)\cong\Z/2$, 
see \cite[\S6]{Wall;lfound},
we find $$L^h(k,1,1)=\Arf^h(k,1,1)\cong\Z/2.$$ 
\item
It is almost immediately clear from the definition of $K_1$,
that
for all rings $A,B$ one has $K_1(A\times B)\cong K_1(A)\oplus K_1(B)$.
\item
The argument here is roughly the same as the one on the `$L$-side'
of the diagram (item 3). 
By Morita theory we have $K_1(\widetilde{D})\cong K_1(\widetilde{k})$.
Alternatively one can see this directly by looking at the 
definition of $K_1$.
Since the projection  $\widetilde{k}\ra k$ splits we have
$$K_1(\widetilde{k})\cong K_1(k)\oplus K_1(\widetilde{k},(T)).$$
It is well-known that $K_1(k)=k^*$, the group of units in $k$ 
and we already saw that
$K_1(\widetilde{k},(T))=1+T\widetilde{k}$.
To decompose $\bigoplus_i K_1(\widetilde{D_i})$ into invariant parts,
we consider the same three possibilities:
\be{itemize}
\item
$D$ is invariant under the involution and the restriction
of $\,\conj\,$ to $k$ is $\hat{\sigma}$.
In this case 
$$H^0(K_1(\widetilde{D}))=H^0(k^*\oplus (1+T\widetilde{k}))
=H^0(k^*)\oplus H^0(1+T\widetilde{k}).$$
But $H^0(k^*)$ vanishes, because $k^*$ has odd order.
And in the third section of chapter II we computed
$H^0(1+T\widetilde{k})=C(k)$, but this also disappears since
$H^0(k;\hat{\sigma})=0$. 
\item
$D=\ol{D}$ and the restriction of $\conj$ to $k$ is the identity.
By the same arguments as in the previous case we obtain
$H^0(K_1(\widetilde{D}))=C(k)$, but now $C(k)$ is precisely $\coker(1+\sigma)$.
\item
$D\neq\ol{D}$. 
Here the involution interchanges the summands
$K_1(\widetilde{D})$ and $K_1(\widetilde{\ol{D}})$, so $H^0$ clearly dies.
\e{itemize}
Thus only the summands of the second type survive.
\e{enumerate}
This completes the proof of the first part of what the theorem asserts.
To prove the second part
let $\rho\colon \mbbf\,_2[E]\ra D=M_n(k)$ be a representation of the special
kind, i.e. $D$ is of the second type, and assume that
$\rho(h)=H$ and $\rho(z)=Z$.
Then 
$$\omega_1^h\left(\wp(h,hz)\right)=
\omega_1^h\left(\left[\pmatrix{h&1\cr0&hz\cr}\right]-
      \left[\pmatrix{0&1\cr0&0\cr}\right]\right).$$
Now
$$\pmatrix{h&1\cr0&hz\cr}+(1+T)\pmatrix{h&1\cr0&hz\cr}^\alpha=
\pmatrix{h+(1+T)h^{-1}&1\cr(1+T)&hz+(1+T)(hz)^{-1}\cr}$$
and
$$\pmatrix{0&1\cr0&0\cr}+(1+T)\pmatrix{0&1\cr0&0\cr}^\alpha=
\pmatrix{0&1\cr(1+T)&0\cr},$$
so
\be{eqnarray*}
\lefteqn{\omega_1^h\left(\wp(h,hz)\right)=}\\
&&\left[\pmatrix{h+(1+T)h^{-1}&1\cr(1+T)&hz+(1+T)(hz)^{-1}\cr}
\pmatrix{0&1\cr(1+T)&0\cr}^{-1}\right]=\\
&&\left[\pmatrix{1&h(1+T)^{-1}+h^{-1}\cr hz+(1+T)(hz)^{-1}&1\cr}\right]=\\
&&\left[\pmatrix{1+h^2z(1+T)^{-1}+h^{-2}z(1+T)&h(1+T)^{-1}+h^{-1}\cr
0&1\cr}\right]=\\
&&\left[\pmatrix{1+h^2z(1+T)^{-1}+h^{-2}z(1+T)&0\cr
0&1\cr}\right]=\\
&&\left[\left(1+h^2z(1+T)^{-1}+h^{-2}z(1+T)\right)\right]\in 
H^0(K_1(\widetilde{R})).
\e{eqnarray*}
The image of this element in $H^0(1+T\widetilde{k})$ equals
\be{eqnarray*}
\lefteqn{\left[\det\left(1+H^2Z(1+T+T^2)+H^{-2}Z(1+T)\right)\right]} \\
&=&\left[\det\left(1+H^2ZT^2\right)\right] 
\left[\det\left(1+(H^2+H^{-2})Z(1+T)\right)\right]\\
&= &\left[\det\left(1+H^2ZT^2\right)\right] \\
&=&[1+\tr(H^2Z)T^2]\\
&=&[1+\tr(Z)T^2] 
\e{eqnarray*}
Here $$\det\left(1+(H^2+H^{-2})Z(1+T)\right)=1 \quad\mbox{ and }\quad
\tr(H^2Z)=\tr(Z)$$ 
by lemma~\ref{lemmanilp},
because $(H^2+H^{-2})Z$ and $1+H^2$ are nilpotent.
Finally, the image of $[1+\tr(Z)T^2]$ in $\coker(1+\sigma)$ equals
$[\tr(Z)]=[\tr(\rho(z))]$ according to the computations in chapter II.
\e{proof}
\be{lemma}\label{lemmanilp}
If $V$ is a finite dimensional $k$-vectorspace, $N\colon V\ra V$ is a
nilpotent linear map and $s$ is an indeterminate, 
then $$\tr(N)=0 \mbox{ \ and \ } \det(1+sN)=1.$$
\e{lemma}
\be{proof}
Suppose $N^n=0$.
We apply induction on $n$. \nl
If $n=1$ the matter is clear.\nl
If $n>1$ consider the diagram
$$\diagram{0&\lra& NV        &\lra&V          &\lra&V/NV       &\lra&0\cr
            &    &\mapdown{N}&    &\mapdown{N}&    &\mapdown{0}&    & \cr
           0&\lra& NV        &\lra&V          &\lra&V/NV       &\lra&0\cr}$$
The first vertical map in this diagram has nilpotency degree $n-1$
and $N\colon V\ra V$ takes the form $\pmatrix{*&*\cr0&0\cr}$.
This proves the assertions.
\e{proof}
\underline{End intermezzo}.

\be{thm}\label{thm2stuks}
If $\plane{a,b}+\plane{c,d}\in\ker(\Upsilon)$, then
$\plane{a,b}=\plane{c,d}$ in $\Arf^h(G)$.
\e{thm}
\be{proof}
Note that $\plane{a,b}+\plane{c,d}\in\ker(\Upsilon)$ if and only if 
$$[ab]=[cd]\in \kl(G)\quad\mbox{ and }\quad
[bd]=1\in\L([ab]).$$
Again the relations in $\Arf^h(G)$ and the structure of $\J(G)$
allow us to assume that $ab=cd$.
Elements of order two in $G$ either have the form
$(1,e)$ with $e^2=1$ and $\hat p(e)=1$ or the form
$(ST^i,e)$ with $e^2=1$ and $\hat p(e)=p(ST^i).$
Thus we may assume that
$$\plane{a,b}+\plane{c,d}=
((\Delta,e),(\Delta T^i,ez))+((\Delta T^j,ex),(\Delta T^{i+j},exz)),$$
where 
\be{itemize}
\item[$\cdot$]
$\Delta=ST^\nu$: if $\Delta=1$ we are through by lemma~\ref{lemmaeindord}
\item[$\cdot$]
$e,x,z\in E$ satisfy $e^2=(ex)^2=(ez)^2=1$ and $xz=zx$
\item[$\cdot$]
$[(T^j,x)]=1\in \L([(T^i,z)])$. 
\e{itemize}
Lemma~\ref{lem2powz} permits us to
replace $(T^i,z)$ by any power-of-two power
of $(T^i,z)$. Hence we may assume that $z$ has odd order, let's say order
$l_0=2l-1$.
\be{nitel}{Case 1$\colon$ $i\neq0$}\nl
Write $m=2^\mu m_0$ and $i=2^\tau i_0$ with $m_0$ and $i_0$ odd. 
Since $(T^i,z)\in G$ and $z$ has order $l_0$ in $E$, we have $m|il_0$, i.e.
$\mu\leq\tau$ and $m_0|i_0l_0$. 
If $\mu<\tau$, then 
$$((\Delta,e),(\Delta T^i,ez))=((\Delta,e),(\Delta T^{i/2},ez^l)),$$
by lemma~\ref{lem2powz}, 
where
\be{itemize}
\item[$\cdot$]
$(T^{i/2},z^l)\in G$, \ because $il\equiv i/2\pmod{m}$
\item[$\cdot$]
$(T^{i/2},z^l)^2=(T^i,z).$
\e{itemize}
So we may assume that $\mu=\tau$.\nl
Further, conjugation by a suitable power of $(T^j,x)$ allows us 
to replace $(T^j,x)$ by any odd power of $(T^j,x)$.
Thus we may assume that $x$ has order a power of 2, let's say $2^k$.
\be{itemize}
\item[$\diamond$]
If necessary we conjugate by $(T^m,1)$ to achieve that $0\leq j<2m$.
\item[$\diamond$]
If $m<j<2m$, then conjugation by $(T^{j-m},x)$ yields
$$((\Delta T^j,ex),(\Delta T^{i+j},exz))=
((\Delta T^{2m-j},ex^{-1}),(\Delta T^{i+2m-j},ex^{-1}z)),$$
so we may assume that $0\leq j\leq m$. \nl
It is important to note that these changes do not affect the order of $x$.
\item[$\diamond$]
If $j=0$, then lemma~\ref{lemmaeindord} gives the desired result.
\item[$\diamond$]
If $j=m$, then conjugation by $(T^{(m+i)/2},z^l)$ yields
\be{eqnarray*}
((\Delta T^m,ex),(\Delta T^{i+m},exz))&=&
((\Delta T^{-i},exz^{-2l}),(\Delta ,exz^{1-2l}))\\
&=&((\Delta ,e),(\Delta T^i,ez))
\e{eqnarray*}
The second identity follows from lemma~\ref{lemmaeindord}.
Note that $(T^{(m+i)/2},z^l)\in G$ if and only if
$(i+m)/2\equiv il\pmod{m}$.  But this condition is satisfied
because 
$$(i+m)/2=2^\mu(i_0+m_0)/2\equiv 2^\mu i_0l=il \pmod{2^\mu m_0}.$$
This finishes the proof in the case that $j=m$.
\item[$\diamond$]
If $0<j<m$, write $j=2^\nu j_0$ with $j_0$ odd.\nl 
We know that $(T^j,x)\in G$ and $x$ has order $2^k$ in $E$,
hence $m|j2^k$, i.e. $\mu\leq k+\nu$ and $m_0|j_0$.
Taking the fact that $j<m$ into account this implies $\nu<\mu$.
\be{itemize}
\item[$\cdot$]
Choose $r\in \N$ such that $r>k+\nu$ and $l_0|2^r-1$.
\item[$\cdot$]
Define $w\isdef z^{l^{\mu-\nu}}$.
\item[$\cdot$]
Choose an $\epsilon$ which satisfies the congruence
$$j\epsilon+il^{\mu-\nu}\equiv j+i_02^\nu\pmod{m}.$$
This is possible:\nl
mod $m_0$ it reads
\be{eqnarray*}
i_02^\nu2^{\mu-\nu} l^{\mu-\nu}&\equiv&i_02^\nu\\
i_02^\nu((2l)^{\mu-\nu}-1)&\equiv&0\\
i_02^\nu((l_0+1)^{\mu-\nu}-1)&\equiv&0,
\e{eqnarray*}
but since $m_0|i_0l_0$, this is automatically true;\nl
mod $2^\mu$ it reads
$$2^\nu j_0(\epsilon-1)\equiv 2^\nu i_0\pmod{2^\mu}$$
which is equivalent to
$$j_0(\epsilon-1)\equiv i_0\pmod{2^{\mu-\nu}},$$
but since $j_0$ is odd this is solvable.
\item[$\cdot$]
$ex(exex^\epsilon w)^{2^{r-\nu}}=
ex(x^\epsilon w)^{2^{r-\nu}}=
exz^{2^{r-\mu}}.$
\item[$\cdot$]
Define $\tilde{j}\isdef j-(2^r-1)2^\nu i_0$ and 
$\tilde{x}\isdef x^\epsilon$.
\e{itemize}
These definitions and facts support the following computation.
\be{eqnarray*}
((\Delta T^j,ex),(\Delta T^{i+j},exz))&=&\\
\mbox{ \ ($r-\mu$ times lemma~\ref{lem2powz})}&=&
((\Delta T^j,ex),(\Delta T^{i2^{r-\mu}+j},exz^{2^{r-\mu}}))\\
\mbox{ \ ($r-\nu$ times lemma~\ref{lem2powz})}&=&
((\Delta T^j,ex),(\Delta T^{i_02^\nu+j},e\tilde{x}w))\\
\mbox{ \ ($r$ times lemma~\ref{lem2powz})}&=&
((\Delta T^{j+2^\nu i_0-2^{r+\nu}i_0},e\tilde{x}),
(\Delta T^{i_02^\nu+j},e\tilde{x}w))\\
&=&((\Delta T^{\tilde{j}},e\tilde{x}),
(\Delta T^{i_02^{r+\nu}+\tilde{j}},e\tilde{x}w))\\
\mbox{ \ ($\mu-\nu$ times lemma~\ref{lem2powz})}&=&
((\Delta T^{\tilde{j}},e\tilde{x}),(\Delta T^{i2^r+\tilde{j}},e\tilde{x}w))\\
\mbox{ \ ($r$ times lemma~\ref{lem2powz})}&=&
((\Delta T^{\tilde{j}},e\tilde{x}),(\Delta T^{i+\tilde{j}},e\tilde{x}z))
\e{eqnarray*}
Observe that  $\tilde{j}$ is a multiple of $2^{\nu+1}m_0$.
Thus we replace the old $(T^j,x)$ by a new one.
Apply one of the preceding steps if $j\geq m$ or $j\leq0$.
\e{itemize}
Repeat this process until $\mu=\nu$, which implies $m|j$.
 
This completes the proof in this first case.
We did not need the fact that
$[(T^j,x)]=1\in \L([(T^i,z)])$!
This means that the primary Arf invariant is already good enough to detect
the Arf-elements in this case.
\e{nitel}
\be{nitel}{Case 2$\colon$ $i=0$}\nl
Our purpose is to show that 
$$((\Delta ,e),\Delta ,ez))=((\Delta T^j,ex),(\Delta T^j,exz))\in\Arf^h(G).$$
We apply induction on $j$, as follows.
\be{itemize}
\item[$\diamond$] If $j<0$ or $j>2m$,\nl
we conjugate by a suitable power of $(T^m,1)$, 
to attain $0\leq j\leq 2m$.
\item[$\diamond$] If $m<j\leq2m$,\nl
we conjugate by $(\Delta T^m,e)$, to achieve $0\leq j\leq m$.
\item[$\diamond$] If $j=0$,\nl
lemma~\ref{lemmaeindord} does the job.
\item[$\diamond$] If $j\not|\,m$,\nl
we define $d\isdef\gcd(j,m)$. 
Obviously $(T^d,x^{n_0})\in\ol{G_{(1,z)}}$ for some $n_0\in \Z$.
Now conjugating by $(T^d,x^{n_0})$ allows us to replace $j$ by $j-2d$.
Notice that $j-2d>0$.
\item[$\diamond$] If $j|m$,\nl 
there are two possibilities.
If there exists $(T^c,y)\in\ol{G_{(1,z)}}$ with $0<c<j$, then
we conjugate by $(T^c,y)$ to replace $j$ by $j-2c$. We have
$$-j+2\leq j-2c\leq j-2.$$  
Conjugating by $(\Delta ,e)$, 
if necessary, yields $$0\leq j-2c\leq j-2.$$
If there is not such a $c$, then the elements of 
$\ol{G_{(1,z)}}$ either have the
form $(T^{jv},\cdot\cdot)$ or
$(\Delta T^{jv},\cdot\cdot)$.
Since any element of $$\ker(\ol{G_{(1,z)}}\ra F(1,z))$$ 
is a product of
squares and $2$-power roots of $(1,z)$, so is $(T^j,x)$.
A little examination reveals that this can only happen when there exist
$2$-power roots $y_1$ and $y_2$ of $z$ such that
$(\Delta ,y_1),(\Delta T^j,y_2)\in G$\nl
Consider the pull-back diagram
$$\diagram{G&\lra&E\cr
\mapdown{ }&&\mapdown{\hat p}\cr
D&\lra&D_{m}\cr}$$
and define 
$$\be{array}{ll} 
F_1\isdef \hat{p}^{-1}(<\sigma\tau^\nu>)&\\
j_1\colon F_1\ra G \mbox{ \ by \ }
&j_1(f)\isdef\cases{(1,f)& if $\hat p(f)=1$\cr (\Delta,f)&otherwise\cr}\\
F_2\isdef \hat{p}^{-1}(<\sigma\tau^{\nu+j}>)&\\
j_2\colon F_2\ra G \mbox{ \ by  \ }
&j_2(f)\isdef\cases{(1,f)& if $\hat p(f)=1$\cr (\Delta T^j,f)&otherwise\cr}\\
F_0\isdef F_1\cap F_2=\ker(\hat p)&\vspace{1mm}\\
E_0\isdef\ker(w\colon E\ra \aut(C))&
\e{array}$$
Now $z\in F_0$, $e\in F_1$, $ex\in F_2$ and in the  diagram
$$\diagram{F_0&\subset&E_0\cr\bigcap&&\bigcap\cr F_1&\subset&E\cr}$$
$[F_1:F_0]=[E:E_0]=2$ and $[E:F_1]=[E_0:F_0]=m$.\nl
We know there exist $y_1\in F_1\setminus F_0$ such that $y_1z=zy_1$.
Then we have $ey_1\in F_0$ and $ey_1z(ey_1)^{-1}=\zm$, so
~\ref{defpseu} guarantees the existence of a pseudo-arf element
$\wp(f_1,f_1z)\in L^h(F_0)$.
Analogous, the existence $y_2\in F_2\setminus F_0$, satisfying
$y_2z=zy_2$, yields a pseudo-arf element $\wp(f_2,f_2z)\in L^h(F_0)$.
through the element $exy_2\in F_0$.\nl
But these pseudo-arf elements must coincide by theorem~\ref{thmwinjend}.
$$\diagram{
&\wp(f_2,f_2z)&&(ex,exz)\cr
\wp(f_1,f_1z)&L^h(F_0)&{\hbox to 30pt{\rightarrowfill}}&L^h(F_2)\cr
&\mapdown{}&
\be{picture}(20,20)
\put(-8,15){\vector(3,-2){32}}
\e{picture}
&\mapdown{j_{2*}}\cr
(e,ez)&L^h(F_1)&
{\buildrel j_{1*}\over {\hbox to 30pt{\rightarrowfill}}} &L^h(G)\cr
}$$
Therefore we may conclude that
$$j_{1*}((e,ez))=j_{2*}((ex,exz))\in L^h(G).$$ 
\e{itemize}
\e{nitel}
This completes the proof of theorem~\ref{thm2stuks}.
\e{proof}
\be{nitel}{proof of theorem~\ref{thmin2end}}
Suppose we have an element of $\Arf^h(G)$ which is killed by $\Upsilon$.
As before we may assume that it has the form $\sum\plane{a_i,a_iz}$.
We apply induction on the number of terms occuring in the expression 
for our element in $\Arf^h(G)$.
Recall that we are dealing with terms like 
$((\Delta ,e),(\Delta T^{i},ez))$ and
$((\Delta T^{j},ex),(\Delta T^{i+j},exz))$.
If there are less than three terms theorem~\ref{thm2stuks} does the job.
Thus assume that the number of terms exceeds two.
If a term $((1,\cdot\cdot),\cdots)$ appears,
lemma~\ref{lemmaeindord} enables us
to cancel two terms.
Otherwise there are two cases:
\be{itemize}  
\item[$\cdot$]
$i\neq0$\nl
We can cancel terms by the first case of theorem~\ref{thm2stuks}
without having any information on $(T^j,x)$. 
\item[$\cdot$]
$i=0$\nl
The following terms occur:
\be{eqnarray*}
&&((\Delta ,e),(\Delta T^{i},ez))\\
&&((\Delta T^{j_1},ex_1),(\Delta T^{i+j_1},ex_1z))\\
&&((\Delta T^{j_2},ex_2),(\Delta T^{i+j_2},ex_2z))
\e{eqnarray*}
Now define $j=\gcd(j_1,j_2)$, say $j=a_1j_1+a_2j_2$.
We conjugate the second and third term by a suitable power of
$$(T^j,x_1^{a_1}x_2^{a_2})$$ to obtain
$$((\Delta ,e\tilde{x}),(\Delta ,e\tilde{x}z))$$
or
$$((\Delta T^{j},e\tilde{x}),(\Delta T^{i+j},e\tilde{x}z)).$$
Applying lemma~\ref{lemmaeindord} once more, we can cancel terms.
\e{itemize}
We see that two terms cancel in all cases.
\e{nitel}

\newpage
\noindent{\Large\bf Samenvatting \vspace{7ex}}

In dit proefschrift worden nieuwe invarianten in de algebra\"{\i}sche 
$L$-theorie
bestudeerd. Een interessant bijprodukt van deze studie is de ontwikkeling
van gereduceerde machtsoperaties in cyclische homologie.

$L$-groepen van een ring met anti-structuur zijn in essentie
Grothendieck en Whitehead groepen van de categorie van kwadratische
modulen over die ring. Om meetkundige redenen is men vooral 
ge\"{\i}nteresseerd
in $L$-groepen van groepenringen.

Hoofdstuk I behelst een eigenzinnige poging de relevante $K$- en $L$-theorie
te etaleren en daarmee een geschikte setting voor de rest van het
proefschrift te cre\"{e}ren.
In hoofdstuk II worden $K$-theoretische invarianten geconstrueerd
door de anti-structuur op de gegeven ring, op een nogal exotische
manier, uit te breiden tot de ring van formele machtreeksen over die
ring. Deze invarianten generaliseren klassieke invarianten zoals de
Arf invariant.
Uitvoerige berekeningen aan de waardengroepen van deze
nieuwe invarianten leiden vervolgens tot de bestudering van homologische
invarianten in hoofdstuk III. Hiertoe worden gereduceerde machtsoperaties
in Hochschild, cyclische en quaternionische homologie gefabriceerd.
Hoofdstuk IV tenslotte, is gewijd aan de toepassing van het voorafgaande
op groepenringen. Het belangrijkste resultaat is hier dat de homologische
invariant `goed genoeg' is voor elke groep met twee einden
(Stelling 4.1).

\newpage
\noindent{\Large\bf Curriculum \vspace{7ex}Vitae}\newline
De schrijver van dit proefschrift werd geboren op 1 juli 1961 te
Heerlen. In 1979 behaalde hij het VWO diploma aan het Eykhagencollege
te Schaesberg. Daarna studeerde hij wiskunde aan de Katholieke
Universiteit te Nijmegen. In 1985 slaagde hij cum laude voor zijn
doctoraalexamen voorzien van een pedagogisch-didactische aantekening.
In oktober 1985 trad hij als wetenschappelijk assistent in dienst
van de Katholieke Universiteit. 
Tot maart 1990 verrichtte hij promotie-onderzoek onder leiding 
van Dr. F.J.-B.J. Clauwens. De meeste resultaten van dat onderzoek
zijn te vinden in dit proefschrift.  

\end{document}